\documentclass[12pt]{article}
\usepackage{amsfonts,amsmath,amsxtra}
\usepackage{latexsym}
\usepackage[matrix,arrow,curve]{xy}
\usepackage{amssymb}
\def\hybrid{\topmargin 0pt      \oddsidemargin 0pt
        \headheight 0pt \headsep 0pt
        \textwidth 16.5cm
        \textheight 23cm
        \marginparwidth 0.0in
        \parskip 5pt plus 1pt   \jot = 1.5ex}
\catcode`\@=11
\def\marginnote#1{}
\newcount\hour
\newcount\minute
\newtoks\amorpm
\hour=\time\divide\hour by60 \minute=\time{\multiply\hour by60
\global\advance\minute by-\hour}
\edef\standardtime{{\ifnum\hour<12 \global\amorpm={am}%
        \else\global\amorpm={pm}\advance\hour by-12 \fi
        \ifnum\hour=0 \hour=12 \fi
      \number\hour:\ifnum\minute<10 0\fi\number\minute\the\amorpm}}
\edef\militarytime{\number\hour:\ifnum\minute<10 0\fi\number\minute}

\def\draftlabel#1{{\@bsphack\if@filesw {\let\thepage\relax
   \xdef\@gtempa{\write\@auxout{\string
      \newlabel{#1}{{\@currentlabel}{\thepage}}}}}\@gtempa
   \if@nobreak \ifvmode\nobreak\fi\fi\fi\@esphack}
        \gdef\@eqnlabel{#1}}
\def\@eqnlabel{}
\def\@vacuum{}
\def\draftmarginnote#1{\marginpar{\raggedright\scriptsize\tt#1}}

\def\draft{\oddsidemargin -0.1truein
        \def\@oddfoot{\sl preliminary draft \hfil
        \rm\thepage\hfil\sl\today\quad\militarytime}
        \let\@evenfoot\@oddfoot \overfullrule 3pt
        \let\label=\draftlabel
        \let\marginnote=\draftmarginnote
\def\@eqnnum{{\rm (\theequation)}
\rlap{\kern\marginparsep\tt\@eqnlabel}%
\global\let\@eqnlabel\@vacuum}  }

\newfont{\Bbbb}{msbm7 scaled 1\@ptsize00}
\newcommand{\zs}{\raise-1pt\hbox{$\mbox{\Bbbb Z}$}}

\def\numberbysection{\@addtoreset{equation}{section}
        \def\theequation{\thesection.\arabic{equation}}}
\numberbysection

\renewcommand{\theequation}{\thesection.\arabic{equation}}

\def\titlepage{\@restonecolfalse\if@twocolumn\@restonecoltrue\onecolumn
     \else \newpage \fi \thispagestyle{empty}\c@page\z@
\def\thefootnote{\fnsymbol{footnote}} }
\def\endtitlepage{\if@restonecol\twocolumn \else  \fi
        \def\thefootnote{\arabic{footnote}}
        \setcounter{footnote}{0}}  
\relax
\hybrid
\makeatletter
\newdimen\normalarrayskip            
\newdimen\minarrayskip               
\normalarrayskip\baselineskip \minarrayskip\jot
\newif\ifold             \oldtrue            \def\new{\oldfalse}
\def\arraymode{\ifold\relax\else\displaystyle\fi}
\def\eqnumphantom{\phantom{(\theequation)}} 
\def\@arrayskip{\ifold\baselineskip\z@\lineskip\z@
     \else
     \baselineskip\minarrayskip\lineskip1\baselineskip\fi}
\def\@arrayclassz{\ifcase \@lastchclass \@acolampacol \or
\@ampacol \or \or \or \@addamp \or
   \@acolampacol \or \@firstampfalse \@acol \fi
\edef\@preamble{\@preamble
  \ifcase \@chnum
     \hfil$\relax\arraymode\@sharp$\hfil
     \or $\relax\arraymode\@sharp$\hfil
     \or \hfil$\relax\arraymode\@sharp$\fi}}
\def\@array[#1]#2{\setbox\@arstrutbox=\hbox{\vrule
     height\arraystretch \ht\strutbox
     depth\arraystretch \dp\strutbox
width\z@}\@mkpream{#2}\edef\@preamble{\halign \noexpand\@halignto
\bgroup \tabskip\z@ \@arstrut \@preamble \tabskip\z@ \cr}%
\let\@startpbox\@@startpbox \let\@endpbox\@@endpbox
  \if #1t\vtop \else \if#1b\vbox \else \vcenter \fi\fi
  \bgroup \let\par\relax
  \let\@sharp##\let\protect\relax
  \@arrayskip\@preamble}
\def\eqnarray{\stepcounter{equation}%
              \let\@currentlabel=\theequation
              \global\@eqnswtrue
              \global\@eqcnt\z@
              \tabskip\@centering              
              \let\\=\@eqncr
              $$%
            \halign to \displaywidth  \bgroup
             \eqnumphantom \@eqnsel
      \hskip\@centering                               
    $\displaystyle  \tabskip\z@ {##}$%
    &\global\@eqcnt\@ne \hskip 2\arraycolsep
         $ \displaystyle  \arraymode{##}$\hfil
    &\global\@eqcnt\tw@ \hskip 2\arraycolsep
         $\displaystyle\tabskip\z@{##}$\hfil
         \tabskip\@centering
    &{##}\tabskip\z@\cr}
\makeatother

\def\IC{\mathbb{C}}

\def\IH{\mathbb{H}}

\def\IR{\mathbb{R}}
\def\IZ{\mathbb{Z}}

\def\CC {\mathcal{C}}

\def\CS {\mathcal{S}}
\def\CT {\mathcal{T}}
\def\CU {\mathcal{U}}

\def\CZ {\mathcal{Z}}
\def\Fg{{\frak g}}

\def\fg{\mathfrak{g}}
\def\fh{\mathfrak{h}}

\def\fS{\mathfrak{S}}
\def\fM{\mathfrak{M}}

\def\a {{\alpha}}

\def\g {{\gamma}}
\def\s {{\sigma}}

\def\la{\lambda}
\def\La{\Lambda}
\def\ve{\varepsilon}

\def\vp{\varpi}
\def\e{\epsilon}

\def\fh{\mathfrak{h}}
\def\so{\mathfrak{so}}
\def\ssp{\mathfrak{sp}}
\def\gl{\mathfrak{gl}}

\def\wt {{\widetilde{\omega}}}

\def\Br{{\rm Br}}

\def\Id{{\rm Id}}

\def\Mat{{\rm Mat}}
\def\End{{\rm End}}

\def\Lie{{\rm Lie}}
\def\N{{\rm Nm}}
\def\Gal{{\rm Gal}}

\def\Pin{{\mathop{\rm Pin}}}

\def\Spin{{\mathop{\rm Spin}}}
\def\Sp{{\mathop{\rm Sp}}}
\def\SO{{\mathop{\rm SO}}}

\def\frak{\mathfrak}

\def\<{\langle}
\def\>{\rangle}
\def\s{\sigma}

\def\ad{{\rm ad}}
\def\wt{\widetilde}

\newtheorem{te}{Theorem}[section]
\newtheorem{de}{Definition}[section]
\newtheorem{prop}{Proposition}[section]           
\newtheorem{cor}{Corollary}[section]
\newtheorem{lem}{Lemma}[section]

\newcommand\bqa{\begin{eqnarray}}
\newcommand\eqa{\end{eqnarray}}
\def\be{\begin{eqnarray}\new\begin{array}{cc}}
\def\ee{\end{array}\end{eqnarray}}
\def\beq{\begin{equation}}
\def\eeq{\end{equation}}
\def\bse{\begin{subequations}}                
\def\ese{\end{subequations}}
\def\bp{\begin{pmatrix}}
\def\ep{\end{pmatrix}}

\newcommand{\proof}{\noindent {\it Proof} }
\newcommand{\ssl}{{\mathfrak{sl}}}
\newcommand{\Ad}{{\mathrm{ Ad}}}
\newcommand{\Aut}{{\mathrm{ Aut}}}
\newcommand{\Int}{{\mathrm{ Int}}}
\newcommand{\Hom}{{\mathrm{ Hom}}}
\newcommand{\Out}{{\mathrm{ Out}}}
\newcommand{\diag}{{\mathrm{ diag}}}
\newcounter{pac}[section]

\newcounter{pacc}[subsection]

0

\begin{document}

\draft

\title{\bf On normalizers of maximal tori \\ in classical Lie groups}
\author{A.A. Gerasimov, D.R. Lebedev and S.V. Oblezin}
\date{}
\maketitle

\renewcommand{\abstractname}{}

\begin{abstract}
\noindent {\bf Abstract}. The normalizer $N_G(H_G)$ of a maximal
torus  $H_G$ in  a  semisimple complex Lie group $G$  does   not in
general allow a presentation as a semidirect product of $H_G$ 
and  the corresponding Weyl group $W_G$. Meanwhile,  splitting holds
for classical groups corresponding to the root systems $A_\ell$,
$B_\ell$, $D_\ell$. For the remaining classical groups corresponding
to the root systems $C_\ell$  there still exists  an embedding of the
Tits extension of $W_G$  into normalizer $N_G(H_G)$.
We provide explicit unified construction of the  lifts
of the Weyl groups into normalizers of  maximal tori for classical
Lie groups corresponding to the root systems  $A_\ell$, $B_\ell$,
$D_\ell$ using  embeddings into  general linear Lie groups. For
symplectic series of classical Lie groups we provide an explanation
of impossibility of embedding of the Weyl group into the symplectic
group. The explicit formula for
 adjoint action of the lifts of the Weyl groups on $\Fg={\rm
 Lie}(G)$ are given. Finally some examples of the groups closely
associated with classical Lie groups are considered.

\end{abstract}


\vspace{5 mm}

\section{Introduction}

Normalizer $N_G(H_G)$ of a  maximal torus $H_G$ in a semisimple
complex Lie group $G$ allows presentation as an extension of the corresponding Weyl group
$W_G$  by $H_G$
 \be\label{NGH0}
  1\longrightarrow H_G\longrightarrow N_G(H_G)
  \stackrel{p}\longrightarrow W_G\longrightarrow 1.
 \ee
This extension does not split in general \cite{CWW}, \cite{AH}.
For the Weyl groups there is a well-known  presentation
via generators and relations. To obtain the
corresponding description of $N_G(H_G)$ one should pick a section of
the projection $p$. Universal solution to  this problem was given by
Demazure \cite{D} and Tits \cite{T2} in terms 
of the Tits extension $W_G^T$ of the Weyl group
$W_G$  (for recent discussions of the Tits groups see e.g. \cite{N},
\cite{DW}, \cite{AH}). One way to understand the  nature of the Tits
 extensions  is via  consideration of the  maximal
split real form $G(\IR)\subset G(\IC)$ of a complex semisimple Lie
group $G(\IC)$ \cite{GLO}.

In this note we consider the special class of classical
Lie groups corresponding to the root systems $A_\ell$, $B_\ell$,
$C_\ell$ and $D_\ell$ and provide an  explicit description of
the group $N_G(H_G)$ and lifts of the Weyl group which differs from
the one proposed in \cite{D}, \cite{T2}. We restrict ourselves to  the
Lie groups of classical type  due to the fact that 
by simple reasoning the exact sequence \eqref{NGH0} is split for all
classical groups except symplectic  ones.  

It is well-known that classical Lie groups may be defined as fixed point subgroups of
general Lie groups under appropriate involutions. 
In particular this allows to express generators of Weyl groups of the classical Lie
groups via generators of the Weyl groups of general linear Lie
groups. The explicit lift of the latter Weyl group into the general
linear Lie groups compatible with the action of the involution of
  the Lie algebra root date  define then the corresponding lift of
the Weyl groups for the classical Lie groups. This way  we  
obtain sections of $p$ in \eqref{NGH0} for all
classical Lie groups. In the case of orthogonal Lie groups
corresponding to series  $B_\ell$ and $D_\ell$ this provides
an embedding of the Weyl groups  into the corresponding Lie groups
while for symplectic groups we 
obtain an embedding of the Tits extension. In any case  this provides a
unified construction of the lifts of the Weyl groups of all
classical Lie groups.

The fact that for  $C_\ell$-series of classical Lie groups the
considered construction provides an embedding of
the   extension of the Weyl group
by $2$-torsion elements of a maximal torus, called the Tits
group but not  an embedding of 
the Weyl group itself looks rather surprising in this
context. We propose an explanation for this  phenomenon  based on
non-commutativity of quaternions. The main argument is
based on the known fact that while general linear and  orthogonal
Lie groups are naturally associated with matrix groups over
complex and real numbers  symplectic groups  allow
  description in terms of matrix groups over quaternions. This leads
to a modification of the standard notion of maximal torus by  taking into
account the non-commutative nature of quaternions. As a
consequence the notions of normalizer of maximal torus and the Weyl
group are also modified. The standard Weyl groups of symplectic Lie
groups appear as  subgroups of thus defined Weyl groups. The Tits
extension of the symplectic Weyl group then arises via construction of
universal cover of ${\rm Aut}(\IH)=SO_3$ by the three-dimensional
spinor group $\Spin_3$. Concretely this might be traced back to the fact
that that additional quaternionic unit $\jmath$ squares to minus
one.

It is worth mentioning that the standard algorithm of 
the  Gelfand-Zetlin construction of bases in finite-dimensional
irreducible representations  fails in the case of symplectic Lie groups.  
This  fact obviously reverberates with the absence of the  splitting of \eqref{NGH0} 
in symplectic case. We believe that this
is not accidental coincidence and the issue may be clarified using
quaternionic geometry.

Let us note in this respect that description of  classical Lie group in
terms of matrix algebras over division algebras might be
generalized to other, non-classical Lie groups (see e.g. \cite{B} and
references therein). We expect that the
(im)possibility of embedding of the Weyl groups  into the
corresponding Lie groups may be elucidated in all these cases
generalizing our considerations for symplectic Lie groups and
quaternionic matrix algebras.

In this note we also consider the problem of  construction of 
 a section of  $p$  in  \eqref{NGH0} for  unimodular linear
groups  $SL_{\ell+1}$, special orthogonal groups $SO_{\ell+1}$,
 pinor/spinor  groups  $\Pin_{\ell+1}$  and $\Spin_{\ell+1}$. Although these
groups are not classical Lie groups in  strict sense they are
related to classical Lie groups  either via central extensions or via taking unimodular
subgroups. In the case of central extension the resulting
section is apparently given by a central extension of the
corresponding Weyl group.   Moreover the case of unimodular
subgroups is also covered by central extension of the Weyl group. We
provide explicit description of maximal torus normalizers in terms
of generators and relations in all these cases.

The plan of the paper is as follows. In Section 2 we recall required
facts on  semisimple Lie algebras and groups including normalizers
of maximal tori (see also Appendix \ref{APP} for detailed discussion
of classical Lie groups). In Section 3  we recall  constructions of
classical Lie groups as fixed point subgroups of appropriate
involutions of  general Lie groups after E.Cartan. In Section 4 the
structure of  normalizers of maximal tori of general linear Lie
groups is considered in details including  explicit lifts of the
Weyl groups. In Section 5 we present construction of the lifts of
the Weyl groups of the classical Lie groups in terms of those for
general Lie groups (see Theorem \ref{MainTheorem}). In Section 6 we
further clarify the explicit formulas of the
previous Section 5 by establishing connections between maximal tori
normalizers, the  Weyl groups and their lifts for general linear
groups and its classical subgroups.  In Section 7 the underlying
reasons for the special properties of symplectic groups are
 considered. In Section 8 we
extend our analysis to the groups $SL_{\ell+1}$, $SO_{\ell+1}$,
$\Pin_{\ell+1}$  and $\Spin_{\ell+1}$. In all these cases we
construct explicit sections of $p$ in  \eqref{NGH0} realized  by central
extensions of the corresponding  Weyl groups. Finally in Section 9
we compute the adjoint action of the lifts of the Weyl groups on Lie
algebra $\Fg={\rm  Lie}(G)$ for all classical Lie groups. Various
  technical details of the construction are provided in   Appendix. 

{\it  Acknowledgments:} The research of the second author  was
supported by RSF grant 16-11-10075. The work of the third author was
partially supported by the EPSRC grant EP/L000865/1.

\section{Preliminaries on semisimple Lie algebras and groups}

  Let $\Fg$ be a complex semisimple Lie algebra and let
$\mathfrak{h}\subset\Fg$ be the Cartan subalgebra of $\fg$, so that
$\dim(\fh)={\rm rank}(\Fg)=\ell$. Let $\Pi=\{\a_i,\,i\in
I\}\subset\fh^*$ be the set of simple roots of $\fg$, indexed by the
set $I$ of vertexes of the Dynkin diagram $\Gamma$. Let
$\Pi^\vee=\{\a_i^{\vee},\,i\in I\}\subset\fh$ be the set of
corresponding coroots of $\fg$. Let $\Phi$ be the set of roots and
let $\Phi^{\vee}$ be the set of co-roots of $\fg$, then let
$(\Pi,\Phi;\,\Pi^{\vee},\Phi^{\vee})$ be the root system associated
with $\fg$, supplied with a non-degenerate $\IZ$-valued pairing
 \be
  \<\,,\,\>\,:\quad\Phi\times\Phi^\vee\,\longrightarrow\,\IZ\,.
 \ee

Introduce the weight lattice,
 \be
  \Lambda_W\,=\,\bigl\{\g\in\fh^*\,:\,
  \<\g,\,\a^{\vee}\>\in\IZ\,,\forall\a^{\vee}\in\Phi^{\vee}\bigr\}\subset\fh^*\,,
 \ee
and define the fundamental weights $\vp_i,\,i\in I$ by
 \be
 \<\vp_i,\,\a_j^{\vee}\>\,=\,\delta_{ij}\,.
 \ee
Fundamental weights provide a  basis of the weight lattice
$\La_W$ and we have
 \be
  \a_j\,=\,\sum_{i=1}^{\ell}a_{ij}\vp_i\,,\qquad \label{CartanMatrix}
  a_{ij}\,=\,\<\a_j,\,\a_i^{\vee}\>\,,
 \ee
where $A=\|a_{ij}\|$ is the Cartan matrix of $\Fg$. The root lattice
$\Lambda_R\subset\fh^*$ is generated by simple roots
$\a_i\in\Pi,\,i\in I$ and appears to be a sublattice of the weight
lattice $\Lambda_W$, so that $\La_R\,\subset\,\La_W$. The
 quotient group $\Lambda_W/\Lambda_R$ is a finite group of order
 \bqa\label{FundamentalGroup}
  |\Lambda_W/\Lambda_R|\,=\,\det(A)\,.
 \eqa

The co-weight and co-root lattices
$\Lambda_W^{\vee},\,\Lambda_R^{\vee}\subset\fh$ in the Euclidean
space $\bigl(\fh=\IC^{\ell};\,\<\cdot,\cdot\>\bigr)$ are generated
by the fundamental co-weights and simple co-roots, given by
 \be\label{coweights}
 \<\a_i,\,\vp_j^{\vee}\>=\delta_{ij}\,,\qquad
  \vp_j^{\vee}\,=\,\sum_{i=1}^{\ell}c_{ji}\a_i^{\vee}\,,\qquad
  \a_j^{\vee}\,=\,\sum_{k=1}^{\ell}a_{jk}\,\vp_k^{\vee}\,,
 \ee
where $C=\|c_{ij}\|=(A^T)^{-1}$ is the inverse transposed Cartan
matrix $A$.

The corresponding Weyl group $W(\Phi)$ is generated by simple roots
reflections,
 \be\label{WeylRoots}
  s_i(\a_j)=\a_j-\<\a_j,\a_i^\vee\>\a_i\,=\,\a_j-a_{ij}\a_i\,,
 \ee
and its action in $\fh^*$ preserves the root system
$\Phi\subset\fh^*$. Each symmetry of the Dynkin diagram
$\Gamma=\Gamma(\Phi)$ induces an automorphism of $\Phi$, and the
group $\Aut(\Phi)$ of all automorphisms of the root system
$\Phi\subset\fh^*$ contains $W(\Phi)$ as a normal subgroup.
Moreover, the following holds:
 \be\label{AutRoots}
  \Aut(\Phi)=W(\Phi)\rtimes\Out(\Phi)\,,
 \ee
with $\Out(\Phi)$ being a group of symmetries of $\Gamma(\Phi)$.

Let $\{h_i=h_{\a_i},\,e_i=e_{\a_i},\,f_i=e_{-\a_i},\,i\in I\}$ be
the standard set of generators of $\fg$:
 \be\label{ChSe}
  [h_i,h_j]\,=\,0\,,\\
  \bigl[h_i,\,e_j\bigr]\,=\,a_{ij}e_j\,,\qquad
  \bigl[h_i,\,f_j\bigr]\,=\,-a_{ij}f_j\,,\\
  \bigl[e_i,f_i\bigr]\,=\,h_i\,,\qquad
  \bigl[e_i,f_j\bigr]\,=\,0\,,\quad i\neq j\,;
 \ee
  \be\label{Serre2}
  {\rm ad}_{e_i}^{1-a_{ij}}(e_j)\,=\,{\rm
    ad}_{f_i}^{1-a_{ij}}(f_j)\,=\,0.
 \ee
 In the following a slightly modified presentation of \eqref{ChSe}
will be useful. Namely, the Lie algebra $\fg$ can be generated by
$\{\vp_i^{\vee}\,,\,e_i\,,\, f_i\,: i\in I\}$ subjected to the
following relations:
 \be\label{LieAlgebra}
  [\vp_i^{\vee},\,e_j]\,=\,e_j\delta_{ij}\,,\qquad
  [\vp_i^{\vee},\,f_j]\,=\,-f_j\delta_{ij}\,,\qquad
  [e_i,\,f_j]\,=\,\delta_{ij}\sum_{k=1}^{\ell}a_{jk}\,\vp_k^{\vee}\,.
 \ee

Now let $G_c$ be a connected compact semisimple Lie group of rank
$\ell$, such that $\fg_c=\Lie(G_c)$ is the tangent Lie algebra at
$1\in G_c$. Let $G=G_c\otimes\IC$ be the complexification of $G_c$
and let $\fg=\Lie(G)=\fg_c\otimes\IC$ be its Lie algebra. Then the
Lie algebra $\fg_c$ is generated by the following elements:
 \be\label{HJP}
  H_k\,=\,\imath h_k\,,\qquad
  J_k\,=\,f_k-e_k\,,\qquad
  P_k\,=\,\imath(e_k+f_k)\,,\qquad k\in I\,.
  \ee
In particular, the following relations hold:
 \be\label{HJPrel}
  \bigl[J_k,\,\imath\vp_j^{\vee}\bigr]\,=\,\delta_{kj}P_k\,,\quad
  \bigl[\imath\vp_j^{\vee},\,P_k\bigr]\,=\,\delta_{kj}J_k\,,\quad
  [P_k,J_k]\,=\,2\sum_{i=1}^{\ell}a_{ki}(\imath\vp_i^{\vee})\,.
 \ee

Let $H_G\subset G$ be its maximal torus, such that $\fh=\Lie(H_G)$,
and
 let $X^*(H_G)$ be the group of (rational) characters
$\chi:\,H_G\to S^1$. Then the differential $d\chi$ at $1\in G$ of a
character $\chi\in X^*(H_G)$ is a linear form on $\fh$; hence it
provides an embedding $X^*(H_G)\subset\fh^*$ as a discrete subgroup
(lattice), supplied with  a scalar product $\<\,,\,\>$.

The dual $X_*(H_G)\subset\fh$ of the lattice $X^*(H)$ is isomorphic
to $\fh_{\IZ}$. Moreover, using the above notations of the (co)root
and (co)weight lattices we have
 \be
  \La_R\,\subseteq\,X^*(H_G)\,\subseteq\,\La_W\,,\qquad
  \La_R^{\vee}\,\subseteq\,X_*(H_G)\simeq\fh_{\IZ}\,\subseteq\,\La_W^{\vee}\,.
 \ee
The adjoint action of maximal torus $H_G$ on the complex Lie algebra
$\fg$ provides the Cartan decomposition:
 \be
  \fg\,=\,\fh_{\IC}\oplus\bigoplus_{\a\in\Phi}\fg_{\a}\,,\qquad
  \fg_{\a}\,=\,\bigl\{X\in\fg\,:\quad\ad_h(X)\,=\,\a(h)X,\,\,\forall h\in\fh\bigr\}\,;\\
  \fg_{\a}\,=\,\IC e_{\a}\,,\qquad
  \fg_{-\a}\,=\,\IC f_{\a}\,,\qquad\a\in\Phi_+\,.
 \ee

The  root system $\Phi$ and the lattice $X^*(H_G)$ (or its dual
lattice $X_*(H_G)\simeq\fh_{\IZ}$) determine a unique (up to
isomorphism) connected semisimple Lie group $G$. In particular, $G$
is simply connected if and only if $X^*(H)\,\cong\,\Lambda_W$.

The center $\CZ(G)$ of a complex connected Lie group $G$ allows for
the following presentation:
 \be
  \CZ(G)\,\simeq\,\La_W^{\vee}/X_*(H_G)\,\simeq\,X^*(H_G)/\La_R\,,
 \ee
and  the group $\CZ(G)\rtimes\Out(\Phi)$ is isomorphic to the group
of symmetries of the extended Dynkin diagram $\wt{\Gamma}(\Phi)$.
More specifically, in case of simply-connected $G$ its center is
isomorphic to the quotient group \eqref{FundamentalGroup}:
 \be
  \CZ(G)\,\simeq\,\La_W/\La_R\,.
 \ee
In particular, for simply-connected complex Lie groups of types
$A_\ell$, $B_\ell$, $C_\ell$ and $D_\ell$ we have
 \be
  \CZ(A_{\ell})\simeq\IZ/(\ell+1)\IZ\,,\qquad
  \CZ(B_{\ell})\,\simeq\,\IZ/2\IZ\,,\qquad\CZ(C_{\ell})\,\simeq\,\IZ/2\IZ\,.\\
  \CZ(D_{\ell})\simeq\left\{\begin{array}{lc}
  \IZ/2\IZ\oplus\IZ/2\IZ\,, & \ell\in2\IZ\\
  \IZ/4\IZ\,, & \ell\in1+2\IZ \end{array}\right.\,.
\ee

Let $\Aut(G)$ be the group of all automorphisms of a connected
  complex semisimple Lie group $G$. Its connected component can be
identified with the group $\Int(G)$ of inner automorphisms, which is
isomorphic to the adjoint group:
 \be
  \Int(G)\,\simeq\,G/\CZ(G)\,.
 \ee
The quotient group over the connected component is the group of
outer automorphisms:
 \be\label{OutG}
  \Out(G)\,=\,\Aut(G)/\Int(G)\,.
 \ee
In case of simply-connected Lie group $G$, the following holds (see
e.g. \cite{L}, Chapter V Theorem 4.5.B):
 \be\label{OUT}
  \Out(G)\,\simeq\,\Out(\Phi)\,.
 \ee

\subsection{Normalizers of maximal tori}

Given a connected semisimple complex Lie group $G$ with a (fixed)
maximal torus $H_G\subset G$ of finite rank $\ell$, let
$N_G=N_{G}(H_G)$ be the normalizer of the maximal torus and let us
write
 \be\label{NGH}
  1\longrightarrow H_G\longrightarrow N_{G}
  \stackrel{p}\longrightarrow W_G\longrightarrow 1\,.
 \ee
The quotient group $W_G:=N_{G}/H_{G}$ is finite and is isomorphic to
the reflection group $W(\Phi_G)$ associated with the corresponding
root system $\Phi_G$. The group $W_G=W(\Phi_G)$ allows presentation
by simple root reflections $\{s_i,\,i\in I\}$ from \eqref{WeylRoots}
as generators subjected to the following relations:
 \be\label{BR0}
  s_i^2=1,
 \ee
 \be\label{BR}
 \underbrace{s_is_js_i\cdots }_{m_{ij}}\,
 =\,\underbrace{s_js_is_j\cdots }_{m_{ij}},\qquad i\neq j\in I,
 \ee
where $m_{ij}=2,3,4,6$ for $a_{ij}a_{ji}=0,1,2,3$, respectively.
Here $a_{ij}=\<\a_j,\,\a_i^{\vee}\>$ is the Cartan matrix entry
\eqref{CartanMatrix}. Equivalently these relations may be written in
the Coxeter form:
 \be\label{CoxeterRel}
  s_i^2=1, \qquad  (s_is_j)^{m_{ij}}=1\,, \qquad i\neq j\in I\,.
 \ee
The exact sequence \eqref{NGH} defines the canonical action of $W_G$
on $H_{G}$ so that the corresponding action on the Lie algebra
$\mathfrak{h}=\Lie(H_G)$ is provided by \eqref{WeylRoots}:
 \be\label{WeylGroupAction}
  s_i(h_j)\,=\,h_j-\<\a_i,\a_j^\vee\>h_i\,=\,h_j-a_{ji}h_i\,.
 \ee
This action preserves the scalar product $\<\,,\,\>$ in $\fh$, which
allows to identify Weyl group $W_G=W(\Phi_G)$ with a subgroup in the
orthogonal group $O\bigr(\fh,\,\<\,,\,\>\bigr)$.

Let us stress  that the situation is bit different in the case of non-connected
groups. Thus in the following we encounter an example of non-connected
group, the orthogonal group $O_{2\ell}$ having two connected
components. In this case the quotient $N_{O_{2\ell}}(H)/H$ is larger
then the Weyl group $W_{O_{2\ell}}$  defined by the generators and
relations \eqref{BR0}, \eqref{BR} and contains the outer
automorphism of the root system of the simple connected Lie groups  $SO_{2\ell}$.

The important fact is that  the exact sequence \eqref{NGH} does not
split in general i.e. $N_{G}$ is not necessarily isomorphic to the
semi-direct product $W_G\ltimes H_{G}$ (for various details see
\cite{D},\cite{T2}, \cite{CWW}, \cite{AH}).
 Thus in general we may only pick a suitable section of the
projection map $p$ in \eqref{NGH}. One such universal section for
arbitrary reductive Lie group we constructed by Tits \cite{T2}
saying that for general $G$ there exists a larger subgroup
$W^T_G\subset N_{G}$, containing $W_G$ as a quotient and fitting
into the following exact sequence:
 \be\label{Titsgroup}\xymatrix{
  1\ar[r] & H_{G}^{(2)}\ar[r] & W^T_G\ar[r] & W_G\ar[r] & 1}\,.
 \ee
 Here $H_{G}^{(2)}\simeq (\IZ/2\IZ)^{\ell}$
 is the 2-torsion subgroup  of the
maximal torus $H_{G}$ in the reductive complex Lie group $G$.
In case the group $G$ is semisimple the Tits group $W^T_G$ allows
for the following explicit presentation by the following generators
\cite{AH} (see the notations of \eqref{HJP})\footnote{In the
following expressions we write the group elements as exponential
   of the linear combinations of Lie algebra generators using
the canonical exponential map $\exp:\,\Lie(G)\to
G$. Note that map has a non-trivial kernel so that  for instance in case $G=GL_{\ell+1}(\IC)$ we
have
 $e^{2\pi J_i}\,=\,e^{2\pi P_i}\,=\,e^{2\pi\imath h_i}\,=\,1$ for each $i\in I$.}
 \be\label{Titsgen}
  \dot{s}_i\,=\,e^{\frac{\pi}{2}J_i}\,,\qquad i\in I\,,
 \ee
and relations (compare with \eqref{BR0}, \eqref{BR}):
  \be\label{Tits}
  (\dot{s}_i)^2\,=\,e^{\pi\imath h_i}\,,\qquad
  \underbrace{\dot{s}_i\dot{s}_j\cdots }_{m_{ij}}\,
  =\,\underbrace{\dot{s}_j\dot{s}_i\cdots }_{m_{ij}}\,\,,\quad i\neq
  j\,,\\
  \Ad_{\dot{s}_i}(h)\,=\,s_i(h)\,,\quad\forall h\in\fh\,.
  \ee
Using this presentation we readily observe that the Tits group
$W^T_G$
  is not only a subgroup of the
normalizer $N_{G}$, but is a subgroup of the corresponding compact
group $G_c\subset G$, and more precisely of the normalizer subgroup
$N_{G_c}(H_c)$ of the maximal torus $H_c$ in the compact group
$G_c$.

The true meaning of this construction
is rather elusive; in  \cite{GLO} we have proposed  some underlying
reasons for the existence of this Tits construction. Below we will
follow another direction and consider only the case of classical Lie
groups $G$. Recall that complex classical Lie group is a simple
reductive Lie group allowing embedding in
  the general Lie group as a fixed point subgroup of an
  involution. Equivalently the classical groups may be defined as
  stabilizer  subgroup of bilinear forms.  Essentially they are
  exhausted by the following families of the groups
 \be\label{ClassicalGroups0}
  GL_{\ell+1}(\IC), \qquad O_{2\ell+1}(\IC), \qquad  \Sp_{2\ell}(\IC),
  \qquad O_{2\ell}(\IC),
 \ee
corresponding to the series  $A_\ell$, $B_{\ell}$, $C_{\ell}$ and
$D_\ell$ of Dynkin diagrams. Let us stress that we distinguish the
classical Lie groups {\it per se} from  their various cousins like
$SL_{\ell+1}(\IC)$, $PSL_{\ell+1}(\IC)$, $\Spin_{\ell+1}(\IC)$, the
metaplectic group ${\rm Mp}_{2\ell}(\IC)$ et cet.

In the case of the classical Lie groups the results of
\cite{D}, \cite{T2}, \cite{CWW}, \cite{AH} might be formulated as the
following statement:

\begin{itemize}
\item \emph{In the case  of general linear group $GL_{\ell+1}(\IC)$ there is
a section to the exact sequences \eqref{NGH} and \eqref{Titsgroup},
so that $N_{GL_{\ell+1}(\IC)}$ (and the Tits subgroup
$W^T_{GL_{\ell+1}}$) contains the Weyl group  $W(A_{\ell})=W_{GL_{\ell+1}}$ associated
with the root system $A_{\ell}$;}

\item \emph{In the case of orthogonal groups $O_{\ell+1}(\IC)$ there
exits a section to the corresponding exact sequence
\eqref{Titsgroup}, so that $W^T_{O_{\ell+1}}$ contains a subgroup
isomorphic to $W_{O_{\ell+1}}$;}

\item \emph{In the case of the symplectic group $\Sp_{2\ell}(\IC)$
    no section of the corresponding exact sequence  \eqref{Titsgroup}
exists.}

\end{itemize}

\section{Classical Lie groups via involutive
automorphisms}\label{ClGrps}

The definition of classical Lie groups as simple Lie groups
isomorphic to fixed point subgroups with respect to certain
involutive automorphism in $GL_{\ell+1}(\IC)$ goes back
to E. Cartan and a relevant expositions of the subject can
be found in \cite{H} and \cite{L}.  Below we recall basics of this
construction.

Given the general linear group $GL_{\ell+1}(\IC)$ and a maximal torus
$H_{\ell+1}$ identified with the subgroup of diagonal elements, let
us describe its group of outer automorphisms $\Out(GL_{\ell+1})$.
Given Dynkin diagram $\Gamma(A_{\ell})$, let us choose natural
ordering of the set of its vertices $I=\{1,2,\ldots,\ell\}$, and
consider the outer automorphism of the root system induced by the
Dynkin diagram symmetry:
\be\label{OutDynkin}
\iota\,:\quad I\,\longrightarrow\,I\,,\qquad
 i\,\longmapsto\,\ell+1-i\,.
\ee
It defines the automorphism of the set $\Pi_{A_{\ell}}$ of simple roots, and
therefore can be extended to an (outer) automorphism of the root
system $\Phi(A_{\ell})$. Thus extended automorphism $\iota$ obviously preserves
the lattice $\fh_{\IZ}$ and can be lifted to the following
involutive automorphism of the Lie group $GL_{\ell+1}(\IC)$ (we keep
the same notation $\iota$):
 \be\label{CanonicalInv}
  \iota\,:\quad g\,\longmapsto\,(g^{\tau})^{-1}\,,\qquad g\in GL_{\ell+1}(\IC)\,,
 \ee
where $\tau$ is the reflection at the opposite diagonal:
 \be\label{TranspOpp}
  (g^\tau)_{ij}=g_{\ell+2-j,\ell+2-i}, \qquad g=\|g_{ij}\|\in
  GL_{\ell+1}(\IC).
  \ee
Note that the automorphism $\iota$ respects the maximal torus
$H_{\ell+1}\subset GL_{\ell+1}(\IC)$ given by invertible diagonal
matrices.  Given the involutive automorphism $\iota\in\Out(GL_{\ell+1})$ we have
 a family of involutive  automorphisms still respecting the maximal torus
 $H_{\ell+1}$. Such automorphisms  are obtained by combining $\iota$
with arbitrary inner automorphisms $\Ad_{\CS}\in\Int(GL_{\ell+1})$,
where $\CS\in N_{GL_{\ell+1}(\IC)}(H_{\ell+1})$ subjected
to the condition $\CS\CS^{\iota}\in\CZ(GL_{\ell+1})$.  Each such
$\theta=\Ad_{\CS}\circ\iota$
defines the corresponding fixed points subgroups
 \be
  G_{\theta} \,=\,(GL_{\ell+1}(\IC))^{\theta}\,.
 \ee
Classification of the  corresponding  subgroups is due to E. Cartan (see \cite{H}
Chapter IX, and \cite{L} Chapter VII for details). Crucial fact
following from  his theory is that  the requirement of  $G$ being a
simple reductive Lie group imposes strong restrictions on  
possible choice of the inner automorphism $\Ad_{\CS}$. The
corresponding list of  simple reductive groups is exhausted by
the following cases
\be\label{ClassicalGroups}
  O_{2\ell+1}\,=\,\bigl\{g\in GL_{2\ell+1}\,:\quad\theta_B(g)=g\,\bigr\}\,,\\
  \Sp_{2\ell}\,=\,\bigl\{g\in
  GL_{2\ell}\,:\quad\theta_C(g)=g\,\bigr\}\,,\\
  O_{2\ell}\,=\,\bigl\{g\in GL_{2\ell}\,:\quad\theta_D(g)=g\,\}\,,
 \ee
where $\theta_G$ are the involutive automorphisms of general linear
group
 \be\label{CanonicalInv1}
  \theta_G\,:\quad g\,\longmapsto\,\CS_G\,g^{\iota}\,\CS_G^{-1}\,,
 \ee
with the diagonal matrices $S_G$ given by
 \be\label{CanonicalInv2}
  \CS_B\,=\,\diag(1,-1,\,\ldots,\,1)\,\in\,GL_{2\ell+1}\,,\qquad
  \det(\CS_B)\,=\,(-1)^{\ell}\,;\\
  \CS_C\,=\,\diag(1,-1,\,\ldots,\,1,-1)\,\in\,GL_{2\ell}\,,
  \qquad\det(\CS_C)\,=\,(-1)^{\ell}\,;\\
  \CS_D\,=\,\CS_C\cdot\eta\,,\qquad\eta\,=\,\diag(\Id_{\ell},\,-\Id_{\ell})\,,
 \qquad   \det(\CS_D)\,=\,1\,.
 \ee
There is also trivial case of
$\theta_G={\rm id}$ corresponds to the general linear Lie group
$GL_{\ell+1}(\IC)$ itself.
In the following  we will use the term ``\emph{classical group}'' in the
narrow sense by excluding the trivial case of general linear groups.

Let us note that the construction of classical Lie groups via
involutive automorphims described above   differs from  more traditional one where
instead of the reflection at the opposite diagonal \eqref{TranspOpp}
the standard transposition is used. For instance  elements of the
orthogonal  subgroup $O_{\ell+1}\subset GL_{\ell+1}(\IC)$  are
defined by   the condition $g=(g^t)^{-1}$ where
\be\label{TRST}
(g^t)_{ij}=g_{j,i}, \qquad g=\|g_{ij}\|\in   GL_{\ell+1}(\IC).
\ee  
Although two constructions are equivalent (one may be obtained
from another by adjoint action of a certain group element of
the general linear group) the construction \eqref{ClassicalGroups}
is more convenient when dealing with symmetries of the Dynkin diagrams
and abstract root data \cite{DS}. Thus in the following we use
\eqref{ClassicalGroups} as our basic definition.

\section{The normalizer of maximal torus in
$GL_{\ell+1}$}\label{GLNORM}

In the case of  general linear Lie groups $GL_{\ell+1}$ the exact
sequence \eqref{NGH} is split by simple reasons. Indeed each element
$g\in N_{\ell+1}$ induces an (inner) automorphism $\Ad_g$ of the Lie
algebra $\gl_{\ell+1}(\IC)$, preserving the Cartan subalgebra $\fh$
as well as the scalar product $\<\,,\,\>$ in it. This defines a
group homomorphism from $N_{\ell+1}$ into the orthogonal group
$O\bigl(\fh_{\IC},\,\<\,,\,\>\bigr)$, so that its image is the Weyl
group:
 \be\label{splitProjGL}
  N_{\ell+1}/\CZ(GL_{\ell+1})\,\longrightarrow\,W_{GL_{\ell+1}}\subset
  O_{\ell+1}(\IC)\,\subset\,GL_{\ell+1}(\IC)\,.
 \ee
This yields a canonical embedding of the Weyl group into the
orthogonal group, which gives rise to homomorphism
$W_{GL_{\ell+1}}\to GL_{\ell+1}(\IC)$.

This  argument  implies the following explicit presentation of the
lifts of the simple root generators $s_i\in
W(A_{\ell})=W_{GL_{\ell+1}},\,1\leq i\leq\ell$ as elementary
permutation matrices
 \be\label{permat}
  s_i\, \longmapsto\, S_i\,=\,\left(\begin{smallmatrix}
  \Id_{i-1}&&&\\&0&1&\\&1&0&\\
  &&&\Id_{\ell-i}
  \end{smallmatrix}\right)\,\in\,N_{\ell+1}\,,\quad
  \qquad S_i^2\,=\,1\,,
  \ee
  satisfying the standard Coxeter relations
 \be\label{WeylPerm}
  (S_iS_j)^2\,=\,1\,,\quad|i-j|>1\,,\qquad(S_iS_j)^3\,=\,1\,,\quad|i-j|=1\,,
 \ee
of the generators of permutation group $\mathfrak{S}_{\ell+1}$.

 The lifts $S_i$ of the generators $s_i$ of the Weyl group $W_{GL_{\ell+1}}$
  should  be compared with the generators of the corresponding Tits group
$W^T_{GL_{\ell+1}}$ defined as follows.
 For the reductive group $GL_{\ell+1}(\IC)$ the Tits group
$ W^T_{GL_{\ell+1}}\subset N_{\ell+1}$ is given by the extension
 \be\label{TitsgroupGL} \xymatrix{
   1\ar[r] & H_{\ell+1}^{(2)}(\IC)\ar[r] &
   W^T_{GL_{\ell+1}}\ar[r] & W_{GL_{\ell+1}}\ar[r] &
  1},
\ee
where  $H^{(2)}_{\ell+1}$ is the 2-torsion subgroup of the
maximal torus $H_{\ell+1}$. Let $\{e_{kk},\,1\leq k \leq
(\ell+1)\}$ be  the bases of the coweight lattice of
$GL_{\ell+1}(\IC)$ corresponding the diagonal matrices with only  one
non-zero element being a unit on the diagonal (it might be expressed through the
fundamental coweights as follows $e_{kk}=p^{\vee}_k-p_{k-1}^{\vee}$,
see \eqref{GLcoweights}). Then $H^{(2)}_{\ell+1}$ is generated
by the elements
 \be\label{2torsGL}
  T_k\,:=\,e^{\pi\imath e_{kk}}\,,\quad1\leq k\leq (\ell+1)\,.
 \ee
 The group  $W^T_{GL_{\ell+1}}$ is generated by
 \eqref{Titsgen} for  the roots system $A_{\ell}$
 \be
 \dot{s}_i\,=\,e^{\frac{\pi}{2}J_i}\,
  =\,\left(\begin{smallmatrix}
  \Id_{i-1}&&&\\&0&-1&\\&1&0&\\
  &&&\Id_{\ell-i}
  \end{smallmatrix}\right)\in N_{\ell+1}\,,\quad
  (\dot{s}_i)^2\,=\,T_iT_{i+1}^{-1}\,,
 \ee
 together with an additional central element
 given by product of all $T_k$, $k=1,\ldots ,\ell+1$.
Here we write down  matrix form of generator $\dot{s}_i$ using
standard faithful representation. Comparing the matrix forms
we arrive at the following relations:
 \be\label{permat0}
  S_i\,=\,T_i\dot{s}_i\,=\,\dot{s}_iT_{i+1}\,\,,\qquad 1\leq  i\leq\ell\,.
  \ee
This leads to another
  presentation of the Tits group.

\begin{lem}\label{LM1} The elements
 \be\label{2torsGLN}
  S_i\,,\quad1\leq i\leq \ell\qquad\text{and}\qquad
  T_k\,,\quad1\leq k\leq (\ell+1)\,,
  \ee
satisfying  the relations
 \be\label{WeylExt1}
  S_i^2\,=\,T_k^2\,=\,1\,,\qquad T_kT_n\,=\,T_nT_k\,,\\
  (S_iS_j)^2\,=\,1\,,\quad|i-j|>1\,,\qquad(S_iS_j)^3\,=\,1\,,\quad|i-j|=1\,,\\
  S_iT_k\,=\,T_{s_i(k)}S_i\,,
  \ee
provide a presentation of the Tits group $W^T_{GL_{\ell+1}}$.

\end{lem}
\proof: The only non-obvious relation  in
the last line of \eqref{WeylExt1} follows from explicit matrix
computation:
 \be
  T_iS_iT_i^{-1}\,
  =\,T_{i+1}S_iT_{i+1}^{-1}\,
  =\,S_iT_iT_{i+1}\,=\,T_iT_{i+1}S_i\,,\qquad1\leq i\leq\ell\,.
 \ee
$\Box$

In Section \ref{ClGrps} we have introduced involutive
automorphisms $\theta_G\in\Aut(GL_{\ell+1}(\IC))$
\eqref{ClassicalGroups}, \eqref{CanonicalInv1},
\eqref{CanonicalInv2} defining the classical Lie groups. The
involutions $\theta_G$ do not act on the lift of $W_{GL_{\ell+1}}$
defined by the generators $S_i$ but do act on its extension given by
the Tits group. Explicitly we have
 \be
  \theta_B\,:\qquad\dot{s}_i\,\longmapsto\,\dot{s}_{2\ell+1-i}\,,\qquad
  T_i\,\longmapsto\,T_{2\ell+2-i}\,,\\
  \theta_C\,:\qquad\dot{s}_i\,\longmapsto\,\dot{s}_{2\ell-i}\,,\qquad
  T_i\,\longmapsto\,T_{2\ell+1-i}\,,\\ \label{CINV}
  \theta_D\,:\qquad\dot{s}_i\,\longmapsto\,\dot{s}_{2\ell-i}\,,\qquad
  T_i\,\longmapsto\,T_{2\ell+1-i}\,.
 \ee
Therefore, thus defined automorphism
$\theta_G\in\Aut(W^T_{GL_{\ell+1}})$ preserves the normal
subgroup $H^{(2)}_{\ell+1}=\<T_k,\,1\leq k\leq\ell+1\>$.

Let us describe the action of the involutions $\theta_G$ in terms of
the another set of generators of $W^T_{GL_{\ell+1}}$; introduce new
elements given by
 \be\label{permat1}
  \overline{S}_i\,
  =\,T_iS_iT_i^{-1}\,
  =\,T_{i+1}S_iT_{i+1}^{-1}\,
  =\,\left(\begin{smallmatrix}
  \Id_{i-1}&&&\\&0&-1&\\&-1&0&\\
  &&&\Id_{\ell-i}
  \end{smallmatrix}\right)\,,\quad
  \overline{S}_i\,=\,(\overline{S}_i)^{-1}\,\,.
  \ee

\begin{lem} The elements $S_i,\,\overline{S}_i\in
  W^T_{GL_{\ell+1}}$ satisfy the following relations:
 \be\label{WeylExt2}
 \begin{array}{c}
  S_i^2\,=\,\overline{S}_i^2\,=\,1\,,\qquad
  S_i\overline{S}_i=\overline{S}_iS_i\,=\,T_iT_{i+1}\,,\qquad1\leq i\leq\ell\,;\\
   (S_iS_j)^2\,=\,(S_i\overline{S}_j)^2\,=\,(\overline{S}_i\overline{S}_j)^2\,=\,1\,,
   \qquad
  |i-j|>1\,,\\
  (S_iS_j)^3\,=\,(\overline{S}_i\overline{S}_j)^3\,=\,1\,,\quad
  S_iS_jS_i\,=\,\overline{S}_j\overline{S}_i\overline{S}_j\,\,,\qquad|i-j|=1\,.
 \end{array}
 \ee
\end{lem}
\proof: The relations follow from the similar relations for the
symmetric group \eqref{BR}, \eqref{BR0} and \eqref{CoxeterRel}. The
latter (braid) relation can be checked using
 \be
  \overline{S}_i\,=\dot{s}_iT_i\,=\,T_{i+1}\dot{s}_i\,,
 \ee
and the braid relation \eqref{Tits} in the Tits group. $\Box$

The action of the involutions $\theta_G$ \eqref{CINV}  may be
written as follows
 \be\label{permatAction}
 \theta_B(S_i)\,
  =\,\overline{S}_{2\ell+1-i}\,,\qquad
 \theta_C(S_i)\,
  =\,\overline{S}_{2\ell-i}\,,\quad1\leq i\leq\ell\,,\\
  \theta_D(S_i)\,=\,\overline{S}_{2\ell-i}\,,\quad1\leq i<\ell\,,\qquad
  \theta_D(S_{\ell})\,=\,S_{\ell}\,.
  \ee
Note that the two sets $\{S_i\}$ and $\{\overline{S}_i\}$ of
generators define two different embeddings of the Weyl group
$W_{GL_{\ell+1}}$ into the Tits group $W^T_{GL_{\ell+1}}$ (and thus
in the normalizer $N_{\ell+1}$) which are interchanged by
involutions $\theta_G$. Let us stress that for example it is possible to find
  such embedding $W_{GL_{\ell+1}}\subset GL_{\ell+1}(\IC)$ that the image
  of $W_{GL_{\ell+1}}$ is invariant with respect to
  $\theta_{O_{\ell+1}}$. Indeed this fact is obvious in the case of the
  realization of the orthogonal subgroup $O_{\ell+1}\subset
  GL_{\ell+1}(\IC)$ based on \eqref{TRST} because  for $S_i$ defined
  by \eqref{permat} we have $S_i\in O_{\ell+1}$. The same fact for the
  realization   \eqref{ClassicalGroups} follows from the
  equivalence of the two realizations. However the resulting matrix
  expressions for the lifts of $s_i$ are more involved and require the
   use of  complex numbers while in our realization the matrix entries are
  elements of the set $\{-1,0,+1\}$. Thus in the following we will use
  the pair $\{S_i\}$, $\{\overline{S}_i\}$  of lifts related by    
\eqref{permatAction}.

\section{Classical groups and their Weyl groups}\label{MAINRES}

For each complex classical Lie group $G=GL_{\ell+1}(\IC)^{\theta_G}$ the
corresponding involution $\theta_G\in\Aut(GL_{\ell+1})$ respecting
the embedding $H_{\ell+1}\subset GL_{\ell+1}(\IC)$ of the diagonal maximal
torus $H_{\ell+1}$ naturally acts on the Lie algebra  $\gl_{\ell+1}(\IC)$,
the corresponding root system $\Phi_{GL_{\ell+1}}$, normalizer
$N_{\ell+1}(H_{\ell+1})$  of $H_{\ell+1}$  and the Weyl group
$W_{GL_{\ell+1}}=N_{\ell+1}(H_{\ell+1})/H_{\ell+1}$. Both the
maximal torus $H_G$ of $G$, normalizer $N_G(H_G)$ and its quotient
$W_G=N_G(H_G)/H_G$ can be expressed in terms of the corresponding
objects of the underlying linear group $GL_{\ell+1}$. Let us give
explicit expressions for the generators of the group of automorphysms
${\rm Aut}(\Phi_G)$ of the root systems and the corresponding Weyl
groups in terms of the generators of the Weyl group of the
corresponding general linear Lie
group.

\begin{prop} The following explicit description of the generators of
  the group $\Aut(\Phi_G)$  of automorphisms of the classical Lie
  algebra root systems $\Phi_G$   and the corresponding  Weyl groups holds:

  \begin{itemize}
\item For $G=O_{2\ell+1}$ one has
$\Aut(\Phi(B_{\ell}))=W(\Phi(B_{\ell}))$, and the simple root
generators $s^{B_{\ell}}_i\in W(\Phi(B_{\ell}))$ can be expressed in
terms of $s_i^{A_{2\ell}}\in W(\Phi(A_{2\ell}))$ as follows:
 \be\label{BWeyl0}
  s_1^{B_{\ell}}\,=\,s_{\ell}^{A_{2\ell}}s_{\ell+1}^{A_{2\ell}}s_{\ell}^{A_{2\ell}}\,
  =\,s_{\ell+1}^{A_{2\ell}}s_{\ell}^{A_{2\ell}}s_{\ell+1}^{A_{2\ell}}\,,\\
  s_k^{B_{\ell}}\,=\,s_{\ell+1-k}^{A_{2\ell}}s_{\ell+k}^{A_{2\ell}}\,
  =\,s_{\ell+k}^{A_{2\ell}}s_{\ell+1-k}^{A_{2\ell}}\,,\quad1<k\leq\ell\,.
 \ee

\item For $G=\Sp_{2\ell}$ one has
$\Aut(\Phi(C_{\ell}))=W(\Phi(C_{\ell}))$, and the simple root
generators $s^{C_{\ell}}_i\in W(\Phi(C_{\ell}))$ can be expressed in
terms of $s_i^{A_{2\ell-1}}\in W(\Phi(A_{2\ell-1}))$ as follows:
 \be\label{CWeyl0}
  s_1^{C_{\ell}}\,=\,s_{\ell}^{A_{2\ell-1}}\,,\qquad
  s_k^{C_{\ell}}\,=\,s_{\ell+1-k}^{A_{2\ell-1}}s_{\ell-1+k}^{A_{2\ell-1}}\,
  =\,s_{\ell-1+k}^{A_{2\ell-1}}s_{\ell+1-k}^{A_{2\ell-1}}\,,\quad1<k\leq\ell\,.
 \ee

\item For $G=O_{2\ell}$ one has
  $\Aut(\Phi(D_{\ell}))=W(\Phi(D_{\ell}))\rtimes\Out(\Phi(D_{\ell}))$.
  The simple root generators $s^{D_{\ell}}_i\in W(\Phi(D_{\ell}))$
together with the generator $R\in\Out(\Phi(D_{\ell}))=\IZ/2\IZ$ can
be expressed in terms of $s_i^{A_{2\ell-1}}\in W(\Phi(A_{2\ell-1}))$
as follows:
 \be\label{DWeyl0}
 s_1^{D_{\ell}}\,=\,s_{\ell}^{A_{2\ell-1}}
 s_{\ell-1}^{A_{2\ell-1}}s_{\ell+1}^{A_{2\ell-1}}s_{\ell}^{A_{2\ell-1}}\,
  =\,s_{\ell}^{A_{2\ell-1}}s_{\ell+1}^{A_{2\ell-1}}
  s_{\ell-1}^{A_{2\ell-1}}s_{\ell}^{A_{2\ell-1}}\,,\\
  s_k^{D_{\ell}}\,=\,s_{\ell+1-k}^{A_{2\ell-1}}s_{\ell-1+k}^{A_{2\ell-1}}\,
  =\,s_{\ell-1+k}^{A_{2\ell-1}}s_{\ell+1-k}^{A_{2\ell-1}}\,,\quad1<k\leq\ell\,;\\
  R\,=\,s_{\ell}^{A_{2\ell-1}}\,.
 \ee
\end{itemize}
\end{prop}
\proof: We verify the assertion via case by case study in Lemmas
\ref{BAutroots}, \ref{CAutroots} and \ref{DAutroots} in the
following subsections. $\Box$

The Weyl group $W(A_{\ell})$ allows embedding into the Lie group
$GL_{\ell+1}$ as was discussed in  Section \ref{GLNORM}.
Therefore it is now  natural to ask for a possibility to lift the
presentations \eqref{BWeyl0}, \eqref{CWeyl0}, \eqref{DWeyl0} of Weyl
groups $W_G$ into the corresponding classical groups
$G=GL_{\ell+1}(\IC)^{\theta_G}$ via the  pair of lifts of the Weyl group
$W_{GL_{\ell+1}}$ into the general linear group $GL_{\ell+1}(\IC)$
considered previously.  For the Tits
extension $W^T_{GL_{\ell+1}}$ the following  result is easily checked.  

\begin{prop} The explicit presentation for  the Tits group $W^T_G$
of the classical Lie group $G=(GL_{\ell+1})^{\theta_G}$
may be provided by replacement of all $s_i$'s by $\dot{s}_i$'s in \eqref{BWeyl0},
\eqref{CWeyl0} and \eqref{DWeyl0}.

\end{prop}

\proof: We give the case by case verification in Lemmas
\ref{BWeylLift}, \ref{CWeylLift} and \ref{DWeylLift} by
straightforward computation using standard faithful representations
\eqref{ClassicalGroups}. $\Box$

For the Weyl groups we can not expect such simple answer. Indeed  
we know that although for the groups $O_{2\ell}$ and $O_{2\ell+1}$
the corresponding Weyl groups allow embedding  in the corresponding classical group
this is not so  for $\Sp_{2\ell}$.  So for the case of the Weyl
groups we have the following result.

\begin{te}\label{MainTheorem}
Given a classical group $G$ with its maximal torus $H_G\subset G$,
let $S_i,\,\overline{S}_i,\,1\leq i\leq\ell$ be the generators
\eqref{permat}, \eqref{permat1} of the Tits group
$W^T_{GL_{\ell+1}}$ and let $T_k,\,1\leq k\leq\ell+1$ be the
elements \eqref{2torsGL}. Then the following holds:
\begin{itemize}
\item In case $G=O_{2\ell+1}$ the generators $S^B_i$ defined by
 \be\label{2torsB0}
  S^B_1\,
  =\,S_{\ell+1}S_{\ell}S_{\ell+1}\,
  =\,S_{\ell}S_{\ell+1}S_{\ell}\,,\\
  S^B_k\,=\,S_{\ell+1-k}\overline{S}_{\ell+k}\,,\qquad1<k\leq\ell\,,
  \ee
are $\theta_B$-invariant and generate a finite group isomorphic to
$W(B_{\ell})=W_{O_{2\ell+1}}$.

\item In case $G=O_{2\ell}$ the generators $S^D_i$ defined by
 \be\label{2torsD0}
  S^D_1\,=\,S_{\ell}S_{\ell-1}\overline{S}_{\ell+1}S_{\ell}\,,\\
  S^D_k\,=\,S_{\ell+1-k}\overline{S}_{\ell-1+k}\,,\qquad1<k\leq\ell\,,
 \ee
are $\theta_D$-invariant and generate a finite group isomorphic to
$W(D_{\ell})=W_{O_{2\ell}}$.

\item In case $G=\Sp_{2\ell}$ the Tits group $W^T_{Sp_{2\ell}}$ is
  generated by $\theta_C$-invariant  generators
 \be
 S^C_1\,=\,T_{\ell}S_{\ell}\,,\\
  S^C_k\,=\,S_{\ell+1-k}\overline{S}_{\ell-1+k}\,,\qquad1<k\leq\ell\,,
  \ee
while the group of $\theta_C$-invariant combinations of
$S_i,\,\overline{S}_i,\,1\leq i\leq\ell$  generated by
 \be\label{2torsC0}
  \tilde{S}^C_1\,=\,S_{\ell}\overline{S}_{\ell}\,=\,T_{\ell}T_{\ell+1}^{-1}\,,\\
  S^C_k\,=\,S_{\ell+1-k}\overline{S}_{\ell-1+k}\,,\qquad1<k\leq\ell\,,
 \ee
 is isomorphic to a proper subgroup of  the Weyl group
$W(C_{\ell})$.
\end{itemize}
\end{te}
\proof: We propose the detailed proof the assertion via case by case
study in Propositions \ref{BWeylLift}, \ref{CWeylLift} and
\ref{DWeylLift} below. $\Box$

\section{Description of  $\theta_G$-invariants}

The explicit expressions for the generators the Weyl groups/Tits
groups of the classical groups presented above implies some general
relations between maximal torus normalizers, Weyl groups and Tits
groups for general linear groups $GL_{\ell+1}(\IC)$
and their classical subgroups $G=(GL_{\ell+1}(\IC))^{\theta_G}$.
In this Section we prove a set of such relations. First let us note
the following property of the involutions $\theta_G$ defining
classical Lie groups $G$.

\begin{lem}\label{PROPC} The centralizer of the fixed point subset
$H_{\ell+1}^{\theta_G}$ in $GL_{\ell+1}(\IC)$ is $H_{\ell+1}$.
\end{lem}

\proof \, We shall solve the set of equation for
  $g\in GL_{\ell+1}(\IC)$
  \be\label{CENTH}
  gh\,=\,hg, \qquad  h\in H_{\ell+1}^{\theta_G}.
  \ee
  Due to obvious isomorphism $H_{\ell+1}^{\theta_G}=H_{\ell+1}^{\iota}$
  we may use the basic autmomorphims $\iota$ \eqref{TranspOpp}
instead of $\theta_G$ in \eqref{CENTH}.  We fix $H_{\ell+1}\subset
GL_{\ell+1}$ be the diagonal subgroup. Explicitly elements  $g\in
C_{GL_{\ell+1}}(H_{\ell+1}^{\iota})$ of the centralizer
$C_{GL_{\ell+1}}(H_{\ell+1}^{\iota})$ of invariant subtorus
$H_{\ell+1}^{\iota}\subset H_{\ell+1}$ satisfy the conditions
\be
x_ig_{ij}\,=\,g_{ij}x_j\,,
\ee
where
\be
g=\|g_{ij}\|, \qquad
  h\,=\,\diag(x_1,\,x_2,\ldots,x_{\ell+1})\,,
\ee
with the additional condition $x_{\ell+2-i}\,=\,x_i^{-1}$ reflecting
$\iota$-invariance of $h$. Since $g$ does not depend on $x_i$ and
there are no universal $x$-independent relations  between entries of $h$
(although they subjected quadratic relations $x_ix_{\ell+2-i}=1$) we
infer $g_{ij}=0,\,i\neq j$ i.e. $g\in H_{\ell+1}$.\,\, $\Box$

As we already notice the definition of the Weyl group as a quotient
of the normalizer $N_G(H_G)$ by $H_G$ is not equivalent to the
definition of the Weyl group via  generators and relations
\eqref{BR0}, \eqref{BR} for non-connected classical Lie groups. Precisely the
group  $N_G(H_G)/H_G$ for all classical Lie groups obtained as fixed
point subgroups of the general linear group (we exclude the case of
trivial involution) may be identified with the group of automorphims
of the corresponding root system
\be
N_G(H_G)/H_G\simeq {\rm
  Aut}(\Phi_G)\,.
\ee
 This group fits the following exact sequence
 \be
 1\longrightarrow W_G\longrightarrow {\rm
   Aut}(\Phi_G)\longrightarrow {\rm Out}(\Phi_G)\longrightarrow 1
 \ee
 and we have
 \be
 {\rm Aut}(B_{\ell})\simeq W(B_{\ell}), \qquad {\rm
Aut}(C_{\ell})\simeq W(C_{\ell}), \qquad {\rm Aut}(D_{\ell})\simeq
W(D_{\ell})\times \IZ_2. \ee Let us stress that notation
$\Out(\Phi_G)$ is natural for connected simple groups where indeed
this group is the factor of all automorphims over inner
automorphisms. In the case of non-connected group (see detailed
discussion of the group $O_{2\ell}$ in Section \ref{APPD}) the whole
group ${\rm Aut}(\Phi_G)$ is realized by inner automorphims.

\begin{prop}\label{PROP1} Let $\theta$ be the  involution in $GL_{\ell+1}$
defining  a classical Lie group $G$ with maximal torus $H_G$ and the
Weyl group $W_{G}$.  Then the following isomorphisms hold
  \be\label{P1}
  H_{\ell+1}^{\theta_G}=H_G,
  \ee
  \be\label{P2}
\left(N_{\ell+1}(H_{\ell+1})\right)^{\theta_G}=N_G(H_G),
  \ee
  \be \label{P3}
  \left(W_{GL_{\ell+1}}\right)^{\theta_G}={\rm Aut}(\Phi_G).
  \ee
\end{prop}

\proof:   Any two $\theta$-fixed elements of $H_{\ell+1}$ define a
pair of mutually commuting elements of $G$. Thus we have an
embedding $H_{\ell+1}^\theta\subset H_G$. By Lemma \ref{PROPC} any
element commuting  with   $H_{\ell+1}^\theta$ shall be in
$H_{\ell+1}$ i.e. we have embedding  $H_G\subset H_{\ell+1}$ and
taking into account that elements of $H_G$ are $\theta_G$-invariant
we obtain $H_G\simeq H_{\ell+1}^{\theta_G}$.

To establish \eqref{P2}  first note that  elements $g\in
\left(N_{\ell+1}(H_{\ell+1})\right)^{\theta_G}$  are
$\theta$-fixed elements of $GL_{\ell+1}$ stabilizing $H_{\ell+1}$
i.e. $ghg^{-1}\in H_{\ell+1}$ for any $h\in H_{\ell+1}$. It is clear
that for any $h\in H_{\ell+1}^{\theta_G}$ we have $ghg^{-1}\in
H_{\ell+1}^{\theta_G}$. This defines an  embedding
$\left(N_{\ell+1}(H_{\ell+1})\right)^{\theta_G}\subset
N_{G}(H)$. In the following we use the fact that the involution
$\theta_G$ on the commutative group $H_{\ell+1}$ defines the
decomposition $H_{\ell+1}=H_+H_-$ into a product of commutative
subgroups, so that for each element $h\in H_{\ell+1}$ there exist
unique $h_{\pm}\in H_{\pm}$ such that
\be
h=h_+h_-=h_-h_+, \qquad
h_+^{\theta_G}=h_+,\qquad h_-^{\theta_G}=h_-^{-1}.
\ee
Using
\eqref{P1} we can make identification  $H_+\simeq H_G$.

Now we shall construct an  embedding $N_{G}(H_G)\subset
N_{\ell+1}(H_{\ell+1})^{\theta_G}$ i.e. prove that any
$\theta$-invariant element $g\in GL_{\ell+1}$ stabilizing $H_G$ also
transforms $h\in H_{\ell+1}$ into some other element $h'\in
H_{\ell+1}$. This is trivial for elements of $H_+\simeq H_G$ and
thus we need to prove that $gh_-g^{-1}$ is an element of
$H_{\ell+1}$. The element $gh_-g^{-1}$ commutes with all elements in
$H_{\ell+1}^{\theta_G}\subset H_{\ell+1}$ due to the fact that $H_+$
and $H_-$ are mutually commute and $g$ stabilizes $H_+$. Moreover
for a given $g$ all elements $gh_-g^{-1}$ mutually commute. By Lemma
\ref{PROPC} this entails that  $gh_-g^{-1}\subset H_{\ell+1}$.
Combining these two embeddings we obtain  \eqref{P2}.
To prove \eqref{P3} let us compare   two groups
\be\label{ISO}
N_{\ell+1}(H_{\ell+1})^\theta/H_{\ell+1}^\theta\,, \qquad
\left(N_{\ell+1}(H_{\ell+1})/H_{\ell+1}\right)^\theta.
\ee
There is a natural map from the first  group to the second  one  as
each element in the l.h.s. coset defines an element of $g\in
N_{\ell+1}(H_{\ell+1})$ modulo element of $H_{\ell+1}^\theta$
and therefore modulo $H_{\ell+1}$. Note that elements of the second
group in \eqref{ISO}  are defined by the condition
\be\label{THETA}
g^\theta=g\cdot h, \qquad g\in N_{\ell+1}(H_{\ell+1}), \qquad
h\in H_{\ell+1},
\ee
modulo right action of $H_{\ell+1}$ on $g$. Due
to  involutivity of $\theta_G$ we have
\be
g=\left(g^\theta\right)^\theta=(gh)^\theta=ghh^\theta,
\ee and thus
\be
h^\theta=h^{-1}.
\ee
Changing representative  $g\to
\tilde{g}=g\tilde{h}$, $\tilde{h}\in H_{\ell+1}$ for $g$ in the
coset $N_{\ell+1}(H_{\ell+1})/H_{\ell+1}$ we may also modify
\eqref{THETA} as follows
\be
\tilde{g}^\theta=g\cdot
h\tilde{h}^\theta=\tilde{g} \cdot h\tilde{h}^\theta \tilde{h}^{-1}.
\ee
Thus one can get rid of $h$ in \eqref{THETA} if we manage to
solve the equation
\be\label{EQ2}
\tilde{h}^\theta=h\tilde{h},
\ee
for any $h$. Using the decomposition
\be
\tilde{h}_-=\tilde{h}_+\tilde{h}_-\,,
\ee
and taking into account that
$h\in H_-$ we can easily check that the equation \eqref{EQ2} can be
always solved. Hence  we conclude that  two groups \eqref{ISO} are
isomorphic and thus the last statement \eqref{P3} follows from the
fact that the Weyl groups are given by the quotients
\be\nonumber
W_{GL_{\ell+1}}^{\theta_G}=
\left(N_{\ell+1}(H_{\ell+1})/H_{\ell+1}\right)^\theta
=N_{\ell+1}(H_{\ell+1})^\theta/H_{\ell+1}^\theta=N_G(H_G)/H_G={\rm
  Aut}(\Phi_G).
\ee
$\Box$

\begin{prop}\label{TitsFixedSubgroup}
Given a classical complex Lie group $G=(GL_{\ell+1}(\IC))^{\theta_G}$, the
following holds:
 \be\label{TitsG}
 \bigl(W^T_{GL_{\ell+1}}\bigr)^{\theta_G}\,
 =\,W^T_G\,.
 \ee
The explicit presentation for \eqref{TitsG} may be provided by
replacement of all $s_i$'s by $\dot{s}_i$'s in \eqref{BWeyl0},
\eqref{CWeyl0} and \eqref{DWeyl0}.

\end{prop}
\proof: To prove the assertion we use the following fact (see
\cite{GLO} Proposition 3.1). Let $G(\IR)$ be the totally  split real
form of $G(\IC)$. Then
 \be\label{TSF}
  \pi_0(N_{G(\IR)}(H(\IR)))=W^T_{G(\IC)}(H(\IC))\,.
 \ee

The connected component of the trivial element of
$W_{G(\IR)}^T\subset N_{G(\IR)}(H(\IR))$ is isomorphic to
$\IR_{>0}^{{\rm rank}(G)}$ and we have the split exact sequence
 \be
  1\longrightarrow \IR_{>0}^{{\rm rank}(G)}\longrightarrow
  N_{G(\IR)}(H_G(\IR))\longrightarrow W^T_{G(\IC)}(H_G(\IC))\longrightarrow 1\,.
 \ee
 Thus we have
 \be
 \left(W^T_{GL_{\ell+1}(\IC)}(H_{\ell+1}(\IC))\right)^{\theta_G}=
 \left(N_{\ell+1}(\IR)(H_{\ell+1})/\IR_{>0}^{\ell+1}\right)^{\theta_G}\,,
 \ee
while on the other hand
 \be\label{ISO1}
  W^T_{G(\IC)}(H_G(\IC))\,
  =\,\left(N_{\ell+1}(\IR)(H_{\ell+1})\right)^{\theta_G}\bigl/
  \left(\IR_{>0}^{\ell+1}\right)^{\theta_G}\,.
  \ee
  There is an obvious map
  \be\label{ISO22}
\left(N_{\ell+1}(\IR)(H_{\ell+1})\right)^{\theta_G}\bigl/
  \left(\IR_{>0}^{\ell+1}\right)^{\theta_G}\,\longrightarrow\,
  \left(N_{\ell+1}(\IR)(H_{\ell+1})/\IR_{>0}^{\ell+1}\right)^{\theta_G}\,,
 \ee
 and we have to show that it is actually  isomorphism.
 The obstruction to the isomorphism \eqref{ISO22} is
given by elements  in $N_{\ell+1}(\IR)(H_{\ell+1})$ invariant
with respect to $\theta_G$ only up to multiplication by an element
$\Lambda\in \IR_{>0}^{\ell+1}$ and we shall consider such elements
up to the right action of $H_{\ell+1}$. As $\theta_G$ has order two
we have the equation
 \be\label{EQ1}
  \Lambda\cdot \Lambda^{\theta_G}=1\,.
 \ee
Now consider another representative $\tilde{g}:=gQ$  in the same
$H_{\ell+1}$-coset that
 \be
  \tilde{g}^{\theta_G}=g^{\theta_G}Q^{\theta_G}=g\Lambda Q^{\theta_G}\,.
 \ee
Thus we can chose $\theta_G$-invariant representative if we can
solve the equation
 \be
  Q=\Lambda Q^{\theta_G}\,.
 \ee
and consistency of this equation follows from \eqref{EQ1}.
Explicitly we have
 \be
  q_{\ell+2-i}=\lambda_iq_i, \qquad Q={\rm diag}(q_1,q_2, \cdots,
  q_{\ell+1}), \quad \Lambda={\rm diag}(\lambda_1,\lambda_2, \cdots,
  \lambda_{\ell+1})\,,
 \ee
where $Q,\Lambda\in \IR_{>0}^{\ell+1}$. Direct check shows that such
equation is always solvable and thus we have the isomorphism
 \be\label{ISO2}
  \left(N_{\ell+1}(\IR)/(H_{\ell+1})\right)^{\theta_G}/
  \left(\IR_{>0}^{\ell+1}\right)^{\theta_G}\simeq
  \left(N_{\ell+1}(\IR)(H_{\ell+1})/\IR_{>0}^{\ell+1}\right)^{\theta_G}\,.
  \ee
This completes the proof of \eqref{TitsG}. Note that the calculations
above actually prove that the first cohomology group of $\IZ/2\IZ$
generated by $\theta_G$ acting  in $\IR_{>0}^{\ell+1}$ is trivial.

 The explicit presentation for the generators of $W^T_G$ is
obtained in Lemmas \ref{BTitsLift}, \ref{CTitsLift} and
\ref{DTitsLift} in the following subsections. $\Box$

\section{ Why $\Sp_{2\ell}$ is different}

According to Theorem \ref{MainTheorem} all classical Lie groups $G$
except of symplectic type allow a lift of the Weyl
group $W_G$ into the corresponding Lie group $G$. In other words the
exact sequence \eqref{NGH} splits, so that the map $p$ has a
section. The splitting in these cases allows a simple explanation:  the corresponding
Weyl group action in $\fh_{\IR}$ preserves
the standard inner product, and therefore the Weyl group may be
identified with a subgroup of $O(\fh)$, which immediately implies
the splitting of the sequence \eqref{NGH0} in cases of $GL_{\ell+1}$
and $O_{\ell+1}$. In contrast, for symplectic groups the exact
sequence \eqref{NGH} does not split and as a section of $p$ we encounter a non-trivial
extension $W^T_{\Sp_{2\ell+2}}$ of the Weyl group
$W_{\Sp_{2\ell+2}}=W(C_{\ell+1})$ introduced by Tits. Obviously  the
peculiarity of the symplectic series of classical Lie groups
 begs for explanation. In this section we propose an
explanation of this phenomena relying on the properties of
quaternion numbers $\IH$, the unique  non-commutative associative normed
division algebra over $\IR$.

Recall that in the case of the
classical Lie groups along with the standard approach to the classification
of the complex Lie algebras via root data  we may use another (Cartan's) approach based on
the analysis of subalgebras of the general linear algebras fixed by
certain involutions. Closely related approach is based on the
classification of \emph{normed division associative algebras}. The
list of normed division associative algebras $A$ is exhausted by the
algebras of real, complex and quaternion numbers;
to these algebras  we associate three series of general
linear groups:
 \be\label{NDA}
  GL_{\ell+1}(\IR)\,, \qquad GL_{\ell+1}(\IC)\,, \qquad GL_{\ell+1}(\IH)\,.
  \ee
The corresponding maximal compact subgroups
 \be
  O_{\ell+1}(\IR)\,, \qquad U_{\ell+1}\,,\qquad U\Sp_{\ell+1}\,,
 \ee
may be defined as stabilizer subgroups of the standard  quadratic
form over the corresponding division algebra $A$:
 \be
  (x,x)\,=\,\sum_{i=1}^{n} x_i^\dagger x_i\,=\,\sum_{i=1}^{n}
  \|x_i\|^2\,,
  \qquad x_i\in A\,.
 \ee
Further complexification of the compact groups provides a complete
list of the classical complex Lie groups:
 \be
  O_{\ell+1}(\IR)\otimes_{\IR}\IC\,=\,O_{\ell+1}(\IC)\,, \qquad
  U_{\ell+1}\otimes_\IR \IC\,=\,GL_{\ell+1}(\IC)\,, \\
  U\Sp_{\ell+1}\otimes_{\IR}\IC\,=\,\Sp_{2\ell+2}(\IC)\,.
 \ee

Actually the point of view based on  matrix groups over non-commutative
fields implies some adjustment of standard definitions in  the
structure theory of Lie groups. An obvious instance is the notion of
the maximal torus: for matrix groups over non-commutative fields it
is more natural to consider  diagonal subgroups of  general
linear groups over $A$  rather than maximal commutative subgroup.
Namely, given a normed division algebra $A$, let $A^*$ be its
multiplicative group of invertible elements; then
$H_A=(A^*)^{\ell+1}$ will be the diagonal subgroup of
$GL_{\ell+1}(A)$. In special case of symplectic groups the
underlying normed division algebra $A=\IH$ is non-commutative and
the corresponding diagonal subgroup of $GL_{\ell+1}(\IH)$,
 \be
  H_{\IH}\,
  =\,\underbrace{\bigl(\IH^*\times \cdots\times \IH^*\bigr)}_{\ell+1}\,
  \subset\,GL_{\ell+1}(\IH)\,,
 \ee
is a non-commutative Lie group. The modification of the notion
 of maximal torus implies the following  modification of the definition of the
 Weyl group.

\begin{de}\label{WeylH} Define the Weyl group $W_{GL_{\ell+1}}(A)$ of
  $GL_{\ell+1}(A)$ to be the group of
  inner automorphisms of the diagonal subgroup $H_A\subset
  GL_{\ell+1}(A)$.
\end{de}
Note that while in the commutative case the group of inner
automorphisms of the diagonal subgroup $H_{\ell+1}\subset GL_{\ell+1}$ is
isomorphic to the quotient group $N_{\ell+1}(H_{\ell+1})/H_{\ell+1}$ for
non-commutative $A$ this is not true.

\begin{lem} The group $W_{GL_{\ell+1}}(\IH)$ allows the following description:
 \be\label{HW}
  W_{GL_{\ell+1}}(\IH)\,=\,\mathfrak{S}_{\ell+1}\ltimes(SO_3)^{\ell+1}\,.
 \ee
\end{lem}

\proof: The group $GL_{\ell+1}(\IH)$ acts on its diagonal subgroup
$H_{\IH}$ by conjugation, so let us find the normalizer subgroup
$N_{\IH}=N_{\ell+1}(\IH)(H_{\IH})$ which preserves $H_{\IH}$.
The diagonal group $H_{\IH}$ is generated by one-parametric
subgroups
 \be
  H_{\IH}^{(k)}\,
  =\,\Big\{\left(\begin{smallmatrix}
  \Id_{k-1}&0&0\\0&d&0\\0&0&\Id_{\ell-k}
  \end{smallmatrix}\right)\,,\quad d\in\IH^*\Big\}\,
  \subset\,H_{\IH}\,,\qquad1\leq k\leq\ell+1\,.
 \ee
Let us pick a diagonal element
$D_k=\diag(\Id_{k-1}\,,\,d_k,\,\Id_{\ell-k})\in H^{(k)}_{\IH}$.
Then $M\in GL_{\ell+1}(\IH)$ belongs to $N_{\IH}$ if and only if it
satisfies the following equations for each $1\leq k\leq\ell+1$:
 \be
  MD_k\,=\,D'M\,,\quad\text{for some}\quad
  D'=\diag(d'_1,\ldots,d'_{\ell+1})\in H_{\IH}\,.
 \ee
Explicitly using the matrix notation $M=\|g_{ij}\|$, the above
equation reads
 \be\label{NormHid}
  g_{ij}\,=\,d_i'g_{ij}\,,\quad j\neq k\,;\qquad
  g_{ik}d_k\,=\,d_i'g_{ik}\,,\quad1\leq i,\,j\,,k\leq\ell+1\,.
 \ee
Now for each $1\leq i\leq\ell+1$ the former relation in
\eqref{NormHid} implies either $g_{ij}=0,\,\forall j\neq k$, or
$d_i'=1$. In case $g_{ij}=0,\,\forall j\neq k$ we necessarily obtain
from the latter relation in \eqref{NormHid} that $g_{ik}\neq0$,
otherwise $M$ has the $i$-th zero row and $M\notin
GL_{\ell+1}(\IH)$. In case $d_i'=1$ it follows from the latter
relation in \eqref{NormHid} that $g_{ik}d_k=g_{ik}$, which entails
$g_{ik}=0$, since we assume $d_k\neq1$ so that
$D_k\neq\Id_{\ell+1}$. Taking into account that $g_{ij}$ are
subjected to \eqref{NormHid} for each $1\leq k\leq\ell+1$, we deduce
from the above that there is only one non-zero element in the $i$-th
row of $M$.

Since the matrix $M$ is invertible its rows are linearly
independent, so for different $i$ the non-zero entries $g_{ij}$ have
different $j$'s. This defines the subgroup of monomial matrices
$\fM_{\ell+1}(\IH)\subset GL_{\ell+1}(\IH)$ consisting of matrices
with only one non-zero element in each column and each row, and
yields $N_{\IH}=\fM_{\ell+1}(\IH)$. Clearly, the quotient group
$N_{\IH}/H_{\IH}$ is isomorphic to the permutation group:
 \be
  \fM_{\ell+1}(\IH)/H_{\IH}\,=\,\fS_{\ell+1}\,,
 \ee
so that Definition \ref{WeylH} reads
 \be\label{HWeyl}
  W_{GL_{\ell+1}}(\IH)\,=\,\fS_{\ell+1}\ltimes\Int(\IH)^{\ell+1}\,.
 \ee
where $\Int(\IH)=\IH^*/\IR^*$ is the group of inner automorphisms of
$\IH$ (note that by Skolem-Noether theorem all automorphisms of
$\IH$ are inner). The norm homomorphism induces the following exact
sequence:
 \be\label{NormH}\xymatrix{
  1\ar[r]& SU_2\ar[r]& \IH^*\ar[r]^{\N}& \IR_{>0}\ar[r]&
  1},\\
  \N(q)\,=\,q\bar{q}\,,\qquad \N(\IH)\,=\,\IR_{>0}\,.
 \ee
Taking into account $\IR^*=\CZ(\IH)\cap\IH^*=\IR_{>0}\times\mu_2$
with $\mu_2=\{\pm1\}$ being a group of roots of 1 in $\IR$, we
obtain
 \be\label{IntH}
  \Aut(\IH)\,=\,\IH^*/\IR^*\,
  =\,SU(2)/\mu_2\,\simeq\, SO_3\,.
 \ee
Thus $\Aut(\IH)$ is identified with the adjoint group $SO_3$ of the
group $SU_2$ of unite quaternions. This complete the proof of
\eqref{HW}. $\Box$

Let $\Aut_{\IC}(\IH)\subset\Aut(\IH)$ be the subgroup of
automorphisms preserving the subalgebra
$\IC=(\IR\oplus\IR\imath)\subset\IH$.
\begin{lem} The following holds:
 \be \label{AUT}
  \Aut_{\IC}(\IH)\,=\,(\IC^*\sqcup\IC^*\jmath\,)/\IR^*\,,
 \ee
and the following central extension splits
 \be\label{AUT1}
  1\longrightarrow\IC^*/\IR^*\longrightarrow\Aut_{\IC}(\IH)
  \longrightarrow\Gal(\IC/\IR)\longrightarrow1\,.
 \ee
\end{lem}
\proof: Consider the tautological faithful two-dimensional representation
$\IH\to\Mat_2(\IC),\,q\mapsto\hat{q}$:
 \be
  \hat{1}\,=\,\Id_2\,,\qquad
  \hat{\imath}\,
  =\,\Big(\begin{smallmatrix}
  \imath&&0\\&&\\0&&-\imath
  \end{smallmatrix}\Big)\,,\qquad
  \hat{\jmath}\,
  =\,\Big(\begin{smallmatrix}
  0&&1\\&&\\-1&&0
  \end{smallmatrix}\Big)\,,\qquad
  \hat{k}\,
  =\,\Big(\begin{smallmatrix}
  0&&\imath\\&&\\ \imath&&0
  \end{smallmatrix}\Big)\,.
 \ee
The norm homomorphism \eqref{NormH} is given by the matrix
determinant:
 \be
  q\,=\,q_0+\imath q_1+\jmath q_2+kq_3 \in \IH\,, \qquad
  |q|^2\,=\,q\bar{q}\,=\,q_0^2+q_1^2+q_2^2+q_3^2\,,\\
  \hat{q}\,=\,\begin{pmatrix} q_0+\imath q_1 & q_2+\imath q_3\\
  -(q_2-\imath q_3) & q_0-\imath q_1
  \end{pmatrix}, \qquad |q|^2=\det \hat{q}\,\in\,\IR_{>0}\,.
 \ee
Then elements $z$ of the subalgebra $\IC=\IR\oplus\IR\imath\subset
\IH$ are identified with the diagonal  matrices
 \be
  z\,\longmapsto\,\hat{z}\,
  =\,\Big(\begin{smallmatrix}
  z&&0\\&&\\0&&\bar{z}
  \end{smallmatrix}\Big)\,,\qquad
  |z|^2\,=\,z\bar{z}\,=\,\det\hat{z}\,,
  \ee
while the general quaternion has representation for $a,\,b\in\IC$
 \be
  q\,\longmapsto\,\hat{q}\,
  =\,\Big(\begin{smallmatrix}
  a&&b\\&&\\-\bar{b}&&\bar{a}
  \end{smallmatrix}\Big)\,,\qquad
\det\hat{q}=|a|^2+|b|^2\,.
  \ee
Now a quaternion $q$ belongs to centralizer of subalgebra
$\IC\subset\IH$ if and only if $\hat{q}\hat{z}=\hat{z}'\,\hat{q}$
for any $z\in\IC$, which reads:
  \be
   \Big(\begin{smallmatrix}
  az&&b\bar{z}\\&&\\-\bar{b}z&&\bar{a}\bar{z}
  \end{smallmatrix}\Big)\,
  =\,\Big(\begin{smallmatrix}
  z'a&&z'b\\&&\\-\bar{z}'\bar{b}&&\bar{z}'\bar{a}
  \end{smallmatrix}\Big)\,.
  \ee
  Since $z\in\IC$ is arbitrary, either $z=z'$ and $b=0$, or
  $z=\bar{z}'$ and $a=0$.
  Obviously, the case $b=0$ corresponds to
$q\in\IC^*\subset\IH^*$, and in case $a=0$ we have
$q\in\IC^*\jmath\subset\IH^*$. Taking into account that $q\in \IR^*$
implies \eqref{AUT}. The second statement \eqref{AUT1} follows from
the fact that $\hat{\jmath}$ acts by complex conjugation:
 \be\label{Jaction}
  \jmath\,z\,\jmath^{-1}\,
  =\,\bar{z}\,,
  \ee
and its image in the quotient group
$(\IC^*\sqcup\IC^*\jmath\,)/\IR^*$ is precisely the generator of
$\Gal(\IC/\IR)$ given by the complex conjugation. $\Box$

Recall that  standard Weyl group of the symplectic Lie group
$\Sp_{2\ell+2}$ is given by
 \be
 W_{\Sp_{2\ell+2}}\,=\,W(C_{\ell+1})\,=\,
 \mathfrak{S}_{\ell+1}\ltimes(\IZ/2\IZ)^{\ell+1}\,.
  \ee
Here $\IZ/2\IZ$ in r.h.s. may be identified with the Galois group ${\rm
  Gal}(\IC/\IR)$ and thus allowing via splitting of \eqref{AUT1} an
embedding
 \be
W_{\Sp_{2\ell+2}}\subset   W_{GL_{\ell+1}}(\IH).
\ee

The quaternionic analog  $W_{GL_{\ell+1}}(\IH)$ of the standard Weyl
group has obvious obstruction for
the embedding into the general Lie group $GL_{\ell+1}(\IH)$. Indeed
this lift implies in particular transition from the adjoint action
of quaternions on itself to its left action. On the other hand we have
the following standard exact sequence
\be\label{SU20}
\xymatrix{1\ar[r]& \IR^*\ar[r]&
  \IH^*\ar[r]^{\pi}& {\rm Aut}(\IH)\ar[r]&  1},
 \ee
or taking into account the identifications $SU(2)=\IH^*/\N(\IH)$
(for the notations see  \eqref{NormH}) and ${\rm Aut}(\IH)=SO_3$
\be\label{SU}
\xymatrix{1\ar[r]& \mu_2\ar[r]&
  SU_2\ar[r]^{\pi}& SO_3(\IR)\ar[r]&  1}.
 \ee
The generator of $\mu_2$ above   might be identified with the square
$\jmath^2$ of the quaternionic unity. Indeed, the
group ${\rm Aut}(\IH)\simeq SO_3$ acts in the
subspace $\IR^3$ of purely imaginary  quaternions
 \be
  \hat{x}\,=\,\imath x_1\,+\,\jmath x_2\,+\,kx_3\,,
 \ee
via standard three-dimensional rotations. Then the section of
projection $\pi$ may be identified with the lift of the orthogonal
rotations to the conjugation action of the unit norm quaternions
$q\in SU_2$:
 \be
  \Ad_q\,:\quad\hat{x}\longrightarrow q\hat{x}q^{-1}\,.
 \ee
Let is pick a rotation given by the diagonal element
 \be\label{SAMPL}
  g\,=\,\diag(-1,+1,-1)\,\in\, SO_3\,.
 \ee
Solving the equation $\widehat{gx}\,=\,q\widehat{x}q^{-1}$ for
\eqref{SAMPL}, that is
 \be
 -\imath x_1\,+\,\jmath x_2\,-\,kx_3\,=
 \,q\bigl(\imath x_1\,+\,\jmath x_2\,+\,kx_3\bigr)\bar{q}\,,
 \ee
 we derive that $q=\jmath$ and thus
 \be
  \hat{g}\,=\,\jmath\,,\qquad \hat{g}^2=-1\,.
 \ee
The fact that the square is equal minus unite element  implies
that indeed we have a central extension with generator which may be
identified with $\jmath^2$.

Now taking into account \eqref{HWeyl} we arrive at the following result.

\begin{prop}
  The quaternionic Weyl group $W_{GL_{\ell+1}}(\IH)$ allows  the extension
$\widetilde{W}_{GL_{\ell+1}}(\IH)$
\be\label{QWEYL}
1\longrightarrow
  (\mu_2)^{\ell+1}\longrightarrow \widetilde{W}_{GL_{\ell+1}}(\IH)
  \longrightarrow W_{GL_{\ell+1}}(\IH)\longrightarrow 1\,.
 \ee
The extended group $\widetilde{W}_{GL_{\ell+1}}(\IH)$ allows a
natural embedding
 \be
  \widetilde{W}_{GL_{\ell+1}}(\IH)\subset GL_{\ell+1}(\IH)\,.
 \ee
The Tits group extension
\be\label{QTITS}
1\longrightarrow
  (\mu_2)^{\ell+1}\longrightarrow W^T_{\Sp_{2\ell+2}}
  \longrightarrow W_{\Sp_{2\ell+2}}\longrightarrow 1\,,
 \ee
is then obtained by restriction of
\eqref{QWEYL} to the subgroup $W_{\Sp_{2\ell+2}}\subset
W_{GL_{\ell+1}}(\IH)$.
\end{prop}

Thus the underlying  reason for the appearance of the Tits groups in the case
of the symplectic Lie groups may be traced back to the
extensions \eqref{SU20},  \eqref{SU}. The non-triviality of this
extension is closely related with the basic $\mu_2$-extension of the
Galois group ${\rm Gal}(\IC/\IR)$ characterizing quaternions.
Actually the cohomology class measuring non-triviality of this extension
is directly related with the basic invariant characterizing the real
central simple division algebras.  Indeed  such algebras are characterized by the
invariant taking values in the 2-torsion of the Brauer group (see
e.g.  \cite{S}, \cite{PR}), which coincides with the whole Brauer group in case of $\IR$:
\be\label{BRR}
\Br(\IR)\,=\,H^2(\Gal(\IC/\IR),\IC^*)\simeq H^2(\Gal(\IC/\IR),\mu_2) \simeq \mu_2 \,.
\ee
The corresponding cohomology group describes the extensions of the
form
\be\label{AEXT}
1\longrightarrow \IC^*\longrightarrow
 A\longrightarrow\Gal(\IC/\IR)\longrightarrow 1\,,
 \ee
and the case of the  semidirect product $\IC\rtimes\Gal(\IC/\IR)$
with the natural action of $\Gal(\IC/\IR)$ via complex conjugation
corresponds to the trivial cocycle. The only other non-trivial case
corresponds to the algebra $A=\IH$ of quaternions so that  the lift
of the generator $\sigma\in  \Gal(\IC/\IR)$ may be identified
with  the  quaternionic unite   $\jmath$, $\jmath^2=-1$. This
provides a link between the non-trivial extension \eqref{SU} and the
construction of quaternions via the non-trivial extension
\eqref{AEXT} (for more conceptual explanations of the relation
between universal $\mu_2$-extension of the orthogonal groups and the
second cohomology of the Galois group of the base field see
\cite{BD} and references therein). In particular this directly
relates  the fact that $\jmath^2=-1$ in the algebra of quaternions
with the unavoidable appearance of the Tits extension of the  Weyl
group for classical symplectic Lie groups.

\section{Lie groups closely related with classical Lie groups}

In this section we consider several examples of construction of
suitable  sections of the extensions \eqref{NGH0} for the Lie
groups closely relate with classical Lie groups. Precisely we will
treat the cases of unimodular subgroup $SL_{\ell+1}\subset
GL_{\ell+1}$ of the general linear group, of the  groups ${\rm
Pin}_{2\ell+1}$,  ${\rm Pin}_{2\ell}$,  ${\rm
Spin}_{2\ell+1}$ and ${\rm Spin}_{2\ell}$. In all cases we provide
explicit construction  of  sections of \eqref{NGH0} realized as
central extensions of the corresponding Weyl groups. Let us stress
that these sections
differ from the Tits lifts (which are not central extensions of the
corresponding Weyl groups in general). Although the possibility to
define sections via central extensions is  obvious for Pin groups
(they are central extensions of the corresponding classical
orthogonal groups) it is a bit less obvious for unimodular  and spinor
Lie groups.

\subsection{Construction of a section for $G=SL_{\ell+1}$}

We will freely use the notations for root data of type $A_{\ell}$
defined in the Appendix.   Recall that the unimodular
subgroup $SL_{\ell+1}\subset GL_{\ell+1}$ is defined as a kernel of
the determinant map thus fitting the following exact sequence:
\be\label{DETes} \xymatrix{ 1\ar[r] &
SL_{\ell+1}(\IC)\ar[r]^{\varphi} &
  GL_{\ell+1}(\IC)\ar[r]^{\det} &  \IC^*\ar[r] & 1}.
 \ee
The group center of unimodular group is identified with the cyclic
group: $\CZ(SL_{\ell+1}(\IC))=\IZ/(\ell+1)\IZ$ generated by
$\zeta=e^{2\pi\imath\vp_{\ell}^{\vee}}$. The image $\varphi(\zeta)$
belongs to the center $\CZ(GL_{\ell+1})=\IC^*$ and is given by:
 \be
  \varphi(\zeta)\,=\,e^{\frac{2\pi\imath}{\ell+1}\,p^{\vee}_{\ell+1}}\,
  \in\,\CZ(GL_{\ell+1}(\IC))\,.
 \ee
In Section \ref{GLNORM} we construct a lift of the Weyl  group
$W_{GL_{\ell+1}}=W(A_{\ell})$ corresponding to the root system
$A_{\ell}$ to the group $GL_{\ell+1}$. As the Weyl groups for
$SL_{\ell+1}$ and $GL_{\ell+1}$ coincide to construct a section of
\eqref{NGH0} for $SL_{\ell+1}(\IC)$ we shall invert the homomorphims
$\varphi$ on the subgroup $W_{GL_{\ell+1}}\subset GL_{\ell+1}$.
Formally the inverse of $\varphi$ may be written as follows \be
  \varphi^{-1}(g)\,=\,g\cdot\bigl[\det{}^{(-1)}(\det
  g)\bigr]^{-1}\,,\qquad g\in GL_{\ell+1}(\IC)\,,
 \ee
where $\det^{(-1)}$ is a multi-valued inverse of the determinant map $\det$ in
\eqref{DETes}, which can be defined by
 \be
  \det{}^{(-1)}\,:\quad\IC^*\,\longrightarrow\,\CZ(GL_{\ell+1})\,,\qquad
  \det{}^{(-1)}(z)\,=\,e^{\frac{1}{\ell+1}\,p^{\vee}_{\ell+1}\,\log(z)}\,,\quad
  z\in\IC^*\,.
  \ee
To make this map single-valued we shall choose a branch
of the $\log$-function. Note that all generators $S_i$ of
$W_{GL_{\ell+1}}(\IC)$ defined in \eqref{permat} satisfy the
relation $\det S_i=-1$. Thus we  pick $\log(-1)=\imath \pi$ and
introduce the following lift of generators $S_i$ in
$SL_{\ell+1}(\IC)$:
 \be\label{LIFTSL}
 \sigma_i\,=
 \,e^{\frac{\pi\imath}{\ell+1}\,p^{\vee}_{\ell+1}}\,S_i\,,\qquad1\leq i\leq\ell\,.
 \ee

\begin{prop} The group $\widetilde{W}(A_{\ell})$ generated by the elements
  \eqref{LIFTSL}   provides a central extension
 \be\xymatrix{
  1\ar[r]& \CZ(SL_{\ell+1})\ar[r]& \widetilde{W}(A_{\ell})\ar[r]&
  W(A_{\ell})\ar[r]&1}\,,
 \ee
 of the Weyl group $W(A_{\ell})$   with the defining relations
 \be\label{CR1}
  \s_k^2\,=\,\zeta\,,\quad\text{for each}\quad1\leq k\leq\ell\,;\\
  \s_k\s_j\,=\,\s_j\s_k\,,\quad\text{if}\quad
  m_{kj}=2,\\
  \text{and}\qquad
  \s_k\s_j\s_k\,=\,\s_j\s_k\s_j\,,\quad\text{if}\quad m_{kj}=3\,,
  \ee
where $\zeta=e^{2\pi\imath\omega^\vee_\ell }$ is the generator of
the center $\CZ(SL_{\ell+1})$.
The group $\tilde{W}(A_{\ell})$ allows an embedding in $N_{SL_{\ell+1}}(H_{SL_{\ell+1}})$
compatible with the exact sequence
 \be
  1\longrightarrow H_{SL_{\ell+1}}\longrightarrow
  N_{SL_{\ell+1}}(H_{SL_{\ell+1}})
  \longrightarrow
  W(A_{\ell}) \longrightarrow 1\,,
 \ee
and the  embedding $\CZ(SL_{\ell+1})\subset H_{SL_{\ell+1}}$.

\end{prop}

\proof.  The element $e^{\frac {\imath \pi}{\ell+1} p_{\ell+1}^{\vee}}$ is in
the center of $GL_{\ell+1}(\IC)$ and its square is equal to $\zeta$.
The relations \eqref{CR1} then follow from the relations
\eqref{permat} and \eqref{WeylPerm} for the generators of the Weyl
group $W(A_\ell)$. $\Box$

\subsection{Construction of a  section for $SO_{\ell+1}$}\label{SOWEYL}

The problem of lifting  the Weyl groups $W(B_{\ell})$ and
$W(D_{\ell})$ into special orthogonal subgroups is similar to that
for $SL_{\ell+1}\subset GL_{\ell+1}$. In \eqref{2torsB0},
\eqref{2torsD0} we introduce a lift of Weyl groups $W(B_{\ell})$ and
$W(D_{\ell})$ to the corresponding orthogonal groups. It is natural
to split this construction into two parts depending on the parity of
the rank.

In case of the odd orthogonal group the following sequence splits,
 \be\label{DetOrthogonal} \xymatrix{
  1\ar[r] & SO_{2\ell+1}\ar[r]^{\varphi} & O_{2\ell+1}\ar[r]^{\det} & \mu_2\ar[r] &
  1}\,,
 \ee
so that
 \be
  O_{2\ell+1}\simeq\CZ(O_{2\ell+1})\times SO_{2\ell+1}\,.
 \ee
Here the generator  $z$ of the  center $\CZ(O_{2\ell+1})=\IZ/2\IZ$
is given by (see Section \ref{APPB} for notations):
 \be\label{CENG}
  z\,=\,T^B_0T^B_1\cdots T^B_{\ell}\,
  =\,T^B_0e^{\pi\imath\vp^{\vee}_1}
  =-\Id_{2\ell+1}\in\CZ(O_{2\ell+1})\,,
 \ee
The Weyl group generators $S^B_i$ \eqref{2torsB0} satisfy
 \be
  \det(S^B_1)\,=\,-1\,,\qquad\det(S^B_k)\,=\,1\,,\quad1<k\leq\ell\,,
 \ee
so that $S^B_1\in(O_{2\ell+1}\setminus SO_{2\ell+1})$ and $S^B_k\in
SO_{2\ell+1},\,1<k\leq\ell$.

 \begin{lem} The elements $\s_i\in SO_{2\ell+1}$ defined by
 \be\label{LIFTB}
  \s_1\,=\,zS^B_1\,,\qquad
  \s_k\,=\,S^B_k\,,\quad1<k\leq\ell\,,
 \ee
generate the group isomorphic to the Weyl group $W(B_{\ell})$.
\end{lem}
\proof:  Follows from the fact that $z^2=1$ in $SO_{2\ell+1}$.
$\Box$

Note that here we follow  the analogous  construction for the  case
of unimodular group $SL_{\ell+1}\subset GL_{\ell+1}$ by using
\eqref{CENG}  as a  lift of the generator of $\mu_2$ in
\eqref{DetOrthogonal}.

In the case of the even orthogonal group generators $S^D_i$
\eqref{2torsD0} of $W(D_{\ell})$ already belong to $SO_{2\ell}$
due to $\det(S^D_i)=1$, and therefore provide an embedding
$W(D_{\ell})\subset N_{SO_{2\ell}}$.

\subsection{Construction of a  section for $\Pin_{\ell+1}$}\label{PINWEYL}

The group $\Pin(V)$ is a central extension of $O(V)$, which can be
defined as follows. Let $V$ be real vector space supplied with the
standard quadratic form
 \be\label{QuadraticSpace}
  \,\|v\|^2\,=\sum_{i=1}^{\dim(V)} v_i^2,\qquad v=(v_1,\ldots , v_{\dim(V)})\in V\,.
 \ee
Let $\CC(V)$ be the corresponding Clifford algebra, defined as a
quotient of the tensor algebra $T(V)$ as follows:
 \be\label{CliffordAlgebra}
  \CC(V)\,=\,T(V)\Big/\bigl(v\cdot v\,-\,\|v\|^2\cdot1\bigr)\,.
 \ee
The standard $\IZ$-grading on the tensor algebra induces a
$\IZ/2\IZ$-grading $\alpha: \CC(V)\to \CC(V)$ of the Clifford
algebra, thus splitting it into the direct sum:
 \be
  \CC(V)\,=\,\CC(V)^+\,\oplus\,\CC(V)^-\,.
 \ee
 Let $\top: \CC(V)\to\CC(V)$ be the transposition
anti-automorphism of the Clifford algebra, defined by $(v_1\cdot
v_2)^\top=v_2\cdot v_1$ for any $v_1,\,v_2\in V\subset\CC(V)$.
Define the conjugation anti-automorphism $x\to
\bar{x}=\alpha(x)^{\top}$ on $\CC(V)$. The spinor norm is then given
by:
 \be
  \N\,:\quad\CC(V)\,\longrightarrow\,\IR^*\,,\qquad
  x\,\longmapsto\,\N(x)\,=\,x\,\bar{x}\,,
 \ee
so that for all $u,\,v\in V$ we have
 \be\label{SpinorNorm}
  (u,v)\,=\,\frac{1}{2}\Big(\N(u)\,+\,\N(v)\,-\N(u+v)\Big)\,
  =\,\frac{u\cdot v\,+\,v\cdot u}{2}\,,\\
  \N(v)\,=-\,v\cdot v\,=\,-\|v\|^2\,=\,-(v,v)\,\,.
  \ee
\begin{de}(see \cite{ABS}) The Clifford group $\Gamma(V)$ is the
  subgroup of the group of invertible elements of $\CC(V)$ that
  respects the linear subspace $V\subset \CC(V)$ i.e.
  \be
   \Gamma(V)\,
   =\,\bigl\{x\in(\CC(V)\setminus\{0\})\,|\, xV\alpha(x)^{-1}\subseteq
   V\bigr\}\,\,.
 \ee
\end{de}
The action of $\Gamma(V)$ on $V$ defined above respects the norm
\eqref{SpinorNorm} and thus allows homomorphism to $O(V)$ for all
$u\in V\subset\CC(V)$:
 \be\label{ReflectionAction}
  u\,\longmapsto\,s_u(v)\,:=\,u\cdot v \cdot\a(u)^{-1}\,
  =\,v\,-\,2\frac{(u,v)}{(u,u)}u\,,\qquad
  \forall\,v\in V\,.
  \ee
Moreover, since the orthogonal group is generated  by reflections
with respect to elements $u\in V$, the group $\Gamma(V)$ maps onto
$O(V)$. Clearly, $\IR^*\subset\Gamma(V)$ acts
 trivially and this leads to the following exact sequence:
 \be\label{CliffordGroup}\xymatrix{
  1\ar[r]& \IR^*\ar[r]& \Gamma(V)\ar[r]^{\pi} &
  O(V)\ar[r] & 1}\,.
 \ee
Restriction to the subgroup  of elements with the norm in
$\mu_2=\{\pm 1\}$ results in taking a quotient over the subgroup
$\N(\IR^*)=\IR_{>0}$, which yields the central extension $\Pin(V)$
of $O(V)$:
 \be\label{PinOrthogonalCentralExt}\xymatrix{
  1\ar[r]& \mu_2\ar[r]& \Pin(V)\ar[r]^{\pi}& O(V)\ar[r]& 1}\,,
 \ee
where $\mu_2=\{\pm 1\}$  is  a subgroup of the center $\CZ(\Pin(V))$
of the pinor group $\Pin(V)$.

Taking into account \eqref{ReflectionAction} we may  define a
 properly normalized lift $\hat{u}$ of a simple reflection $s_u$ with respect
to a vector $u\in V$ into $\Pin(V)$ as follows
 \be\label{LIFT}
  \hat{u}:=\frac{u}
 {\sqrt{|\N(u)|}}\,,\quad\hat{u}^2\,=\,1\,,\quad
 \N(\hat{u})\,=\,-1\,,\qquad\forall u\in V\subset \CC(V)\,.
\ee

\begin{lem} Let $\{\e_1,\ldots,\e_{\dim(V)}\}\subset V$ be an orthonormal
basis, let $T_k,\,1\leq k\leq\dim(V)$ be reflections at $\e_k\in V$
and let $S_i,\,1\leq i\leq(\dim(V)-1)$ be reflections at simple root
$\e_i-\e_{i+1}$. Then the following elements of $\Pin(V)$ represent
liftings of the reflections $S_i,\,T_k\in O(V)$:
 \be\label{PINGL0}
  \CT_k\,=\,\hat{\e}_k\,,\qquad(\CT_k)^2\,=\,1,\,\qquad1\leq k\leq\dim(V)\,,\\
  \CS_i\,
  =\,\frac{\hat{\e}_i-\hat{\e}_{i+1}}{\sqrt{2}}\,,\qquad
  (\CS_i)^2\,=\,1\,,\qquad1\leq
  i\leq\ell\,.
 \ee
\end{lem}
\proof: One shall check that the action of the lifts $\CT_i$ and
$\CS_i$ on  $V$ according to \eqref{ReflectionAction} coincides
with the action of $T_i$  and $S_i$. Thus for $\CT_i$ we have the
following:
 \be\nonumber
  \pi(\CT_k)(v)\,
  =\,\CT_k\cdot  v \cdot \alpha(\CT_k)^{-1}\,=\,\hat{\e}_k\cdot
  v\cdot(-\hat{\e}_k)\\
  =\,\hat{\e}_k\cdot(x_1\hat{\e}_1+\ldots+x_k\hat{\e}_k+\ldots
  +x_{\dim(V)}\hat{\e}_{\dim(V)})\cdot(-\hat{\e}_k)\,
  =\,T_kv\,.
 \ee
Similarly the action of $\CS_i$ can be verified. $\Box$

\begin{lem}\label{LIFT23}
  The elements \eqref{PINGL0} satisfy the following
relations (compare with \eqref{WeylExt1}):
 \be\label{TitsExt}
  \CT_i\CT_j\,=\,-\CT_j\CT_i\,,\quad i\neq j\,;\qquad
  \CS_i\CS_j\,=\,-\CS_j\CS_i\,,\quad|i-j|>1\,,\\
  \CS_i\CS_j\CS_i\,=\,\CS_j\CS_i\CS_j\,,\qquad|i-j|=1\,,\\
  \CT_i\CS_i\,=\,-\CS_i\CT_{i+1}\,,
 \ee
and therefore provide a central extension of the Tits group
$W^T_{GL(V)}$ (and of its subgroup $W_{GL(V)}$ generated by $\CS_i$,
$1\leq i\leq\dim(V)-1$) by $\mu_2\subseteq\CZ(\Pin(V))$.
\end{lem}

\proof: The relations can be easily checked by straightforward
computation. For example let us verify the 3-move braid relation: one has
 \be
  \CS_i\CS_{i+1}\,
  =\,\frac{\hat{\e}_i\hat{\e}_{i+1}-\hat{\e}_i\hat{\e}_{i+2}+
    \hat{\e}_{i+1}\hat{\e}_{i+2}-1}{2}\,,\\
  (\CS_i\CS_{i+1})^2\,
  =\,\frac{-\hat{\e}_i\hat{\e}_{i+1}+\hat{\e}_i\hat{\e}_{i+2}-
    \hat{\e}_{i+1}\hat{\e}_{i+2}-1}{2}\,.
 \ee
This implies the Coxeter relation $(\CS_i\CS_j)^3=1$ for
$|i-j|=1$, which is equivalent to the 3-move braid relation due to
$\CS_i^{-1}=\CS_i$.\,\,  $\Box$

Now let us construct a lift of the generators of Weyl group
$W_{O(V)}$ into $\Pin(V)$.  We use a version of
the  embedding that is based on the expressions given in Theorem
\ref{MainTheorem} for the lift of $W_{O(V)}$ to $O(V)$. 
Note that the conjugated generators
$\overline{S}_i^A$ are expressed thorough $S_i^A$ and $T_i^A$
according to \eqref{permat1} and the lift of the generators $S_i^A$
and $T_i^A$ is provided by Lemma \ref{LIFT23}. We will 
consider the cases of $W(B_{\ell})$  and $W(D_{\ell})$ separately.

In case odd orthogonal group we have $O(V)=O_{2\ell+1}$. By
\eqref{2torsB0} Weyl group $W(B_\ell)$ can be imbedded into the
$\theta_B$-invariant subgroup $(W^T_{GL_{2\ell+1}})^{\theta_B}$ via
 \be
  S_1^B\,=\,S_{\ell}S_{\ell+1}S_{\ell}\,
  =\,S_{\ell+1}S_{\ell}S_{\ell+1}\,,\qquad
  S_k^B\,=\,S_{\ell+1-k}\overline{S}_{\ell+k}\,\,,\quad1<k\leq\ell\,.
 \ee
 Using \eqref{permat1} and  \eqref{PINGL0}  we
obtain the following:
 \be\label{PinWEYLgens1}
  \CS_1^B\,=\,\CS_{\ell}\CS_{\ell+1}\CS_{\ell}\,,\qquad
  \CS_k^B\,
  =\,\CS_{\ell+1-k}\CT_{\ell+k}\CS_{\ell+k}\CT_{\ell+k}\,,
  \quad1<k\leq\ell\,.
 \ee

In case odd orthogonal group we have $O(V)=O_{2\ell}$. By
\eqref{2torsD0} Weyl group $W(D_\ell)$ can be imbedded into the
$\theta_D$-invariant subgroup $(W^T_{GL_{2\ell}})^{\theta_D}$ via
 \be
  S_1^D\,=\,S_{\ell}S_{\ell-1}\overline{S}_{\ell+1}S_{\ell}\,\,,\qquad
  S_k^D\,=\,S_{\ell+1-k}\overline{S}_{\ell-1+k}\,\,,\quad1<k\leq\ell\,\,.
 \ee
Therefore, using \eqref{permat1} and \eqref{PINGL0} for the lift to
$\Pin_{2\ell}$ we obtain
 \be\label{PinWEYLgens2}
  \CS_1^D\,
  =\,\CS_{\ell}\CS_{\ell-1}\CT_{\ell+1}
  \CS_{\ell+1}\CT_{\ell+1}\CS_{\ell}\,,\\
  \CS_k^D\,
  =\,\CS_{\ell+1-k}\CT_{\ell-1+k}\CS_{\ell-1+k}
  \CT_{\ell-1+k}\,,\quad1<k\leq\ell\,\,.
 \ee

The  subgroups in $\Pin_{2\ell+1}$ and $\Pin_{2\ell}$  generated by
the elements $\CS_i^B$ and $\CS^D_i$, respectively, appears to be
central group extensions of the corresponding Weyl groups.

\begin{prop}\label{LIFTPin} The elements of pinor group $\Pin(V)$
introduced in
  \eqref{PinWEYLgens1}, \eqref{PinWEYLgens2} satisfy the
following relations for any $i,\,j\in I$:
 \be\label{PinWEYLBrels}
  (\CS^B_1)^2\,=\,1\,,\qquad(\CS^B_k)^2\,=\,-1\,,\quad1<k\leq\ell\,,\\
  \CS^B_i\CS^B_j\,=\,\CS^B_j\CS^B_i\,,\qquad a_{ij}=0\,,\\
  \CS^B_i\CS^B_j\CS^B_i\,=\,\CS^B_j\CS^B_i\CS^B_j\,,\qquad a_{ij}a_{ji}=1\,,\\
  \CS^B_i\CS^B_j\CS^B_i\CS^B_j\,=\,\CS^B_j\CS^B_i\CS^B_j\CS^B_i\,,\qquad a_{ij}a_{ji}=2\,,
 \ee
and
 \be\label{PinWEYLDrels}
  (\CS^D_i)^2\,=\,-1\,,\quad\CS^D_1\CS^D_2\,=\,-\CS^D_2\CS^D_1\,,\quad
  \CS^D_1\CS^D_3\CS^D_1\,=\,-\CS^D_3\CS^D_1\CS^D_3\,,\\
  \CS^D_i\CS^D_j\,=\,\CS^D_j\CS^D_i\,,\quad a_{ij}=0\,,\qquad i,j\neq (1,2)\,,\\
  \CS^D_i\CS^D_j\CS^D_i\,=\,\CS^D_j\CS^D_i\CS^D_j\,,\qquad a_{ij}a_{ji}=1\,,\quad i,j\neq (1,3)\,.
 \ee
This results in a central extension of the Weyl group $W_{O(V)}$ by
$\mu_2\subseteq\CZ(\Pin(V))$:
 \be\label{PinWEYL}\xymatrix{
  1\ar[r]& \mu_2\ar[r]& \widetilde{W}_{O(V)}\ar[r]^{\pi}& W_{O(V)}\ar[r]& 1}\,.
 \ee\end{prop}
\proof: Relations in the first lines of \eqref{PinWEYLBrels},
\eqref{PinWEYLDrels} follow from $\hat{\e}_i^2=1$ and
$\hat{\e}_i\cdot\hat{\e}_j\,=\,-\hat{\e}_j\cdot\hat{\e}_i$ for
$i\neq j$. The other relations may be verified by writing explicitly
from \eqref{PINGL0}:
 \be
  \CS^B_1\,=\,\frac{\hat{\e}_{\ell}-\hat{\e}_{\ell+2}}{\sqrt{2}}\,,\qquad
  \CT_i\CS_i\CT_i^{-1}\,
  =\,\frac{\hat{\e}_i+\hat{\e}_{i+1}}{\sqrt{2}}\,,\\
  \CS_1^D\,
  =\,\frac{\hat{\e}_{\ell-1}\hat{\e}_{\ell}+\hat{\e}_{\ell-1}\hat{\e}_{\ell+2}
  +\hat{\e}_{\ell}\hat{\e}_{\ell+1}-\hat{\e}_{\ell+1}\hat{\e}_{\ell+2}}{2}\,.
 \ee
In particular, one can verify in straightforward way that
 \be
  \CS^D_1\CS^D_3\CS^D_1\,
  =\,\frac{\hat{\e}_{\ell-2}\hat{\e}_{\ell}+\hat{\e}_{\ell}\hat{\e}_{\ell+1}
  +\hat{\e}_{\ell+1}\hat{\e}_{\ell+3}-\hat{\e}_{\ell-2}\hat{\e}_{\ell+3}}{2}\,
  =\,-\CS^D_3\CS^D_1\CS^D_3\,.
 \ee
The other 3-move braid relations from the last lines of
\eqref{PinWEYLBrels}, \eqref{PinWEYLDrels} follow by
\eqref{TitsExt}.
 $\Box$

\subsection{Construction of a section for $\Spin_{\ell+1}$}\label{SPINWEYL}

In this section we describe a lift of the Weyl group $W_{O(V)}$ into
the spinor group $\Spin(V)$ combining the results of previous
Sections \ref{SOWEYL} and \ref{PINWEYL} with the constructions of
Section \ref{MAINRES}. The construction will be compatible with the
following commutative diagram:
 \be\label{SPin}\xymatrix{
  & \mu_2\ar[d]\ar@{=}[r] & \mu_2\ar[d] &&\\
  1\ar[r] & \Spin(V)\ar@{->>}[d]^{\pi}\ar[r]^{\varphi} &
  \Pin(V)\ar@{->>}[d]^{\pi}\ar[r]^{\det}
   &\mu_2\ar@{=}[d]\ar[r]& 1\\
  1\ar[r] & \SO(V)\ar[r]^{\varphi} & O(V)\ar[r]^{\det} & \mu_2\ar[r] & 1}.
 \ee
Note that the elements of $\Spin(V)\subset \Pin(V)$ are
single out by the condition that they are represented by a product
of even number of elements $\hat{u}\in\Pin(V)$, or equivalently,
$\Spin(V)=\Pin(V)\cap \CC(V)^+$. Moreover we have epimorphism
$\Spin(V)\to SO(V)$. As we already manage to construct inverse image
under $\varphi: SO(V)\to O(V)$ of the generators of the
corresponding Weyl groups we can use \eqref{LIFT} to  lift
generators \eqref{LIFTB} and generators $S^D_i$ of $W(D_\ell)$ into
the corresponding spinor groups.
We again consider the cases of odd and even orthogonal groups
separately.

\begin{lem}\label{SPINLIFTB} Introduce the following elements
  $\wt{\CS}^{B}_i\in\Spin_{2\ell+1}$
 \be\label{LIFTSpin}
  \wt{\CS}^{B}_1\,=\,\hat{z}\,\CS^{B}_1\,,\qquad
  \wt{\CS}^{B}_k\,=\,\CS^{B}_k\,,\quad1<k\leq\ell\,,
  \ee
where $\CS^{B}_i$ are given by \eqref{PinWEYLgens1} and
  \be
  \hat{z}\,
  =\,\hat{\epsilon}_1\cdot\ldots\cdot\hat{\epsilon}_{2\ell+1}\,\in\,\Pin_{2\ell+1}\,,
  \ee
is a lift of the central element $z=T^B_0T^B_1\cdots
T^B_{\ell}\in\CZ(SO_{2\ell+1})$. Then \eqref{LIFTSpin} satisfy the
following relations:
 \be\label{SPINWEYLrels}
  (\wt{\CS}^{B}_i)^2\,=\,(-1)^{\ell}\,,\qquad(\wt{\CS}^{B}_k)^2\,=\,-1\,,\quad1<k\leq\ell\,,\\
  \wt{\CS}^B_i\wt{\CS}^B_j\,=\,\wt{\CS}^B_j\wt{\CS}^B_i\,,\quad a_{ij}=0\,,\\
  \wt{\CS}^{B}_i\wt{\CS}^{B}_j\wt{\CS}^{B}_i\,
  =\,\wt{\CS}^{B}_i\wt{\CS}^{B}_i\wt{\CS}^{B}_j\,,\quad a_{ij}=-1\,,\\
  \wt{\CS}^{B}_1\wt{\CS}^{B}_2\wt{\CS}^{B}_1\wt{\CS}^{B}_2\,
  =\,\wt{\CS}^{B}_2\wt{\CS}^{B}_1\wt{\CS}^{B}_2\wt{\CS}^{B}_1\,.
 \ee
Therefore, \eqref{LIFTSpin} generate a central extension of Weyl
group $W(B_{\ell})$ by $\CZ(\Spin_{2\ell+1})=\mu_2$:
 \be\xymatrix{
  1\ar[r]& \CZ(\Spin_{2\ell+1})\ar[r]& \widetilde{W}(B_{\ell})\ar[r]&
  W(B_{\ell})\ar[r]&1}\,.
 \ee
\end{lem}
\proof: By \eqref{TitsExt} one has
$\hat{z}\CS^B_i\,=\,\CS^B_i\hat{z},\,\forall i\in I$, therefore
using \eqref{PinWEYLBrels} one finds out,
 \be
  (\wt{\CS}^{B}_1)^2\,
  =\,\hat{z}^2(\CS^B_1)^2\,=\,\hat{z}^2\,=\,(-1)^{\ell}\,.
 \ee
The same argument implies the remaining relations in
\eqref{SPINWEYLrels}. $\Box$

The case of even orthogonal group is already covered by the
Proposition \ref{LIFTPin} as the elements $\CS_i^D$ entering
\eqref{PinWEYLDrels} are already in $\Spin_{2\ell}$.

\section{Adjoint action of the Tits groups $W^T_G$}

In this section we compute the action of the elements $S^G_i\in
W^T_G$ in the corresponding Lie algebras $\fg=\Lie(G)$ for classical
groups Lie groups  $G$, including $G=GL_{\ell+1}$. Given a classical group $G$,
the calculation of adjoint action can be done using the appropriate
faithful representation. We provide explicit description of the action using the 
description of the groups $W^T_G$ from Sections 4 and 5 above.  Note
that the resulting formulas easily follow from the explicit
expressions for $\Ad_{\dot{s}_i},\,i\in I$ obtained in \cite{GLO}.
\begin{prop} The adjoint action of the Tits group $W_G^T$ on the Lie
  algebra   $\Fg=\Lie(G)$ via
  homomorphism \eqref{Tits} is given by
  \be\label{ADT1}
  \dot{s}_i\,e_i\,\dot{s}_i^{-1}\,=\,-f_i, \qquad
  \dot{s}_i\,f_i\,\dot{s}_i^{-1}\,=\,-e_i\, ,
  \ee
  \be \label{ADT11}
  \dot{s}_i\,e_j\,\dot{s}_i^{-1}=e_j,
  \qquad \dot{s}_i\,f_j\,\dot{s}_i^{-1}\,=\,f_j,\qquad a_{ij}=0\,,
  \ee
  \be \label{ADT2}
   \dot{s}_i\,e_j\,\dot{s}_i^{-1}\,
   =\,\frac{1}{|a_{ij}|!}
   \underbrace{\bigl[e_i,[\ldots[e_i}_{|a_{ij}|},\,e_j]\ldots]\bigr]\,,\\
   \dot{s}_i\,f_j\,\dot{s}_i^{-1}\,
   =\,\frac{(-1)^{|a_{ij}|}}{|a_{ij}|!}\underbrace{\bigl[f_i,[\ldots[f_i}_{|a_{ij}|},\,f_j]
   \ldots]\bigr],\qquad i\neq j\,.
    \ee
  \end{prop}

Let us emphasize that we keep the ordering $I=\{1,\ldots,\ell\}$ of
the corresponding set of vertices of Dynkin diagram for each
classical $G$, and use the explicit realization of the Cartan-Weyl
generators $\phi(e_i),\,\phi(f_i),\,i\in I$ via the standard
faithful representation $\phi$.

\begin{prop}[\bf The case $A_{\ell}$] The adjoint action of elements $S_i\in
W^T_{GL_{\ell+1}}$ for $i\in\{1,\ldots,\ell\}$ is given by
 \be
  \Ad_{S_i}(e_j)\,
  =\,\left\{\begin{array}{lc}
  f_j\,, & i=j\\
  e_j\,, & |i-j|>1\\
  \ad_{e_i}(e_j)\,, & i-j=-1\\
  -\ad_{e_i}(e_j)\,, & i-j=1
  \end{array}\right.\,,
 \ee
and
 \be
  \Ad_{S_i}(f_j)\,
  =\,\left\{\begin{array}{lc}
  e_j\,, & i=j\\
  f_j\,, & |i-j|>1\\
  -\ad_{f_i}(f_j)\,, & i-j=-1\\
  \ad_{f_i}(f_j)\,, & i-j=1
  \end{array}\right..
 \ee
\end{prop}

\begin{prop}[\bf The case $B_{\ell}$] The adjoint action of elements $S^B_i\in
W^T_{SO_{2\ell+1}}$ for $i\in\{1,\ldots,\ell\}$ is given by
 \be
  \Ad_{S^B_i}(e_j)\,
  =\,\left\{\begin{array}{lc}
  f_j\,, & i=j\\
  e_j\,, & |i-j|>1\\
  -\frac{(-1)^{\delta_{i,1}}}{|a_{ij}|!}\ad^{-a_{ij}}_{e_i}(e_j)\,, & i-j=-1\\
  \ad_{e_i}(e_j)\,, & i-j=1
  \end{array}\right.\,,
 \ee
and
 \be
  \Ad_{S^B_i}(f_j)\,
  =\,\left\{\begin{array}{lc}
  e_j\,, & i=j\\
  f_j\,, & |i-j|>1\\
  \frac{1}{|a_{ij}|!}\ad^{-a_{ij}}_{f_i}(f_j)\,, & i-j=-1\\
  -\ad_{f_i}(f_j)\,, & i-j=1
  \end{array}\right..
 \ee
\end{prop}

\begin{prop}[\bf The case $C_{\ell}$] The adjoint action of elements $S^C_i\in
W^T_{\Sp_{2\ell}}$ for $i\in\{1,\ldots,\ell\}$ is given by
 \be
  \Ad_{S^C_i}(e_j)\,
  =\,\left\{\begin{array}{lc}
  (-1)^{\delta_{i,1}}f_j\,, & i=j\\
  e_j\,, & |i-j|>1\\
  -\ad_{e_i}(e_j)\,, & i-j=-1\\
  \frac{1}{|a_{ij}|!}\ad^{-a_{ij}}_{e_i}(e_j)\,, & i-j=1
  \end{array}\right.,
 \ee
and
 \be
  \Ad_{S^C_i}(f_j)\,
  =\,\left\{\begin{array}{lc}
  (-1)^{\delta_{i,1}}e_j\,, & i=j\\
  f_j\,, & |i-j|>1\\
  \ad^{-a_{ij}}_{f_i}(f_j)\,, & i-j=-1\\
  \frac{(-1)^{|a_{ij}|}}{|a_{ij}|!}\ad^{-a_{ij}}_{f_i}(f_j)\,, & i-j=1
  \end{array}\right..
 \ee
\end{prop}

\begin{prop}[\bf The case $D_{\ell}$] The adjoint action of elements $S^D_i\in
W^T_{O_{2\ell}}$ for $i\in\{1,2;\,3,\ldots,\ell\}$ is given by
 \be
  \Ad_{S^D_i}(e_j)\,
  =\,\left\{\begin{array}{lc}
  f_j\,, & i=j\\
  -e_j\,, & a_{ij}=0\,,\quad\iota(i)=j\\
  e_j\,, & a_{ij}=0\\,,\quad\iota(i)\neq j\
  -\ad_{e_i}(e_j)\,, & a_{ij}=-1\,,\quad i<j\\
  \ad_{e_i}(e_j)\,, & a_{ij}=-1\,,\quad i>j  
 \end{array}\right.
 \ee
and
 \be
  \Ad_{S^D_i}(f_j)\,
  =\,\left\{\begin{array}{lc}
  e_j\,, & i=j\\
  -f_j\,, & \iota(i)=j\,,\quad\iota(i)=j\\
  f_j\,, & a_{ij}=0\,,\quad\iota(i)\neq j\\
  \ad_{f_i}(f_j)\,, & a_{ij}=-1\,,\quad i<j\\
  -\ad_{f_i}(f_j)\,, & a_{ij}=-1\,,\quad i>j  
  \end{array}\right..
 \ee
\end{prop}

\section{Appendix:  General linear and classical Lie groups}\label{APP}

\subsection{General linear group $GL_{\ell+1}(\IC)$}\label{APPA}

Let $\gl_{\ell+1}(\IC)=\gl(V)$ be the Lie algebra induced by  endomorphisms
$\phi:\,V\to V$ of the complex vector space $V\simeq\IC^{\ell+1}$
via the  commutator Lie bracket:
 \be
  [\phi_1,\,\phi_2]\,=\,\phi_1\circ\phi_2\,-\,\phi_2\circ\phi_1\,,
 \ee
where $\circ$ is a composition of linear maps $V\to V$. Let
$\fh\subset\gl_{\ell+1}(\IC)$ be its maximal commutative subalgebra.
Choosing a basis $\{\ve_1,\ldots,\ve_{\ell+1}\}\subset V$ of vector
space $V$ we identify $\gl_{\ell+1}(\IC)$ with the Lie
algebra of matrices with the basis $e_{ij},\,1\leq i,j\leq\ell+1$
defined by the following endomorphisms:
 \be\label{Matunits}
  e_{ij}\,:\quad
  V\,\longrightarrow\,V\,,\qquad\ve_k\,\longmapsto\,\delta_{ik}\ve_j\,,\quad
  1\leq k\leq\ell+1\,.
 \ee
The corresponding Lie brackets is given by
 \be
  \bigl[e_{ij},\,e_{kl}\bigr]\,=\,\delta_{jk}e_{il}\,-\,\delta_{il}e_{kj}\,.
 \ee
The Cartan subalgebra $\fh_{\ell+1}\subset\gl_{\ell+1}(\IC)$ is
identified with a subalgebra of diagonal matrices.

Fix an orthonormal basis $\{\e_1,\ldots,\e_{\ell+1}\}\subset\fh^*$
in the Euclidean space
$\bigl(\fh^*\simeq\IC^{\ell+1};\,\<\,,\,\>\bigr)$ and identify $\fh$
and $\fh^*$ via the standard scalar product $\<\,,\>$ on
$\IC^{\ell+1}$. Thus we can consider the diagonal entries
$\e_k:=e_{kk}$ as functions on $\fh$, so that they provide a
(complex) linear coordinate system on $\fh$. Then any linear
functional $\la\in\fh^*$ has the form
 \be
  \la\,=\,\la_1\e_1\,+\,\ldots\,+\,\la_{\ell+1}\e_{\ell+1}\,.
 \ee

The general linear group $GL_{\ell+1}(\IC)$ is isomorphic to the
group of invertible rank $\ell+1$ matrices, and maximal torus
$H_{\ell+1}\subset GL_{\ell+1}(\IC)$ is identified with the
subgroup of diagonal matrices. The diagonal entries of
$t_k=t_k(g),\,g\in H_{\ell+1}$ provide coordinates on the
maximal torus $H_{\ell+1}$ and generate the group of (rational)
characters
 \be
  X^*(H_{GL_{\ell+1}})\,
  =\,\Hom(H_{\ell+1},\,\IC^*)\simeq\fh^*_{\IZ}=\IZ^{\ell+1}\,.
 \ee
Namely, any $\la\in\fh^*$ gives rise to a homomorphism
 \be\label{WeightGL}
  e^{\la}\,:\quad H_{\ell+1}\,\longrightarrow\,\IC^*\,,\qquad
  e^h\,\longmapsto\,e^{2\pi\la(h)}\,=\,e^{2\pi\<\la,\,h\>}\,,
  \quad h\in\fh\,;\\
  e^{\la}\,=\,t_1^{\la_1}t_2^{\la_2}\cdot\ldots\cdot t_{\ell+1}^{\la_{\ell+1}}\,.
 \ee
The elements $\la\in\fh^*_{\IZ}$ consisting of $\la\in\fh^*$ with
integral components $(\la_1,\ldots,\la_{\ell+1})\in\IZ^{\ell+1}$ are
called the weights; they span the weight lattice:
 \be\label{GLweights}
  \La_W\,=\,\IZ\e_1\oplus\ldots\oplus\IZ\e_{\ell+1}\simeq\,\fh^*_{\IZ}\,.
 \ee

The adjoint action of $H_{\ell+1}$ in Lie algebra
$\gl_{\ell+1}(\IC)$ provides the Cartan decomposition of
$\gl_{\ell+1}(\IC)$ with respect to $H_{\ell+1}$
 \be
  \gl_{\ell+1}(\IC)\,
  =\,\fh\,\oplus\,
  \bigoplus_{\a_{ij}\in\Phi(A_{\ell})}\bigl(\gl_{\ell+1}\bigr)_{ij}\,,\\
  (\gl_{\ell+1})_{ij}\,=\,\bigl\{X\in\gl_{\ell+1}\,:
  \quad\ad_h(X)\,=\,\a_{ij}(h)X,\,\,\forall h\in\fh\bigr\}\,
  =\,\IC e_{ij}\,,\quad i\neq j\,,
 \ee
where the corresponding root system $\Phi(A_{\ell})$ of type
$A_{\ell}$ reads
 \be\label{RootdataGL}
  \Phi(A_{\ell})\,=\,\bigl\{\a_{ij}\,=\,\e_i-\e_j\,,\quad i\neq
  j\bigr\}\,\subset\,\fh^*_{\IZ}\,.
 \ee
The simple root system $\Pi_{A_{\ell}}\subset\Phi_+(A_{\ell})$ is
given by
 \be
 \Pi_{A_{\ell}}\,=\,\bigl\{\a_i\,=\,\e_i-\e_{i+1}\bigr\}\,
 \subset\,\Phi_+(A_{\ell})\,,\\
  \Phi_+(A_{\ell})\,=\,\bigl\{\a_{ij}\,=\,\e_i-\e_j\,,\quad
  i<j\bigr\}\,\subset\,\fh^*_{\IZ}\,,
 \ee
and the root lattice is
 \be\label{GLroots}
  \La_R\,=\,\IZ\a_1\oplus\ldots\oplus\IZ\a_{\ell}\,\subset\,\La_W\simeq\fh^*_{\IZ}\,.
 \ee

  Let
$X_*(H_{\ell+1})=\Hom(\IC^*,H_{\ell+1})$ be the group
of co-characters (i.e. one-parametric subgroups) in the maximal
torus:
 \be\label{GLcochars}
  \chi_a^{\vee}\,:\quad\IC^*\,\longrightarrow\,H_{\ell+1}\,,\qquad
  t\,\longmapsto\,\diag\{t^{a_1},\ldots,t^{a_{\ell+1}}\}\,.
 \ee
 The group of co-characters  $X_*(H_{\ell+1})\simeq\fh_{\IZ}=\IZ^{\ell+1}$
is dual to $X^*(H_{\ell+1})$ via
 \be
  e^{\la}\bigl(\chi_a^{\vee}(t)\bigr)\,
  =\,t^{a_1\la_1+\ldots+a_{\ell+1}\la_{\ell+1}}\,
  =\,t^{\<a^{\vee},\,\la\>}\,.
 \ee
The dual basis in $\fh_{\IZ}$ is referred to as \emph{fundamental
co-weights}:
 \be\label{GLcoweights}
  p_k^{\vee}\,=\,e_{11}+\ldots+e_{kk}\,,\qquad
  \<p_k^{\vee},\,\a_i\>\,=\,\delta_{ki}\,,\qquad1\leq
  k,i\leq\ell+1\,,
 \ee
that span the co-weight lattice:
 \be
  \La_W^{\vee}\,
  =\,\bigl\{h\in\fh\,:\,\la(h)\in\IZ\,,\quad\forall\la\in\La_W\bigr\}\,
  =\,\IZ p_1^{\vee}\oplus\ldots\oplus\IZ p_{\ell+1}^{\vee}\,.
 \ee
The lattice $\La_W^{\vee}$ can be identified with the kernel of the
map $\exp:\fh\to H_{\ell+1}$, so that
$H_{\ell+1}=\fh/\La_W^{\vee}$\,.

The determinant map $\det: GL_{\ell+1}(\IC)\to \IC^*$ may be defined
  as a unique homomorphims such that
  \be
   \det\,:\quad e^{\sum\limits_{i=1}^{\ell+1}t_ie_{ii}}\longrightarrow
     e^{\sum\limits_{i=1}^{\ell+1}t_i}.
  \ee
We have  the following exact sequence
 \be\xymatrix{
  1\ar[r] & SL_{\ell+1}(\IC)\ar[r]^{\pi} &
  GL_{\ell+1}(\IC)\ar[r]^{\det} &  \IC^*\ar[r] & 1}\,,
 \ee
where $SL_{\ell+1}(\IC)$ is the unimodular subgroup. Its Lie algebra
$\ssl_{\ell+1}$ has the standard faithful representation given by
the following homomorphism of associative algebras:
 \be\label{StandardrepA}
  \phi\,:\quad\CU\bigl(\ssl_{\ell+1}(\IC)\bigr)\,\longrightarrow\,
  \End(\IC^{\ell+1})\,,\\
  \phi(X_{\a_i})\,=\,e_{i,i+1}\,,\qquad
  \phi(X_{-\a_i})\,=\,-e_{i+1,i}\,,\qquad
  \phi(h_{\a_i})\,=\,e_{i,i}-e_{i+1,i+1}\,,
 \ee
where $e_{ij}$ are the matrix units \eqref{Matunits}. The Lie algebra
$\ssl_{\ell+1}$ may  be identified with the Lie subalgebra of
matrices with zero trace in $\gl_{\ell+1}$, so that
 \be
  \gl_{\ell+1}\,=\,\ssl_{\ell+1}\oplus\IC\,.
 \ee
The Lie algebras $\ssl_{\ell+1}$ and $\gl_{\ell+1}$ share the Dynkin
diagram $\Gamma_{A_{\ell}}$ with the set of vertices $I$ of size
$\ell$, the reduced root system $\Phi(A_{\ell})$ and the Weyl group
$W_{A_{\ell}}$. However, in case of $\ssl_{\ell+1}$ its rank equals
to the rank of the root lattice, contrary to the $\gl_{\ell+1}$
case. More precisely, let $\{\ve_1,\ldots,\ve_{\ell+1}\}$ be an
orthonormal basis in $\IC^{\ell+1}$. Then the
fundamental weights $p_k$ of $\gl_{\ell+1}$ are linear forms defined
by $p_k(h_{\a_i})=\delta_{ki}$ (see \eqref{GLweights},
\eqref{GLcoweights}):
 \be
  p_k\,=\,\e_1+\,\ldots+\,\e_k\,,\qquad1\leq k\leq\ell+1\,.
 \ee
Define a collection of vectors in $\IR^{\ell+1}$,
  \be
    \e_k^{A_{\ell}}\,=\,\e_k\,-\,\frac{\e_1+\ldots+\e_{\ell+1}}{\ell+1}\,;\qquad
    1\leq k\leq\ell+1\,,
  \ee
which  span  a codimension one Euclidean subspace, due to
$\e_1^{A_{\ell}}+\ldots+\e^{A_{\ell}}_{\ell+1}\,=\,0$. Then the
simple roots and fundamental weights of $\ssl_{\ell+1}$ can be
written as follows:
 \be\label{RootDataA}
  \a_k^{A_{\ell}}\,=\,\e_k^{A_{\ell}}-\e_{k+1}^{A_{\ell}}\,,\qquad
  \vp_k^{A_{\ell}}\,=\,p_k\,
  -\,\frac{k}{\ell+1}\,p_{\ell+1}\,
  =\,\e_1^{A_{\ell}}+\ldots+\e_k^{A_{\ell}}\,.
  \ee
The Cartan matrix $A\,=\,\|a^{A_{\ell}}_{ij}\|$ with
$a_{ij}=\<\a^{A_{\ell}}_j,(\a_i^{\vee})^{A_{\ell}}\>$, and its
inverse $A^{-1}=\|c_{ij}^{A_{\ell}}\|$ take the form:
 \be\label{ACartanMat}
  A\,=\,\left(\begin{smallmatrix}
  2&-1&0&\ldots&0\\-1&2&\ddots&\ddots&\vdots\\
  0&\ddots&\ddots&-1&0\\
  \vdots&\ddots&-1&2&-1\\
  0&\ldots&0&-1&2
  \end{smallmatrix}\right)\,,\quad
  A^{-1}\,
  =\,\frac{1}{\ell+1}\left(\begin{smallmatrix}
  \ell & \ell-1 & \ell-2 & \ldots & 2 & 1\\
  \\
  \ell-1 & 2(\ell-1) & 2(\ell-2) & \ldots & 2\cdot2 & 2\\
  \\
  \ell-2 & 2(\ell-2) & 3(\ell-2) & \ldots & 3\cdot2 & 3\\
  \\
  \ldots & \ldots & \ldots & \ldots & \ldots & \ldots\\
  \\
  2& 2\cdot2 & 3\cdot2 & \ldots & (\ell-1)2 & \ell-1\\
  \\
  1 & 2 & 3 & \ldots & \ell-1 & \ell
 \end{smallmatrix}\right)\,.
 \ee
In these terms one obtains the following expressions for the coroots
and coweights:
 \be\label{Acoweights}
  h^{A_{\ell}}_i\,=\,(\a_i^{A_{\ell}})^{\vee}\,
  =\,\sum_{j\in I}\<\a_j^{A_{\ell}},\,(\a_i^{A_{\ell}})^{\vee}\>(\vp_j^{A_{\ell}})^{\vee}\,
  =\,\sum_{j\in I}a_{ij}^{A_{\ell}}(\vp_j^{A_{\ell}})^{\vee}\,,\\
  (\vp_i^{A_{\ell}})^{\vee}\,=\,\sum_{j\in I}c_{ij}^{A_{\ell}}(\a_j^{A_{\ell}})^{\vee}\,,\qquad i\in I\,.
 \ee

\subsection{The case of type $B_{\ell}$ root system}\label{APPB}

The Lie algebra $\so_{2\ell+1}$ can be identified with the
$\theta_B$-fixed subalgebra:
 \be\label{Bfaithful}
  \so_{2\ell+1}\,=\,\bigl\{X\in\gl_{2\ell+1}(\IC)\,:\quad
  X\,=\,\theta_B(X)\bigr\}\,,\qquad\theta_B(X)\,=\,-S_BX^{\tau}S_B^{-1}\,,\\
  S_B\,=\,\diag\bigl(1,-1,\,\ldots,\,1\bigr)\,\in\,GL_{2\ell+1}\,,
 \ee
where $\tau$ is the transposition along the opposite diagonal
\eqref{TranspOpp}. This provides the faithful representation,
 \be\label{StandardrepB}
  \phi\,:\quad\so_{2\ell+1}\,\longrightarrow\,\gl_{2\ell+1}(\IC)\,,
 \ee
given by the following presentation of the Chevalley-Weyl generators
$h_k=\a_k^{\vee},\,e_k,\,f_k,\,k\in I$ of $\so_{2\ell+1}$ in terms
of the generators of $\gl_{2\ell+1}(\IC)$:
 \be\label{BstandardRep}
  \phi(e_1)\,=\,\sqrt{2}(e_{\ell,\ell+1}+e_{\ell+1,\ell+2})\,,\qquad
  \phi(f_1)\,=\,\sqrt{2}(e_{\ell+1,\ell}+e_{\ell+2,\ell+1})\,;\\
  \phi(e_k)\,=\,e_{\ell+1-k,\ell+2-k}+e_{\ell+k,\ell+1+k}\,,\quad
  \phi(f_k)\,=\,e_{\ell+2-k,\ell+1-k}+e_{\ell+1+k,\ell+k}\,,\\1<k\leq\ell\,;\\
  \phi(h_1)\,=\,2(e_{\ell\ell}-e_{\ell+2,\ell+2})\,,\\
  \phi(h_k)\,=\,(e_{\ell+1-k,\ell+1-k}+e_{\ell+k,\ell+k})\,
  -\,(e_{\ell+2-k,\ell+2-k}+e_{\ell+1+k,\ell+1+k})\,,\\1<k\leq\ell\,.
 \ee
The representation \eqref{BstandardRep} implies the following
presentation of the type $B_{\ell}$ root system in terms of the root
system of type $A_{2\ell}$
 \be\label{BAroots}
  \a_i^{B_{\ell}}\,
  =\,\a_{\ell+1-i}^{A_{2\ell}}\,+\,\a_{\ell+i}^{A_{2\ell}}\,,\qquad i\in
  I\,.
 \ee
Explicitly the symmetry of the Dynkin diagram of type $A_{2\ell}$ is
given by
 \be
  \iota\,:\quad i\,\longmapsto\,2\ell+1-i\,,\qquad
  i\in\{1,2,\ldots,2\ell\}\,,\\
  \a_i^{A_{2\ell}}\,\longmapsto\,\a_{2\ell+1-i}^{A_{2\ell}}\,.
 \ee
Then  \eqref{BAroots} reads
 \be
  \a_i^{B_{\ell}}\,=\,\a_{\ell+1-i}^{A_{2\ell}}\,+\,\iota\bigl(\a_{\ell+1-i}^{A_{2\ell}}\bigr)\,.
 \ee
Introduce the Euclidean space $\IR^{2\ell+1}$ with a standard
orthonormal basis $\bigl\{\e_1,\ldots,\e_{2\ell+1}\bigr\}$, acted on
by the involution $\iota$ as follows:
 \be
  \iota\,:\quad\IR^{2\ell+1}\,\longrightarrow\,\IR^{2\ell+1}\,,\qquad
  \e_i\,\longmapsto\,-\e_{2\ell+2-i}\,.
 \ee
Consider the $\iota$-fixed Euclidean subspace
$\IR^{\ell}\subset\IR^{2\ell+1}$ spanned by
 \be\label{Bfolding}
  \e_k^{B_{\ell}}\,=\,\e_{\ell+1-k}\,+\,\iota\bigl(\e_{\ell+1-k}\bigr)\,
  =\,\e_{\ell+1-k}\,-\,\e_{\ell+1+k}\,,\qquad1\leq
  k\leq\ell\,,
 \ee
then the root data of type $B_{\ell}$ is given by
 \be\label{RootDataB}
  \left\{\begin{array}{l}
  \a_1^{B_{\ell}}\,=\,\e_1^{B_{\ell}}\\
  \a_k^{B_{\ell}}\,=\,\e_k^{B_{\ell}}-\e_{k-1}^{B_{\ell}}\,,\\
  1<k\leq\ell\,;
  \end{array}\right.\qquad
  \left\{\begin{array}{l}
  \vp_1^{B_{\ell}}\,=\,\frac{\e_1^{B_{\ell}}+\ldots+\e_{\ell}^{B_{\ell}}}{2}\\
  \vp_k^{B_{\ell}}\,=\,\e_k^{B_{\ell}}+\ldots+\e_{\ell}^{B_{\ell}}\,,\\
  1<k\leq\ell\,;
  \end{array}\right.\\
\Phi(B_{\ell})\,=\,
\bigl\{\pm\e_i^{B_{\ell}}\pm\e_j^{B_{\ell}}\,,\quad\pm\e_i^{B_{\ell}}\bigr\}\,.
 \ee
The simple co-roots and fundamental co-weights are determined via
 \be
  \<\a_i^{B_{\ell}},\,(\vp_j^{B_{\ell}})^{\vee}\>\,
  =\,\<(\a_i^{B_{\ell}})^{\vee},\,\vp_j^{B_{\ell}}\>\,=\,\delta_{ij}\,,
 \ee
and have the following form:
 \be
  \left\{\begin{array}{l}
  \bigl(\a_1^{B_{\ell}}\bigr)^{\vee}\,=\,2\e_1^{B_{\ell}}\\
  \bigl(\a_k^{B_{\ell}}\bigr)^{\vee}\,=\,\e_k^{B_{\ell}}-\e_{k-1}^{B_{\ell}}\,,\\
  1<k\leq\ell\,;
  \end{array}\right.\qquad
  \left\{\begin{array}{l}
           \bigl(\vp_i^{B_{\ell}}\bigr)^{\vee}\,=\,
           \e_i^{B_{\ell}}+\ldots+\e_{\ell}^{B_{\ell}}\,,\\
  i\in I
  \end{array}\right.\\
\Phi^{\vee}(B_{\ell})\,=\,\bigl\{\pm\e_i^{B_{\ell}}\pm\e_j^{B_{\ell}}\,,
\quad\pm2\e_i^{B_{\ell}}\bigr\}\,.
 \ee
The Cartan matrix
$A=\|a_{ij}^{B_{\ell}}\|=\|\<\a_j^{B_{\ell}},\,\bigl(\a_i^{B_{\ell}}\bigr)^{\vee}\>\|$
and its inverse $A^{-1}=\|c_{ij}\|$ are given by
 \be\label{BCartanMat}
  A\,
  =\,\left(\begin{smallmatrix}
  2&-2&0&\ldots&0\\-1&2&-1&\ddots&\vdots\\
  0&\ddots&\ddots&\ddots&0\\
  \vdots&\ddots&-1&2&-1\\
  0&\ldots&0&-1&2
  \end{smallmatrix}\right)\,,\qquad
  A^{-1}\,
  =\,\left(\begin{smallmatrix}
  \frac{\ell}{2}&\ell-1&\ldots&3&2&1\\
  \frac{\ell-1}{2}&\ell-1&\ldots&3&2&1\\
  \vdots&\vdots&\ddots&\vdots&\vdots\\
  \frac{3}{2}&3&\ldots&3&2&1\\
  1&2&\ldots&2&2&1\\\frac{1}{2}&1&\ldots&1&1&1
  \end{smallmatrix}\right)\,,\\
  h_i^{B_{\ell}}\,=\,\bigl(\a_i^{B_{\ell}}\bigr)^{\vee}\,
  =\,\sum_{j\in
    I}\<\a_j^{B_{\ell}},\,\big(\a_i^{B_{\ell}}\bigr)^{\vee}\>
  \bigl(\vp_j^{B_{\ell}}\bigr)^{\vee}\,
  =\,\sum_{j\in
  I}a_{ij}^{B_{\ell}}\bigl(\vp_j^{B_{\ell}}\bigr)^{\vee}\,,\\
\bigl(\vp_i^{B_{\ell}}\bigr)^{\vee}\,=\, \sum_{j\in
I}c_{ij}\bigl(\a_j^{B_{\ell}}\bigr)^{\vee}\,,\qquad i\in I\,.
 \ee

\begin{lem}\label{BAutroots}
The Weyl group $W(B_{\ell})$ is isomorphic to a subgroup in
$W(A_{2\ell})$ via
 \be\label{BWeyl}
  s_1^{B_{\ell}}\,=\,s_{\ell}^{A_{2\ell}}s_{\ell+1}^{A_{2\ell}}s_{\ell}^{A_{2\ell}}\,
  =\,s_{\ell+1}^{A_{2\ell}}s_{\ell}^{A_{2\ell}}s_{\ell+1}^{A_{2\ell}}\,,\\
  s_k^{B_{\ell}}\,=\,s_{\ell+1-k}^{A_{2\ell}}s_{\ell+k}^{A_{2\ell}}\,
  =\,s_{\ell+k}^{A_{2\ell}}s_{\ell+1-k}^{A_{2\ell}}\,,\quad1<k\leq\ell\,.
 \ee
\end{lem}
\proof: For $1<k\leq\ell$ given
$\la=\la_1\e^{B_{\ell}}_1+\ldots+\la_{\ell}\e^{B_{\ell}}_{\ell}$ by
straightforward computation one has
 \be
  s_k^{B_{\ell}}(\la)\,=\,\la\,
  -\,\<\la,\,(\a_k^{B_{\ell}})^{\vee}\>_{B_{\ell}}\a_k^{B_{\ell}}\,
  =\,\la\,
  -\,\<\la,\,\e_k^{B_{\ell}}-\e_{k-1}^{B_{\ell}}\>_{B_{\ell}}
  (\e_k^{B_{\ell}}-\e_{k-1}^{B_{\ell}})\\
  =\,\la\,
  -\,(\la_k-\la_{k-1})(\e_k^{B_{\ell}}-\e_{k-1}^{B_{\ell}})\\
  =\,\la_1\e^{B_{\ell}}_1+\ldots+\la_k\e^{B_{\ell}}_{k-1}
  +\la_{k-1}\e^{B_{\ell}}_k+\ldots+\la_{\ell}\e^{B_{\ell}}_{\ell}\\
  =\,\la_1\e^{A_{2\ell}}_{\ell}+\ldots+\la_k\e^{A_{2\ell}}_{\ell+2-k}+
  \la_{k-1}\e^{A_{2\ell}}_{\ell+1-k}
  +\ldots+\la_{\ell}\e^{A_{2\ell}}_1\\
  -\,\la_1\e^{A_{2\ell}}_{\ell+2}-\ldots-\la_k\e^{A_{2\ell}}_{\ell+k}-
  \la_{k-1}\e^{A_{2\ell}}_{\ell+1+k}-\ldots
  -\la_{\ell}\e^{A_{2\ell}}_{2\ell+1}\\
  =\,\la\,-\,\<\la,\,(\a_{\ell+1-k}^{A_{2\ell}})^{\vee}\>_{A_{2\ell}}
  \a_{\ell+1-k}^{A_{2\ell}}\,
  -\,\<\la,\,(\a_{\ell+k}^{A_{2\ell}})^{\vee}\>_{A_{2\ell}}\a_{\ell+k}^{A_{2\ell}}\\
  =\,s_{\ell+k}^{A_{2\ell}}s_{\ell+1-k}^{A_{2\ell}}(\la)\,.
 \ee
Similarly, for $k=1$ one has
 \be
  s_1^{B_{\ell}}(\la)\,=\,\la\,
  -\,\<\la,\,(\a_1^{B_{\ell}})^{\vee}\>_{B_{\ell}}\a_1^{B_{\ell}}\,
  =\,\la\,
  -\,\<\la,\,2\e_1^{B_{\ell}}\>_{B_{\ell}}\e_1^{B_{\ell}}\\
  =\,-\la_1\e^{B_{\ell}}_1+\la_2\e^{B_{\ell}}_2+\ldots+\la_{\ell}\e^{B_{\ell}}_{\ell}\\
  =\,\la_{\ell}\e^{A_{2\ell}}_1+\ldots+\la_2\e^{A_{2\ell}}_{\ell-1}-
  \la_1\e^{A_{2\ell}}_{\ell}\,
  +\,\la_1\e^{A_{2\ell}}_{\ell+2}-\la_2\e^{A_{2\ell}}_{\ell+3}-\ldots-
  \la_{\ell}\e^{A_{2\ell}}_{2\ell+1}\\
  =\,s_{\ell}^{A_{2\ell}}s_{\ell+1}^{A_{2\ell}}s_{\ell}^{A_{2\ell}}(\la)\,,
 \ee
since
$s_{\ell}^{A_{2\ell}}s_{\ell+1}^{A_{2\ell}}s_{\ell}^{A_{2\ell}}$
simply swaps
$\e_{\ell}^{A_{2\ell}}\leftrightarrow\e_{\ell+2}^{A_{2\ell}}$.
$\Box$

The Lie group $O_{2\ell+1}$ may be presented as $\theta_B$-invariant
subgroup of $GL_{2\ell+1}$:
 \be\label{FaithulRepB}
  O_{2\ell+1}\,=\,\bigl\{g\in GL_{2\ell+1}\,:\quad
  g\,=\,g^{\theta_B}\,\bigr\}\,\,.
 \ee
Let us describe the corresponding Tits group $W^T_{O_{2\ell+1}}$ in
terms of generators of $W^T_{GL_{2\ell+1}}$.
\begin{lem}\label{BTitsLift}
Let $\dot{s}^B_i=e^{\frac{\pi}{2}J_i},\,i\in I$ be the Tits
generators \eqref{Titsgen} then the following holds:
\begin{enumerate}
\item The Tits group $W^T_{O_{\ell+1}}$ is isomorphic to a
subgroup of $W^T_{GL_{2\ell+1}}$ via
 \be\label{TitsgenB}
  \dot{s}_1^B\,
  =\,\dot{s}_{\ell}^{A_{2\ell}}\dot{s}_{\ell+1}^{A_{2\ell}}\dot{s}_{\ell}^{A_{2\ell}}\,,\qquad
  \dot{s}^B_k\,
  =\,\dot{s}_{\ell+1-k}^{A_{2\ell}}\dot{s}_{\ell+k}^{A_{2\ell}}\,,\quad1<k\leq\ell\,.
 \ee

\item The elements \eqref{TitsgenB} belong to the $\theta_B$-fixed subgroup
$(W^T_{GL_{2\ell+1}})^{\theta_B}$.

\item Presentation \eqref{TitsgenB} matches with \eqref{BWeyl} due
to
 \be
  \Ad_{\dot{s}^B_i}|_{\fh}\,=\,s^B_i\,.
 \ee
\end{enumerate}
\end{lem}

\proof: (1) Since $\dot{s}^B_i\in SO_{2\ell+1},\,i\in I$, one might
verify \eqref{TitsgenB} using the standard faithful representation
$\phi:\,SO_{2\ell+1}\to GL_{2\ell+1}$. For $i=1$, we have
 \be\label{TitsgenB1}
  \phi(\dot{s}^B_1)\,
  =\,\phi(e^{\frac{\pi}{2}J_1})\,
  =\left(\begin{smallmatrix}
  \Id_{\ell-1}&&&&\\
  &0&0&1&\\&0&-1&0&\\&1&0&0&\\
  &&&&\Id_{\ell-1}
  \end{smallmatrix}\right)\,
  =\,\phi(\dot{s}_{\ell}^{A_{2\ell}}\dot{s}_{\ell+1}^{A_{2\ell}}\dot{s}_{\ell}^{A_{2\ell}})\,,
 \ee
and similarly, for $1<k\leq\ell$ we obtain
 \be\label{TitsgenB2}
  \phi(\dot{s}^B_k)\,
  =\,\phi(e^{\frac{\pi}{2}J_k})\,
  =\left(\begin{smallmatrix}
  \Id_{\ell-k}&&&&&&\\
  &0&-1&&&&\\&1&0&&&&\\
  &&&\Id_{2k-3}&&&\\&&&&0&-1&\\
  &&&&1&0&\\&&&&&&\Id_{\ell-k}
  \end{smallmatrix}\right)\,
  =\,\phi(\dot{s}_{\ell+1-k}^{A_{2\ell}}\dot{s}_{\ell+k}^{A_{2\ell}})\,.
 \ee

(2) From \eqref{CINV} one reads
 \be
  \theta_B\,:\qquad\dot{s}^{A_{2\ell}}_i\,\longmapsto\,
  \dot{s}^{A_{2\ell}}_{2\ell+1-i}\,,\qquad1\leq
  i\leq2\ell+1\,,
 \ee
and using the Tits relations \eqref{Tits} one infers that
\eqref{TitsgenB} are $\theta_B$-invariant.

(3) Follows by \eqref{Tits} and \eqref{BWeyl}. $\Box$

By analogy with the case of general linear group there is another
way to lift simple root generators $s^B_i\in W(B_{\ell})$ into the
Tits group $W^T_{O_{2\ell+1}}$. Recall that $W^T_{GL_{2\ell+1}}$
contains the following elements:
 \be\label{TitsGL}
  T_k\,:=\,e^{\pi\imath e_{kk}}\,,\qquad1\leq k\leq (2\ell+1)\,;\\
  S_i\,=\,T_i\dot{s}_i\,=\,\dot{s}_iT_{i+1}\,\,,\qquad
  \overline{S}_i\,=\,T_{i+1}\dot{s}_i\,=\,\dot{s}_iT_i\,\,,\qquad1\leq  i\leq2\ell\,.
 \ee
\begin{lem} The following elements belong to the fixed point
subgroup $(W^T_{GL_{2\ell+1}})^{\theta_B}$:
 \be\label{2torsB}
  S^B_1\,
  =\,S_{\ell+1}S_{\ell}S_{\ell+1}\,,\qquad
  S^B_k\,=\,S_{\ell+1-k}\overline{S}_{\ell+k}\,,\quad1<k\leq\ell\,;\\
  T^B_i\,=\,T_{\ell+1-i}T_{\ell+1+i}\,=\,e^{\pi\imath\e^B_i}\,,\quad1\leq i\leq\ell\qquad
  \text{and}\qquad T^B_0\,:=\,T_{\ell+1}\,,
 \ee
where $T^B_0$ has $\det(T^B_0)=-1$ and generates the center
$\CZ\bigl(W^T_{O_{2\ell+1}}\bigr)=\mu_2$.
\end{lem}
\proof: The explicit action of $\theta_B$ on generators
\eqref{TitsGL} reads \eqref{permatAction}:
 \be
  \theta_B(S_i)\,=\,\overline{S}_{2\ell+1-i}\,,\qquad
  \theta_B(\overline{S}_i)\,=\,S_{2\ell+1-i}\,,\qquad\theta_B(T_k)\,=\,T_{2\ell+2-k}\,.
 \ee
Using \eqref{WeylExt2} this implies that \eqref{2torsB} are
invariant under $\theta_B$.\,\, $\Box$

\begin{cor} The elements \eqref{TitsgenB} and \eqref{2torsB} can be identified via
 \be\label{BTits}
  S^B_1\,=\,T^B_0\,\dot{s}_1^B\,\,,\qquad
  S^B_k\,=\,T^B_k\dot{s}^B_k\,\,,\quad1<k\leq\ell\,.
 \ee
\end{cor}
\proof: By straightforward computation using faithful representation
\eqref{FaithulRepB} one finds
 \be
  \phi(S^B_1)\,
  =\left(\begin{smallmatrix}
  \Id_{\ell-1}&&&&\\
  &0&0&1&\\&0&1&0&\\&1&0&0&\\
  &&&&\Id_{\ell-1}
  \end{smallmatrix}\right)\,,\qquad
  \phi(S^B_k)\,
  =\left(\begin{smallmatrix}
  \Id_{\ell-k}&&&&&&\\
  &0&1&&&&\\&1&0&&&&\\
  &&&\Id_{2k-3}&&&\\&&&&0&-1&\\
  &&&&-1&0&\\&&&&&&\Id_{\ell-k}
  \end{smallmatrix}\right)\,.
 \ee
Identifying this with \eqref{TitsgenB1}, \eqref{TitsgenB2} one
deduces \eqref{BTits}. $\Box$

\begin{prop}\label{BWeylLift}
The elements $S^B_i\in W^T_{O_{2\ell+1}},\,i\in I$ defined in
\eqref{2torsB} satisfy the relations \eqref{BR0}, \eqref{BR} of type
$B_{\ell}$:
 \be\label{BCentExt}
  (S^B_1)^2\,=\,\ldots\,=\,(S^B_{\ell})^2\,
  =\,1\,\,,\\
  S^B_iS^B_j\,=\,S^B_jS^B_i\,,\quad i,j\in I\quad \text{such
  that}\quad a_{ij}=0\,\,;\\
  S^B_iS^B_jS^B_i\,=\,S^B_jS^B_iS^B_j\,,\qquad i,j\in I\quad\text{such
  that}\quad a_{ij}=-1\,,\\
  S^B_1S^B_2S^B_1S^B_2\,=\,S^B_2S^B_1S^B_2S^B_1\,.
 \ee
\end{prop}
\proof: One might prove \eqref{BCentExt} applying faithful
representation \eqref{FaithulRepB}. Alternatively, since the
relations \eqref{BCentExt} are exactly the relations \eqref{BR},
they can be derived from \eqref{Tits} using substitution
\eqref{BTits}. For the first line of \eqref{BCentExt} we have:
 \be
  (S^B_1)^2\,=\,T^B_0\dot{s}^B_1T^B_0\dot{s}^B_1\,
  =\,(T^B_0)^2(\dot{s}^B_1)^2\,=\,e^{\pi\imath h_1}\,=\,e^{2\pi\imath(e_{\ell\ell}-e_{\ell+2,\ell+2})}\,
  =\,1\,,
 \ee
and since $\dot{s}^B_kT^B_k=T^B_{k-1}\dot{s}^B_k$ we derive
 \be
  (S^B_k)^2\,=\,T^B_k\dot{s}^B_kT^B_k\dot{s}^B_k\,
  =\,T^B_kT^B_{k-1}(\dot{s}^B_k)^2\,=\,(T^B_kT^B_{k-1})^2\,=\,1\,.
 \ee
For $i,j\in I$ such that $a_{ij}=0$, since
$\dot{s}^B_iT^B_j=T^B_j\dot{s}^B_i$ we have
 \be
  S^B_iS^B_j\,=\,T^B_i\dot{s}^B_iT^B_j\dot{s}^B_j\,
  =\,T^B_j\dot{s}^B_jT^B_i\dot{s}^B_i\,=\,S^B_jS^B_i\,,
 \ee
and similarly, since $\dot{s}^B_1T^B_2=T^B_2\dot{s}^B_1$, we obtain
 \be
  S^B_1S^B_2S^B_1S^B_2\,
  =\,T^B_0\dot{s}^B_1T^B_2\dot{s}^B_2T^B_0\dot{s}^B_1T^B_2\dot{s}^B_2\,
  =\,(T^B_0T^B_2)^2\dot{s}^B_1\dot{s}^B_2\dot{s}^B_1\dot{s}^B_2\\
  =\,\dot{s}^B_2\dot{s}^B_1\dot{s}^B_2\dot{s}^B_1\,
  =\,S^B_2S^B_1S^B_2S^B_1\,.
 \ee
The 3-move braid relation for $i,\,j=i+1$ we have:
 \be
  S^B_iS^B_{i+1}S^B_i\,
  =\,T^B_i\dot{s}^B_iT^B_{i+1}\dot{s}^B_{i+1}T^B_i\dot{s}^B_i\,
  =\,T^B_iT^B_{i+1}\dot{s}^B_iT^B_{i+1}\dot{s}^B_{i+1}\dot{s}^B_i\\
  =\,T^B_i(T^B_{i+1})^2\dot{s}^B_i\dot{s}^B_{i+1}\dot{s}^B_i\,
  =\,T^B_i(T^B_{i+1})^2\dot{s}^B_{i+1}\dot{s}^B_i\dot{s}^B_{i+1}\,
  =\,T^B_iT^B_{i+1}\dot{s}^B_{i+1}T^B_i\dot{s}^B_i\dot{s}^B_{i+1}\\
  =\,T^B_{i+1}\dot{s}^B_{i+1}T^B_iT^B_{i+1}\dot{s}^B_i\dot{s}^B_{i+1}\,
  =\,S^B_{i+1}S^B_iS^B_{i+1}\,,
 \ee
since $T_{i+1}\dot{s}^B_{i+1}\,=\,\dot{s}^B_{i+1}T^B_i$.\,\, $\Box$

\subsection{The case of type $C_{\ell}$ root system}\label{APPC}

The Lie algebra $\ssp_{2\ell}$ is identified with the
$\theta_C$-fixed subalgebra of $\gl_{2\ell}$:
 \be\label{StandardrepC}
  \frak{sp}_{2\ell}\,=\,\bigl\{X\in\gl_{2\ell}\,:\quad
  X\,=\,\theta_C(X)\bigr\}\,\subseteq\,\gl_{2\ell}\,,\qquad
  \theta_C(X)\,=\,-S_CX^{\tau}S_C^{-1}\,,\\
  S_C\,=\,\diag\bigl(1,\,-1,\,\ldots,\,1,\,-1\bigr)\,,
 \ee
where $\tau$ is the matrix transposition with respect to the
opposite diagonal \eqref{TranspOpp}. This provides the standard
faithful representation:
 \be\label{StandardrepC0}
  \phi\,:\quad\ssp_{2\ell}\,\longrightarrow\,\gl_{2\ell}\,,
 \ee
given by the following presentation of the Chevalley-Weyl generators
$h_i,\,e_i,\,f_i,\,i\in I$:
 \be\label{CFaithulRep}
  \phi(h_1)\,=\,e_{\ell\ell}\,-\,e_{\ell+1,\ell+1}\,,\\
  \phi(h_k)\,=\,(e_{\ell+1-k,\ell+1-k}+e_{\ell+k-1,\ell+k-1})\,
  -\,(e_{\ell+2-k,\ell+2-k}+e_{\ell+k,\ell+k})\,,\\1<k\leq\ell\,;\\
  \phi(e_1)\,=\,e_{\ell,\ell+1}\,,\quad
  \phi(e_k)\,=\,e_{\ell+1-k,\,\ell+2-k}+e_{\ell+k-1,\ell+k}\,,\quad1<k\leq\ell\,,\\
  \phi(f_1)\,=\,e_{\ell+1,\ell}\,,\quad
  \phi(f_k)\,=\,e_{\ell+2-k,\,\ell+1-k}+e_{\ell+k,\ell+k-1}\,,\quad1<k\leq\ell\,\,.
 \ee

The representation \eqref{CFaithulRep} yields the following
presentation of the type $C_{\ell}$ root system:
 \be\label{CAroots}
  \a_1^{C_{\ell}}\,
  =\,2\a_{\ell}^{A_{2\ell-1}}\,,\qquad
  \a_k^{C_{\ell}}
  =\,\a_{\ell+1-k}^{A_{2\ell-1}}\,+\,\a_{\ell-1+k}^{A_{2\ell-1}}\,,\qquad1<k\leq\ell\,\,.
 \ee
Consider the root system of type $A_{2\ell-1}$ endowed with the
automorphism $\iota$ of its Dynkin diagram:
 \be\label{AinvOddC}
  \iota\,:\quad i\,\longmapsto\,2\ell-i\,,\qquad i\in\{1,2,\ldots,2\ell-1\}\,;\\
  \a_i^{A_{2\ell}}\,\longmapsto\,\a_{2\ell-i}^{A_{2\ell-1}}\,,
 \ee
so that \eqref{CAroots} reads
 \be
  \a_1^{C_{\ell}}\,=\,\a_{\ell}^{A_{2\ell-1}}+\iota\bigl(\a_{\ell}^{A_{2\ell-1}}\bigr)\,,\quad
  \a_k^{C_{\ell}}\,=\,\a_{\ell+1-k}^{A_{2\ell-1}}\,+\,\iota\bigl
  (\a_{\ell+1-k}^{A_{2\ell-1}}\bigr)\,\,,\quad1<k\leq\ell\,.
 \ee
Introduce the Euclidean vector space $\IR^{2\ell}$ be with
orthonormal basis $\bigl\{\e_1,\ldots,\e_{2\ell}\bigr\}$, supplied
with an action of involution
 \be
  \iota\,:\quad\IR^{2\ell}\,\longrightarrow\,\IR^{2\ell}\,,\qquad
  \e_i\,\longmapsto\,-\e_{2\ell+1-i}\,,\quad1\leq i\leq2\ell\,.
 \ee
Consider the $\iota$-fixed Euclidean subspace
$\IR^{\ell}\subset\IR^{2\ell}$, spanned by
 \be\label{Cfolding}
 \e_i^{C_{\ell}}\,=\,\e_{\ell+1-i}\,+\,\iota(\e_{\ell+1-i})\,
  =\,\e_{\ell+1-i}\,-\,\e_{\ell+i}\,,\qquad1\leq i\leq\ell\,.
 \ee
Then the root data of type $C_{\ell}$  reads
 \be\label{RootDataC}
  \left\{\begin{array}{l}
  \a_1^{C_{\ell}}\,=\,2\e_1^{C_{\ell}}\\
  \a_k^{C_{\ell}}\,=\,\e_k^{C_{\ell}}-\e_{k-1}^{C_{\ell}}\,,\\
  1<k\leq\ell
  \end{array}\right.
 \hspace{1cm}
  \left\{\begin{array}{l}
  \vp_i^{C_{\ell}}\,=\,\e_i^{C_{\ell}}+\ldots+\e_{\ell}^{C_{\ell}}\,,\\
  i\in I\,;
  \end{array}\right.\\
  \<\a_i^{C_{\ell}},\,(\vp_j^{C_{\ell}})^{\vee}\>\,
  =\,\<(\a_i^{C_{\ell}})^{\vee},\,\vp_j^{C_{\ell}}\>\,=\,\delta_{ij}\,,\\
  \left\{\begin{array}{l}
  (\a_1^{C_{\ell}})^{\vee}\,=\,\e_1^{C_{\ell}}\\
  (\a_k^{C_{\ell}})^{\vee}\,=\,\e_k^{C_{\ell}}-\e_{k-1}^{C_{\ell}}\,,\\
  1<k\leq\ell\,;
  \end{array}\right.\qquad
  \left\{\begin{array}{l}
  (\vp_1^{C_{\ell}})^{\vee}\,=\,\frac{\e_1^{C_{\ell}}+\ldots+\e_{\ell}^{C_{\ell}}}{2}\\
  (\vp_k^{C_{\ell}})^{\vee}\,=\,\e_k^{C_{\ell}}+\ldots+\e_{\ell}^{C_{\ell}}\,,\\
  1<k\leq\ell\,.
  \end{array}\right.
 \ee
The Cartan matrix
$A\,=\,\|\<\a_j^{C_{\ell}},\,(\a_i^{C_{\ell}})^{\vee}\>\|$ and its
inverse are given by
 \be\label{CCartanMat}
  A\,
  =\,\left(\begin{smallmatrix}
  2&-1&0&\ldots&0\\-2&2&-1&\ddots&\vdots\\
  0&\ddots&\ddots&\ddots&0\\
  \vdots&\ddots&-1&2&-1\\
  0&\ldots&0&-1&2
  \end{smallmatrix}\right)\,,\qquad
  A^{-1}\,
  =\,\left(\begin{smallmatrix}
  \frac{\ell}{2}&\frac{\ell-1}{2}&\ldots&\frac{3}{2}&1&\frac{1}{2}\\
  \ell-1&\ell-1&\ldots&3&2&1\\
  \vdots&\vdots&\ddots&\vdots&\vdots\\
  3&3&\ldots&3&2&1\\
  2&2&\ldots&2&2&1\\
  1&1&\ldots&1&1&1
  \end{smallmatrix}\right)\,,\\
  h_i^{C_{\ell}}\,=\,(\a_i^{C_{\ell}})^{\vee}\,
  =\,\sum_{j\in I}\<\a_j^{C_{\ell}},\,(\a_i^{C_{\ell}})^{\vee}\>\vp_j^{\vee}\,
  =\,\sum_{j\in I}a_{ij}^{C_{\ell}}(\vp_j^{C_{\ell}})^{\vee}\,,\qquad
  \vp_i^{\vee}\,=\,\sum_{j\in I}c_{ij}^{C_{\ell}}(\a_j^{C_{\ell}})^{\vee}\,,\qquad i\in I\,.
 \ee

\begin{lem}\label{CAutroots}
The Weyl group $W(C_{\ell})$ is isomorphic to a subgroup of
$W(A_{2\ell-1})$ via
 \be\label{CWeyl}
  s_1^{C_{\ell}}\,=\,s_{\ell}^{A_{2\ell-1}}\,,\qquad
  s_k^{C_{\ell}}\,=\,s_{\ell+1-k}^{A_{2\ell-1}}s_{\ell-1+k}^{A_{2\ell-1}}\,
  =\,s_{\ell-1+k}^{A_{2\ell-1}}s_{\ell+1-k}^{A_{2\ell-1}}\,,\quad1<k\leq\ell\,.
 \ee
\end{lem}
\proof: For $k=1$ given
$\la=\la_1\e^{C_{\ell}}_1+\ldots+\la_{\ell}\e^{C_{\ell}}_{\ell}$ one
has
 \be
  s_1^{C_{\ell}}(\la)\,=\,\la\,
  -\,\<\la,\,(\a_1^{C_{\ell}})^{\vee}\>_{C_{\ell}}\a_1^{C_{\ell}}\,
  =\,\la\,
  -\,\<\la,\,\e_1^{C_{\ell}}\>_{C_{\ell}}2\e_1^{C_{\ell}}\\
  =\,-\la_1\e^{C_{\ell}}_1+\la_2\e^{C_{\ell}}_2+\ldots+\la_{\ell}\e^{C_{\ell}}_{\ell}\\
  =\,\la_{\ell}\e^{A_{2\ell-1}}_1+\ldots+\la_2\e^{A_{2\ell}}_{\ell-1}-\la_1
  \e^{A_{2\ell-1}}_{\ell}\,
  +\,\la_1\e^{A_{2\ell-1}}_{\ell+1}-\la_2\e^{A_{2\ell-1}}_{\ell+2}-\ldots-\la_{\ell}
  \e^{A_{2\ell-1}}_{2\ell}\\
  =\,s_{\ell}^{A_{2\ell-1}}(\la)\,,
 \ee
For $1<k\leq\ell$ the computation reproduces the one from Lemma
\ref{BAutroots}. $\Box$

The symplectic Lie group $\Sp_{2\ell}$ may be identified with the
$\theta_C$-fixed subgroup of $GL_{2\ell}$:
 \be\label{FaithulRepC}
  \Sp_{2\ell}\,=\,\bigl\{g\in GL_{2\ell}\,:\quad
  g\,=\,g^{\theta_C}\,\bigr\}\,\subset\,GL_{2\ell}\,\,.
 \ee
Let us describe the Tits group $W^T_{\Sp_{2\ell}}$, which is the
extension of the Weyl group $W(C_{\ell})$ \eqref{CWeyl}, and
identify it with the fixed point subgroup in $W^T_{GL_{2\ell}}$.

\begin{lem}\label{CTitsLift}
Let $\dot{s}^C_i=e^{\frac{\pi}{2}J_i},\,i\in I$ be the Tits
generators \eqref{Titsgen} then the following holds:
\begin{enumerate}
\item The Tits group $W^T_{\Sp_{2\ell}}$ is isomorphic to a
subgroup of $W^T_{GL_{2\ell}}$ via
 \be\label{TitsgenC}
  \dot{s}^C_1\,
  =\,\dot{s}^{A_{2\ell-1}}_{\ell}\,,\qquad
  \dot{s}^C_k\,
  =\,\dot{s}^{A_{2\ell-1}}_{\ell+1-k}\dot{s}^{A_{2\ell-1}}_{\ell-1+k}\,,\quad1<k\leq\ell\,.
 \ee

\item The elements \eqref{TitsgenC} belong to the $\theta_C$-fixed subgroup
$(W^T_{GL_{2\ell}})^{\theta_C}$.

\item Presentation \eqref{TitsgenC} matches with \eqref{CWeyl} due
to
 \be
  \Ad_{\dot{s}^C_i}|_{\fh}\,=\,s^C_i\,,\qquad i\in I\,.
 \ee
\end{enumerate}
\end{lem}

\proof: (1) Since $\dot{s}^C_i\in\Sp_{2\ell},\,i\in I$, one might
verify \eqref{TitsgenC} using the standard faithful representation
$\phi:\,\Sp_{2\ell}\to GL_{2\ell}$. For $i=1$, we have
 \be\label{TitsgenC1}
  \phi(\dot{s}^C_1)\,
  =\,\phi(e^{\frac{\pi}{2}J_1})\,
  =\,\left(\begin{smallmatrix}
  \Id_{\ell-1}&&&\\
  &0&-1&\\&1&0&\\&&&\Id_{\ell-1}
  \end{smallmatrix}\right)\,
  =\,\phi(\dot{s}_{\ell}^{A_{2\ell-1}}\,,
 \ee
and similarly, for $1<k\leq\ell$ we obtain
 \be\label{TitsgenC2}
  \phi(\dot{s}^C_k)\,
  =\,\phi(e^{\frac{\pi}{2}J_k})\,
  =\left(\begin{smallmatrix}
  \Id_{\ell-k}&&&&&&\\
  &0&-1&&&&\\&1&0&&&&\\
  &&&\Id_{2k-4}&&&\\&&&&0&-1&\\
  &&&&1&0&\\&&&&&&\Id_{\ell-k}
  \end{smallmatrix}\right)\,
  =\,\phi(\dot{s}_{\ell+1-k}^{A_{2\ell-1}}\dot{s}_{\ell-1+k}^{A_{2\ell-1}})\,.
 \ee

(2) From \eqref{CINV} one reads
 \be
  \theta_C\,:\qquad\dot{s}^{A_{2\ell-1}}_i\,\longmapsto\,
  \dot{s}^{A_{2\ell-1}}_{2\ell-i}\,,\qquad1\leq
  i\leq2\ell\,,
 \ee
so that using  \eqref{Tits} the expressions \eqref{TitsgenC} are
$\theta_C$-invariant.

(3) Follows by \eqref{Tits} and \eqref{CWeyl}. $\Box$

By analogy with the case of general linear group there is another
way to lift simple root generators $s^C_i\in W(C_{\ell})$ into the
Tits group $W^T_{\Sp_{2\ell}}$. Recall that $W^T_{GL_{2\ell}}$
contains the following elements:
 \be\label{TitsGL2L}
  T_k\,:=\,e^{\pi\imath e_{kk}}\,,\qquad1\leq k\leq (2\ell)\,;\\
  S_i\,=\,T_i\dot{s}_i\,=\,\dot{s}_iT_{i+1}\,\,,\qquad
  \overline{S}_i\,=\,T_{i+1}\dot{s}_i\,=\,\dot{s}_iT_i\,,\qquad 1\leq  i\leq2\ell\,.
 \ee
\begin{lem} The following elements belong to the fixed point
subgroup $(W^T_{GL_{2\ell}})^{\theta_C}$:
 \be\label{2torsC}
  S^C_1\,=\,T_{\ell}S_{\ell}\,=\,T_{\ell+1}\overline{S}_{\ell}\,\,,\qquad
  S^C_k\,=\,S_{\ell+1-k}\overline{S}_{\ell-1+k}\,,\quad1<k\leq\ell\,;\\
  T^C_i\,=\,T_{\ell+1-i}T_{\ell+i}\,=\,e^{\pi\imath\e^C_i}\,,\qquad1\leq i\leq\ell\,.
 \ee
\end{lem}
\proof: The explicit action of $\theta_C$ on generators
\eqref{TitsGL2L} reads \eqref{permatAction}:
 \be
  \theta_C(S_i)\,=\,\overline{S}_{2\ell-i}\,,\qquad
  \theta_C(\overline{S}_i)\,=\,S_{2\ell-i}\,,\qquad
  \theta_C(T_k)\,=\,T_{2\ell+1-k}\,.
 \ee
This implies that \eqref{2torsC} are invariant under $\theta_C$.\,\,
$\Box$

\begin{cor} The elements \eqref{TitsgenC} and \eqref{2torsC} can be identified via
 \be\label{CTits}
  S^C_1\,=\,\dot{s}_1^C\,\,,\qquad
  S^C_k\,=\,T^C_k\dot{s}^C_k\,\,,\quad1<k\leq\ell\,.
 \ee
\end{cor}
\proof: By straightforward computation using faithful representation
\eqref{FaithulRepC} one finds
 \be
  \phi(S^C_1)\,
  =\,\left(\begin{smallmatrix}
  \Id_{\ell-1}&&&\\
  &0&-1&\\&1&0&\\&&&\Id_{\ell-1}
  \end{smallmatrix}\right)\,
  =\,\phi(\dot{s}_{\ell}^{A_{2\ell-1}})\,,\\
  \phi(S^C_k)\,
  =\left(\begin{smallmatrix}
  \Id_{\ell-k}&&&&&&\\
  &0&1&&&&\\&1&0&&&&\\
  &&&\Id_{2k-4}&&&\\&&&&0&-1&\\
  &&&&-1&0&\\&&&&&&\Id_{\ell-k}
  \end{smallmatrix}\right)\,.
 \ee
Identifying this with \eqref{TitsgenC1}, \eqref{TitsgenC2} one
deduces \eqref{CTits}. $\Box$

\begin{prop}\label{CWeylLift}
The elements $S^C_i,\,i\in I$ of the Tits group $W^T_{\Sp_{2\ell}}$
satisfy the following relations:
 \be\label{CCentExt}
  (S^C_1)^2\,=\,T^C_1\,=\,e^{\pi\imath h^C_1}\,,\qquad
  (S^C_2)^2\,=\,\ldots\,
    =\,(S^C_{\ell})^2\,
  =\,1\,,\\
  S^C_iS^C_j\,=\,S^C_jS^C_i\,,\quad i,j\in I\quad \text{such
  that}\quad a_{ij}=0\,,\,;\\
  S^C_iS^C_jS^C_i\,=\,S^C_jS^C_iS^C_j\,,\qquad i,j\in I\quad\text{such
  that}\quad a_{ij}=-1\,,\\
  (S^C_1S^C_2)^4\,=\,(S^C_2S^C_1)^4\,=\,1\,,
 \ee
and they generate the whole Tits group $W^T_{\Sp_{2\ell}}$.
\end{prop}
\proof: One might verify \eqref{CCentExt} using the standard fathful
representation \eqref{FaithulRepC}. Alternatively, similarly to the
proof of Proposition \ref{BWeylLift} the relations \eqref{CCentExt}
can be deduced from \eqref{Tits} and \eqref{CTits}. $\Box$

\subsection{The case of type $D_{\ell}$ root system}\label{APPD}

The Lie algebra $\so_{2\ell}$ is identified with the a subalgebra of
$\gl_{2\ell}$ as follows:
 \be\label{StandardrepD}
  \frak{so}_{2\ell}\,=\,\bigl\{X\in\gl_{2\ell}\,:\quad
  X\,=\,\theta_D(X)\bigr\}\,\subseteq\,\gl_{2\ell}\,,\qquad
  \theta_D(X)\,=\,-S_DX^{\tau}S_D^{-1}\,,\\
  S_D\,=\,\diag\bigl(1,\,-1,\,\ldots,\,(-1)^{\ell-1};\,(-1)^{\ell-1},\,(-1)^{\ell},\,\ldots,\,1\bigr)\,,
 \ee
where $\tau$ is the matrix transposition with respect to the reverse
diagonal \eqref{TranspOpp}. This provides the standard faithful
representation:
 \be\label{StandardrepD0}
  \phi\,:\quad\so_{2\ell}\,\longrightarrow\,\gl_{2\ell}\,,
 \ee
which implies the following presentation of the Chevalley-Weyl
generators $h_i,\,e_i,\,f_i,\,i\in I$:
 \be\label{DFaithulRep}
  \phi(h_1)\,=\,(e_{\ell-1,\ell-1}+e_{\ell\ell})\,-\,(e_{\ell+1,\ell+1}+e_{\ell+2,\ell+2})\,,\\
  \phi(h_k)\,=\,(e_{\ell+1-k,\ell+1-k}+e_{\ell+k-1,\ell+k-1})\,
  -\,(e_{\ell+2-k,\ell+2-k}+e_{\ell+k,\ell+k})\,,\\1<k\leq\ell\,,\\
  \phi(e_1)\,=\,e_{\ell-1,\ell+1}+e_{\ell,\ell+2}\,,\\
  \phi(e_k)\,=\,e_{\ell+1-k,\,\ell+2-k}+e_{\ell+k-1,\ell+k}\,,\quad1<k\leq\ell\,,\\
 \hspace{-1cm}
  \phi(f_1)\,=\,e_{\ell+1,\ell-1}+e_{\ell+2,\ell}\,,\\
  \phi(f_k)\,=\,e_{\ell+2-k,\,\ell+1-k}+e_{\ell+k,\ell+k-1}\,,\quad1<k\leq\ell\,\,.
 \ee
The representation \eqref{DFaithulRep} yields the following
presentation of the type $D_{\ell}$ root system:
 \be\label{DAroots}
  \a_1^{D_{\ell}}\,
  =\,\a_{\ell-1}^{A_{2\ell-1}}\,+\,2\a_{\ell}^{A_{2\ell-1}}\,+\,
  \a_{\ell+1}^{A_{2\ell-1}}\,,\\
  \a_k^{D_{\ell}}\,
  =\,\a_{\ell+1-k}^{A_{2\ell-1}}\,+\,\a_{\ell-1+k}^{A_{2\ell-1}}\,,\qquad1<k\leq\ell\,\,.
 \ee
Consider the root system of type $A_{2\ell-1}$ endowed with the
automorphism $\iota$ of its Dynkin diagram:
 \be\label{AinvOdd}
  \iota\,:\quad i\,\longmapsto\,2\ell-i\,,\qquad i\in\{1,2,\ldots,2\ell-1\}\,;\\
  \a_i^{A_{2\ell}}\,\longmapsto\,\a_{2\ell-i}^{A_{2\ell-1}}\,.
 \ee
so that \eqref{DAroots} reads
 \be
  \a_1^{D_{\ell}}\,=\,\a_{\ell-1}^{A_{2\ell-1}}+\a_{\ell}^{A_{2\ell-1}}\,
  +\,\iota\bigl(\a_{\ell-1}^{A_{2\ell-1}}+\a_{\ell}^{A_{2\ell-1}}\bigr)\,,\\
  \a_k^{D_{\ell}}\,=\,\a_{\ell+1-k}^{A_{2\ell-1}}\,+\,
  \iota\bigl(\a_{\ell+1-k}^{A_{2\ell-1}}\bigr)\,,\qquad1<k\leq\ell\,.
 \ee
Introduce Euclidean vector space $\IR^{2\ell}$ with orthonormal
basis $\bigl\{\e_1,\ldots,\e_{2\ell}\bigr\}$, supplied with an
action of involution:
 \be
  \iota\,:\quad\IR^{2\ell}\,\longrightarrow\,\IR^{2\ell}\,,\qquad
  \e_i\,\longmapsto\,-\e_{2\ell+1-i}\,,\quad1\leq i\leq(2\ell)\,.
 \ee
Consider the $\iota$-fixed Euclidean subspace
$\IR^{\ell}\subset\IR^{2\ell}$, spanned by
 \be\label{Dfolding}
  \e_i^{D_{\ell}}\,=\,\e_{\ell+1-i}^{A_{2\ell-1}}\,+\,\iota
  \bigl(\e_{\ell+1-i}^{A_{2\ell-1}}\bigr)\,
  =\,\e_{\ell+1-i}^{A_{2\ell-1}}\,-\,\e_{\ell+i}^{A_{2\ell-1}}\,,\qquad
  1\leq i\leq\ell\,.
 \ee
Then the root data of type $D_{\ell}$ reads
 \be\label{RootDataD}
  \left\{\begin{array}{l}
  \a^{D_{\ell}}_1\,=\,\e_2^{D_{\ell}}+\e_1^{D_{\ell}}\\
  \\
  \a^{D_{\ell}}_k\,=\,\e_k^{D_{\ell}}-\e_{k-1}^{D_{\ell}}\,,\\
  1<k\leq\ell
  \end{array}\right.
 \hspace{1cm}
  \left\{\begin{array}{l}
  \vp^{D_{\ell}}_1\,=\,\frac{\e_1^{D_{\ell}}+\e_2^{D_{\ell}}+\ldots+\e_{\ell}^{D_{\ell}}}{2}\\
  \vp^{D_{\ell}}_2\,=\,\frac{-\e_1^{D_{\ell}}+\e_2^{D_{\ell}}+\ldots+\e_{\ell}^{D_{\ell}}}{2}\\
  \vp^{D_{\ell}}_k\,=\,\e_k^{D_{\ell}}+\ldots+\e_{\ell}^{D_{\ell}}\,,\quad2<k\leq\ell\,;
  \end{array}\right.\\
  \<\a_i^{D_{\ell}},\,(\vp_j^{D_{\ell}})^{\vee}\>\,
  =\,\<(\a_i^{D_{\ell}})^{\vee},\,\vp_j^{D_{\ell}}\>\,=\,\delta_{ij}\,,
 \ee
and the Cartan matrix
$A\,=\,\|\<\a^{D_{\ell}}_j,\,(\a_i^{D_{\ell}})^{\vee}\>\|$ and its
inverse are given by
 \be\label{DCartanMat}
  A\,
  =\,\left(\begin{smallmatrix}
  2&0&-1&0&\ldots&0\\
  0&2&-1&0&\ddots&\vdots\\
  -1&-1&2&-1&\ddots&\vdots\\
  0&\ddots&\ddots&\ddots&\ddots&0\\
  \vdots&\ddots&0&-1&2&-1\\
  0&\ldots&\ldots&0&-1&2
  \end{smallmatrix}\right)\,,\quad
  A^{-1}\,=\,\|c_{ij}^{D_{\ell}}\|\,
  =\,\left(\begin{smallmatrix}
  \frac{\ell}{4}&\frac{\ell-2}{4}&\frac{\ell-2}{2}&\ldots&\frac{3}{2}&1&\frac{1}{2}\\
  \frac{\ell-2}{4}&\frac{\ell}{4}&\frac{\ell-2}{2}&\ldots&\frac{3}{2}&1&\frac{1}{2}\\
  \frac{\ell-2}{2}&\frac{\ell-2}{2}&\ell-2&\ldots&3&2&1\\
  \vdots&\vdots&\vdots&\ddots&\vdots&\vdots&\vdots\\
  \frac{3}{2}&\frac{3}{2}&3&\ldots&3&2&1\\
  1&1&2&\ldots&2&2&1\\
  \frac{1}{2}&\frac{1}{2}&1&\ldots&1&1&1
  \end{smallmatrix}\right)\,,\\
  h_i^{D_{\ell}}\,=\,(\a_i^{D_{\ell}})^{\vee}\,
  =\,\sum_{j\in I}\<\a_j^{D_{\ell}},\,(\a_i^{D_{\ell}})^{\vee}\>\vp_j^{\vee}\,
  =\,\sum_{j\in I}a_{ij}^{D_{\ell}}(\vp_j^{D_{\ell}})^{\vee}\,,\\
  \vp_i^{\vee}\,=\,\sum_{j\in I}c_{ij}^{D_{\ell}}(\a_j^{D_{\ell}})^{\vee}\,,\qquad i\in I\,.
 \ee

\begin{lem}\label{DAutroots}
The Weyl group $W(D_{\ell})$ is isomorphic to a subgroup of
$W(A_{2\ell-1})$ via
 \be\label{DWeyl}
 s_1^{D_{\ell}}\,=\,s_{\ell}^{A_{2\ell-1}}
 s_{\ell-1}^{A_{2\ell-1}}s_{\ell+1}^{A_{2\ell-1}}s_{\ell}^{A_{2\ell-1}}\,
  =\,s_{\ell}^{A_{2\ell-1}}s_{\ell+1}^{A_{2\ell-1}}
  s_{\ell-1}^{A_{2\ell-1}}s_{\ell}^{A_{2\ell-1}}\,,\\
  s_k^{D_{\ell}}\,=\,s_{\ell+1-k}^{A_{2\ell-1}}s_{\ell-1+k}^{A_{2\ell-1}}\,
  =\,s_{\ell-1+k}^{A_{2\ell-1}}s_{\ell+1-k}^{A_{2\ell-1}}\,,\quad1<k\leq\ell\,.
 \ee
In particular, one has
$s_1^{D_{\ell}}\,=\,s_{\ell}^{A_{2\ell-1}}s_2^{D_{\ell}}s_{\ell}^{A_{2\ell-1}}$,
so that
$\Out(\Phi(D_{\ell}))=\<s_{\ell}^{A_{2\ell-1}}\>\simeq\IZ/2\IZ$.
\end{lem}
\proof: For $1<k\leq\ell$ the proof literarily follows similar
statement of Lemma 5.1. For $k=1$  given
$\la=\la_1\e^{D_{\ell}}_1+\ldots+\la_{\ell}\e^{D_{\ell}}_{\ell}$ one
has
 \be
  s_1^{D_{\ell}}(\la)\,=\,\la\,
  -\,\<\la,\,(\a_1^{D_{\ell}})^{\vee}\>_{D_{\ell}}\a_1^{D_{\ell}}\,
  =\,\la\,
  -\,\<\la,\,\e_1^{D_{\ell}}+\e_2^{D_{\ell}}\>_{D_{\ell}}
  (\e_1^{D_{\ell}}+\e_2^{D_{\ell}})\\
  =\,-\la_2\e^{D_{\ell}}_1-\la_1\e^{D_{\ell}}_2+\la_3
  \e^{D_{\ell}}_3+\ldots+\la_{\ell}\e^{D_{\ell}}_{\ell}\\
  =\,\la_{\ell}\e^{A_{2\ell-1}}_1+\ldots+\la_3
  \e^{A_{2\ell}}_{\ell-2}-\la_1\e^{A_{2\ell-1}}_{\ell-1}-\la_2\e^{A_{2\ell-1}}_{\ell}\\
  +\,\la_2\e^{A_{2\ell-1}}_{\ell+1}+\la_1\e^{A_{2\ell-1}}_{\ell+2}-
  \la_3\e^{A_{2\ell-1}}_{\ell+3}-\ldots-\la_{\ell}\e^{A_{2\ell-1}}_{2\ell}\\
  =\,s_{\ell}^{A_{2\ell-1}}s_{\ell-1}^{A_{2\ell-1}}
  s_{\ell+1}^{A_{2\ell-1}}s_{\ell}^{A_{2\ell-1}}(\la)\,,
 \ee
since
$s_{\ell}^{A_{2\ell-1}}s_{\ell-1}^{A_{2\ell-1}}s_{\ell+1}^{A_{2\ell-1}}s_{\ell}^{A_{2\ell-1}}$
simply swaps
$\e^{A_{2\ell-1}}_{\ell}\leftrightarrow\e^{A_{2\ell-1}}_{\ell+2}$
and
$\e^{A_{2\ell-1}}_{\ell-1}\leftrightarrow\e^{A_{2\ell-1}}_{\ell+1}$.
$\Box$

The orthogonal group $O_{2\ell}$ may be identified with the fixed
point subgroup of general linear group:
 \be\label{FaithulRepD}
  O_{2\ell}\,=\,\bigl\{g\in GL_{2\ell}\,:\quad
  g\,=\,g^{\theta_D}\bigr\}\,\subset\,GL_{2\ell}\,.
 \ee
The group $O_{2\ell}$ is not simply connected and there is the
following exact sequence:
 \be\label{SOeven}
  \xymatrix{1\ar[r]&SO_{2\ell}\ar[r]&O_{2\ell}\ar[r]^{\det}&\{\pm1\}\ar@/^1pc/[l]^z\ar[r]&1}\,,
 \ee
with the non-trivial action of $\{\pm1\}$ on $SO_{2\ell}$.

\begin{lem} The even orthogonal group allows for the
following decomposition:
 \be
  O_{2\ell}\,=\,SO_{2\ell}\sqcup S_{\ell}\cdot SO_{2\ell}\,,
  \qquad
  S_{\ell}\,
  =\,\left(\begin{smallmatrix}
  \Id_{\ell-1}&&&\\0&0&1&\\
  0&1&0&\\&&&\Id_{\ell-1}
  \end{smallmatrix}\right)\,,\\
  \Ad_{S_{\ell}}\,:\quad\a_1^{D_{\ell}}\,\longleftrightarrow\,\a_2^{D_{\ell}}\,\,.
 \ee
so that $\Ad_{S_{\ell}}$ represents the symmetry of Dynkin diagram
of type $D_{\ell}$. Moreover, $O_{2\ell}$ has a structure of a
semidirect product:
$O_{2\ell}\,=\,SO_{2\ell}\rtimes\Out(SO_{2\ell})$.
\end{lem}
\proof: we have
 \be
  (H_{GL_{2\ell}})^{\theta_D}\,=\,H_{O_{2\ell}}\,=\,H_{\Sp_{2\ell}}\,
  =\,(H_{GL_{2\ell}})^{\theta_C}\,,
 \ee
hence $H_{O_{2\ell}}=H_{\SO_{2\ell}}$. It is convenient to pick
$z(-1)=M_{\ell}\in O_{2\ell}$ as a representative of section
$z:\{\pm1\}\to O_{2\ell}$, since its represents the outer
automorphism of $SO_{2\ell}$. $\Box$

\begin{lem}\label{DTitsLift}
Let $\dot{s}^D_i=e^{\frac{\pi}{2}J_i},\,i\in I$ be the Tits
generators \eqref{Titsgen} then the following holds:
\begin{enumerate}
\item The Tits group $W^T_{\SO_{2\ell}}$ is isomorphic to a
subgroup of $W^T_{GL_{2\ell}}$ via
 \be\label{TitsgenD}
  \dot{s}^D_1\,
  =\,\dot{s}_{\ell-1}^{A_{2\ell-1}}\dot{s}_{\ell}^{A_{2\ell-1}}
  (\dot{s}^{A_{2\ell-1}}_{\ell+1})^{-1}(\dot{s}^{A_{2\ell-1}}_{\ell-1})^{-1}\,,\\
  \dot{s}^D_k\,
  =\,\dot{s}^{A_{2\ell-1}}_{\ell+1-k}\dot{s}^{A_{2\ell-1}}_{\ell-1+k}\,,\quad1<k\leq\ell\,.
 \ee

\item The elements \eqref{TitsgenD} belong to the $\theta_D$-fixed subgroup
$(W^T_{GL_{2\ell}})^{\theta_D}$.

\item Presentation \eqref{TitsgenD} matches with \eqref{DWeyl} due
to
 \be
  \Ad_{\dot{s}^D_i}|_{\fh}\,=\,s^D_i\,,\qquad i\in I\,.
 \ee

\item The Tits group $W^T_{O_{2\ell}}$ fits into the following
exact sequence:
 \be\label{DNorm}\xymatrix{
  1\ar[r]&W^T_{SO_{2\ell}}\ar[r]&
  W^T_{O_{2\ell}}\ar[r] &
  \Out(\Phi(D_{\ell}))\ar[r]&1}\,,
 \ee
that actually splits, so that
$W^T_{O_{2\ell}}=W^T_{SO_{2\ell}}\rtimes\Out(\Phi(D_{\ell}))$.
\end{enumerate}
\end{lem}

\proof: (1) Since $\dot{s}^D_i\in SO_{2\ell},\,i\in I$, one might
verify \eqref{TitsgenD} using the standard faithful representation
$\phi:\,SO_{2\ell}\to GL_{2\ell}$. For $i=1$, we have
 \be\label{TitsgenD1}
  \phi(\dot{s}^D_1)\,
  =\,\phi(e^{\frac{\pi}{2}J_1})\,
  =\,\left(\begin{smallmatrix}
  \Id_{\ell-2}&&&&&&\\
  &0&0&-1&0&\\
  &0&0&0&-1&\\
  &1&0&0&0&\\
  &0&1&0&0&\\
  &&&&&\Id_{\ell-2}
  \end{smallmatrix}\right)\,
  =\,\phi(\dot{s}_{\ell}^{A_{2\ell-1}})\,,
 \ee
and similarly, for $1<k\leq\ell$ we obtain
 \be\label{TitsgenD2}
  \phi(\dot{s}^D_k)\,
  =\,\phi(e^{\frac{\pi}{2}J_k})\,
  =\left(\begin{smallmatrix}
  \Id_{\ell-k}&&&&&&\\
  &0&-1&&&&\\&1&0&&&&\\
  &&&\Id_{2k-4}&&&\\&&&&0&-1&\\
  &&&&1&0&\\&&&&&&\Id_{\ell-k}
  \end{smallmatrix}\right)\,
  =\,\phi(\dot{s}_{\ell+1-k}^{A_{2\ell-1}}\dot{s}_{\ell-1+k}^{A_{2\ell-1}})\,.
 \ee

(2) From \eqref{CINV} one reads
 \be
  \theta_D\,:\qquad\dot{s}^{A_{2\ell-1}}_i\,\longmapsto\,
  \dot{s}^{A_{2\ell-1}}_{2\ell-i}\,,\qquad1\leq
  i\leq2\ell\,,
 \ee
so that \eqref{TitsgenD} are $\theta_D$-invariant.

(3) Follows by \eqref{Tits} and \eqref{DWeyl}.

(4) The element $\dot{R}=\dot{s}_{\ell}$ represents
$R\in\Out(D_{\ell})$ from \eqref{DWeyl}; it is not a reflection at
any root of $\Phi(D_{\ell})$, that is
$\Ad_{\dot{s}_{\ell}}|_{\fh}=s_{\ell}\notin W(D_{\ell})$ due to
\eqref{DWeyl}. Moreover, by \eqref{OUT}
$\Out(\Phi(D_{\ell}))=\Out(SO_{2\ell})$ one obtains \eqref{DNorm}.
$\Box$

By analogy with the case of general linear group there is another
way to lift simple root generators $s^D_i\in W(D_{\ell})$ into the
Tits group $W^T_{\SO_{2\ell}}$. Recall that $W^T_{GL_{2\ell}}$
contains the following elements:
 \be\label{TitsGL2LL}
  T_k\,:=\,e^{\pi\imath e_{kk}}\,,\qquad1\leq k\leq 2\ell\,;\\
  S_i\,=\,T_i\dot{s}_i\,=\,\dot{s}_iT_{i+1}\,,\qquad
  \overline{S}_i\,=\,T_{i+1}\dot{s}_i\,=\,\dot{s}_iT_i\,,\qquad 1\leq  i\leq2\ell\,.
 \ee
\begin{lem} The following elements belong to the fixed point
subgroup
$(W^T_{GL_{2\ell}})^{\theta_D}$:
 \be\label{2torsD}
  S^D_1\,=\,S_{\ell}S_{\ell-1}\overline{S}_{\ell+1}S_{\ell}\,,\\
  S^D_k\,=\,S_{\ell+1-k}\overline{S}_{\ell-1+k}\,,\qquad1<k\leq\ell\,;\\
  T^D_i\,=\,T_{\ell+1-i}T_{\ell+i}\,,\qquad1\leq i\leq\ell\,,
 \ee
\end{lem}
\proof: The explicit action of $\theta_C$ on generators
\eqref{TitsGL2LL} reads \eqref{permatAction}:
 \be
  \theta_D(S_i)\,=\,\overline{S}_{2\ell-i}\,,\qquad
  \theta_D(\overline{S}_i)\,=\,S_{2\ell-i}\,,\qquad\theta_D(T_k)\,=\,T_{2\ell+1-k}\,.
 \ee
This implies that \eqref{2torsD} are invariant under $\theta_D$.\,\,
$\Box$

\begin{cor} The elements \eqref{TitsgenD} and \eqref{2torsD} can be identified via
 \be\label{DTits}
  S^D_1\,=\,T^D_2\dot{s}^D_1\,,\qquad
  S^D_k\,=\,T^D_k\dot{s}^D_k\,\,,\quad1<k\leq\ell\,.
 \ee
\end{cor}
\proof: By straightforward computation using faithful representation
\eqref{FaithulRepD} one finds
 \be
  \phi(S^D_1)\,
  =\,\left(\begin{smallmatrix}
  \Id_{\ell-2}&&&&&&\\
  &0&0&1&0&\\
  &0&0&0&-1&\\
  &1&0&0&0&\\
  &0&-1&0&0&\\
  &&&&&\Id_{\ell-2}
  \end{smallmatrix}\right)\,,\\
  \phi(S^D_k)\,
  =\left(\begin{smallmatrix}
  \Id_{\ell-k}&&&&&&\\
  &0&1&&&&\\&1&0&&&&\\
  &&&\Id_{2k-4}&&&\\&&&&0&-1&\\
  &&&&-1&0&\\&&&&&&\Id_{\ell-k}
  \end{smallmatrix}\right)\,.
 \ee
Identifying this with \eqref{TitsgenC1}, \eqref{TitsgenC2} one
deduces \eqref{CTits}. $\Box$

\begin{prop}\label{DWeylLift}
The elements $S^D_i,\,i\in I$ from \eqref{2torsD} satisfy the
relations \eqref{BR0}, \eqref{BR} of type $D_{\ell}$:
 \be\label{DCentExt}
  (S^D_1)^2\,=\,(S^D_2)^2\,=\,\ldots\,=\,(S^D_{\ell})^2\,
  =\,1\,,\\
  S^D_iS^D_j\,=\,S^D_jS^D_i\,,\quad i,j\in I\quad \text{such
  that}\quad a_{ij}=0\,,\,;\\
  S^D_iS^D_jS^D_i\,=\,S^D_jS^D_iS^D_j\,,\qquad i,j\in I\quad\text{such
  that}\quad a_{ij}=-1\,.
 \ee
\end{prop}
\proof: One might verify \eqref{DCentExt} using the standard
faithful representation \eqref{FaithulRepD}. Alternatively,
similarly to our proof of Proposition \ref{BWeylLift} the relations
\eqref{DCentExt} can be deduced from \eqref{Tits} and \eqref{DTits}.
$\Box$

\noindent {\small {\bf A.A.G.} {\sl Laboratory for Quantum Field Theory
and Information},\\
\hphantom{xxxx} {\sl Institute for Information
Transmission Problems, RAS, 127994, Moscow, Russia};\\
\hphantom{xxxx} {\it E-mail address}: {\tt anton.a.gerasimov@gmail.com}}\\
\noindent{\small {\bf D.R.L.}
{\sl Laboratory for Quantum Field Theory
and Information},\\
\hphantom{xxxx}  {\sl Institute for Information
Transmission Problems, RAS, 127994, Moscow, Russia};\\
\hphantom{xxxx} {\sl Moscow Center for Continuous Mathematical
Education,\\
\hphantom{xxxx} 119002,  Bol. Vlasyevsky per. 11, Moscow, Russia};\\
\hphantom{xxxx} {\it E-mail address}: {\tt lebedev.dm@gmail.com}}\\
\noindent{\small {\bf S.V.O.} {\sl
 School of Mathematical Sciences, University of Nottingham\,,\\
\hphantom{xxxx} University Park, NG7\, 2RD, Nottingham, United Kingdom};\\
\hphantom{xxxx} {\sl
 Institute for Theoretical and Experimental Physics,
117259, Moscow, Russia};\\
\hphantom{xxxx} {\it E-mail address}: {\tt oblezin@gmail.com}}

\end{document}